\documentclass[11pt]{amsart}

\usepackage{amsmath, amssymb, amsthm, amsfonts,bm,mathtools}
\usepackage{booktabs}
\usepackage{array}
\usepackage{caption}

\usepackage[margin=1in]{geometry}

\usepackage[colorlinks=true,linkcolor=blue,
urlcolor=blue,citecolor=blue]{hyperref}

\usepackage{aliascnt}
\usepackage[nameinlink]{cleveref}

\newtheorem{theorem}{Theorem}[section]

\newaliascnt{lemma}{theorem}
\newtheorem{lemma}[lemma]{Lemma}
\aliascntresetthe{lemma}

\newaliascnt{corollary}{theorem}
\newtheorem{corollary}[corollary]{Corollary}
\aliascntresetthe{corollary}

\newaliascnt{proposition}{theorem}
\newtheorem{proposition}[proposition]{Proposition}
\aliascntresetthe{proposition}

\crefname{theorem}{theorem}{theorems}
\Crefname{theorem}{Theorem}{Theorems}

\crefname{lemma}{lemma}{lemmas}
\Crefname{lemma}{Lemma}{Lemmas}

\crefname{corollary}{corollary}{corollaries}
\Crefname{corollary}{Corollary}{Corollaries}

\crefname{proposition}{proposition}{propositions}
\Crefname{proposition}{Proposition}{Propositions}

\newaliascnt{conjecture}{theorem}
\newtheorem{conjecture}[conjecture]{Conjecture}
\aliascntresetthe{conjecture}

\crefname{conjecture}{conjecture}{conjectures}
\Crefname{conjecture}{Conjecture}{Conjectures}

\theoremstyle{definition}
\newtheorem{remark}[theorem]{Remark}

\numberwithin{equation}{section}
\title{Asymptotics for generalized integer partitions}
\author{Armanzhan Onggarbay \and Damir Yeliussizov}
\address{Kazakh-British Technical University, Almaty, Kazakhstan}
\email{\href{mailto:arrmanzhan@gmail.com}{arrmanzhan@gmail.com}}
\email{\href{mailto:yeldamir@gmail.com}{yeldamir@gmail.com}}
\date{}

\begin{document}

\begin{abstract}
We study the asymptotic behavior of weighted partition functions $p_f(n)$, defined by the generating function $\prod_{n \geq 1} (1 - q^{n})^{-f(n)}$ with non-negative real weights $f(n)$. We assume that the associated Dirichlet series $D(s) = \sum_{n \geq 1} f(n) n^{-s}$ admits a meromorphic continuation to a half-plane $\Re(s) \geq -\varepsilon$, with finitely many real poles of arbitrary orders, the largest of which is a positive pole $\alpha$ of order $k$, and a possibly infinite collection of complex poles lying strictly to the left of $\alpha$, with some additional conditions. 
We extend the saddle-point method of Debruyne and Tenenbaum, and show that the pole $\alpha$ dominates the contribution of all other poles. 
As a result, we give an explicit asymptotic formula for $p_f(n)$ described via generalized Lambert functions. 
Furthermore, we 
  show that for $k \geq 2$ no finite combination of powers of $n$, $\log n$, $\log \log n, \ldots$ approximates the saddle point well enough to give an asymptotic formula for $p_f(n)$.
We also discuss a number of applications related to algebra, number theory, and combinatorics. 
\end{abstract}

\maketitle

\setcounter{tocdepth}{1}
\tableofcontents

\section{Introduction}
The study of the asymptotic behaviour of integer partition functions goes back to Hardy and Ramanujan \cite{HardyRamanujan}, who proved that the number $p(n)$ of partitions of $n$ satisfies
\[
p(n) \sim \frac{1}{4\sqrt{3}\, n} e^{\pi \sqrt{2n/3}}, \quad n \to \infty.
\]
Meinardus \cite{Meinardus} obtained asymptotic formulas for coefficients of more general products of the form 
$$\prod_{n \geq 1} (1 - q^{n})^{-f(n)} = \sum_{n = 0}^{\infty} p_f(n) q^n,$$ 
assuming that the Dirichlet series $D(s) := \sum_{n \geq 1} f(n) n^{-s}$ has a single simple pole at its abscissa of convergence, and imposing an analytic condition on the generating function away from the real axis. 
Granovsky, Stark and Erlihson \cite{GranovskyStarkErlihson} gave a probabilistic proof of Meinardus' theorem and extended it to other decomposable structures, and Granovsky and Stark \cite{GranovskyStark2012, GranovskyStark2015} extended it, by the Khintchine--Meinardus probabilistic method, to Dirichlet series with several simple poles on the positive real axis and possibly a simple pole at $0$. 

More recently, Debruyne and Tenenbaum \cite{Tenenbaum2020} developed a saddle-point method in which this analytic condition is replaced by an arithmetic condition on the set $\Lambda$ of parts. In their setting, the weights are the indicator function of $\Lambda$, and the Dirichlet series has a single simple pole in a half-plane $\Re(s) \geq -\varepsilon$. Bridges, Brindle, Bringmann and Franke \cite{Bridges2023} extended the saddle-point approach to Dirichlet series with several simple real poles, and Bringmann, Franke and Heim \cite{BringmannFrankeHeim} used their result to determine the asymptotics of the number of commuting $\ell$-tuples of permutations. Ahmadi, G\'omez-A\'{\i}za and Ward \cite{AhmadiGomezAizaWard} studied the families of products whose Dirichlet series are $\zeta(s - 1)^{i} \zeta(s + 1)^{j} \zeta(s)^{k}$, which in general have poles of higher order; they obtained asymptotic formulas in a few cases, some of them in terms of the Lambert function, and the asymptotics of the logarithm of the coefficients for the general families. For some specific weights, asymptotic formulas were also obtained by the circle method, for instance for multiplicities given by divisor sums \cite{BerndtRoblesZaharescuZeindler}, for weights built from the M\"obius function \cite{BasakRoblesZaharescu}, and for the number of $n$-dimensional representations of $SU(3)$ \cite{Romik, BringmannFranke2023}.

\subsection{Synopsis of main results}
In this paper, we extend the Debruyne--Tenenbaum method to arbitrary non-negative real weights $f(n)$, whose associated Dirichlet series may have finitely many real poles $\alpha_1 > \cdots > \alpha_m$ of arbitrary orders $k_1, \ldots, k_m$, where $\alpha = \alpha_1 > 0$ and the other real poles may be positive, zero or negative, together with a possibly infinite collection of complex poles lying to the left of $\alpha$. One of our main results, \Cref{thm:Explicit saddle point}, is an explicit asymptotic formula for $p_f(n)$ when the other poles lie to the left of a line $\Re(s) = \sigma_0$ with $\sigma_0 < \alpha/2$, on which $D(s)$ satisfies the growth condition of (A2) below. The pole $\alpha$ and its order $k = k_1$ determine the saddle point $\rho_\alpha$ of the residue at $z = \alpha$, which is an explicit power of $n$ for $k = 1$, is given by the Lambert function for $k = 2$, and by a generalized Lambert function for $k \geq 3$, and the other poles enter through their residues. When $k \geq 2$, the order of the pole produces powers of $\log \rho_\alpha^{-1}$ in the formula. The cases $\sigma_0 \geq \alpha/2$ are covered in fruther results as well: for every $\sigma_0 < \alpha$, a fixed number of Newton steps starting from $\rho_\alpha$ gives an asymptotic formula (\Cref{thm:Dominant saddle point}(iii)), and so does an expansion of the exponent around $\rho_\alpha$ to a fixed order (\Cref{thm:Exponent expansion}).

All these results rest on the key saddle-point formula $p_f(n) \sim e^{\rho n} P(\rho)/\sqrt{2 \pi \Psi''(\rho)}$ (\Cref{thm:Main theorem}), in which $\Psi(\rho) = \log P(\rho)$ is the sum of the residues of all poles up to an error $O(\rho^{\varepsilon})$, and $\rho$ is the exact saddle point. This formula has the classical form, and its proof extends the method of Debruyne and Tenenbaum. We state these general results in \Cref{sec:Further results}, where we also 
show that for $k \geq 2$ the saddle point {cannot} be replaced by finitely many terms of its asymptotic expansion (\Cref{thm:Accuracy of saddle point}), 
and as a result  no finite combination of powers of $n$, $\log n$, $\log \log n, \ldots$ approximates the saddle point well enough to give an asymptotic formula for $p_f(n)$.
Furthermore, we conjecture that the generalized Lambert function cannot even be approximated by the class of logarithmico-exponential functions of Hardy up to an error of a negative power of $n$ (\Cref{thm:No closed form}). 

\medskip

First, we introduce the weighted partition functions under consideration.

\subsection{Weighted partitions and main conditions}
Let $f: \mathbb{N} \to \mathbb{R}_{\geq 0}$ be an arithmetic function, and define the generating function 
\begin{equation}
\label{eq:Generating function}
P(\tau) := \prod_{n = 1}^{\infty} (1 - q^{n})^{-f(n)} = \sum_{n = 0}^{\infty} p_f(n) q^n, \quad q = e^{-\tau}, \ \Re(\tau) > 0,
\end{equation}
where each factor $(1 - q^{n})^{-f(n)}$ is expanded by the generalized binomial theorem. The associated Dirichlet series is 
\begin{equation}
D(s) := \sum_{n = 1}^{\infty} \frac{f(n)}{n^s}.
\end{equation}
We consider Dirichlet series satisfying the following conditions: 
\begin{enumerate}
    \item[(A1)] (Poles) The function $D(s)$ admits a meromorphic continuation to the half-plane
    \[
    \Re(s) \geq -\varepsilon
    \]
    for some $0 < \varepsilon < 1$. In this region, $D(s)$ has finitely many distinct real poles
    \[
    \mathcal{A} := \{\alpha_1, \ldots, \alpha_m\}, \quad \alpha_1 > \cdots > \alpha_m, \quad \alpha_1 > 0,
    \]
    where the pole at $s = \alpha_i$ has order (i.e. multiplicity) $k_i \in \mathbb{N}$; we also set $\alpha = \alpha_1$ and $k = k_1$. The poles $\alpha_2, \ldots, \alpha_m$ may be positive, zero or negative. Let $h_0$ be the order of the pole of $D(s)$ at $s = 0$. Together with $\mathcal{A}$, the series $D(s)$ has the collection of complex poles 
    \[
    \mathcal{B} := \{\beta_1, \beta_2, \ldots\},
    \]
    where the pole $\beta \in \mathcal{B}$ has order $\gamma_\beta \in \mathbb{N}$. We assume that the complex poles lie strictly to the left of $\alpha$:
    \begin{equation}
    \label{eq:Complex poles to the left}
    \sup_{\beta \in \mathcal{B}} \Re(\beta) < \alpha.
\end{equation}
    Moreover, in a neighborhood of $s = p$, where $p$ is a pole of order $h$ of $D(s)$, we have the expansion 
    \begin{equation}
    \label{eq:Principal part}
    D(s) = \sum_{j = 1}^{h} \frac{A(p, j)}{(s - p)^j} + H_p(s),
    \end{equation}
    where $A(p,h) \neq 0$, and $H_p(s)$ is holomorphic near $s = p$.
    
    \item[(A2)] (Bounds) There are $a < \pi/2$ and a positive unbounded sequence $T_j \to \infty$ with the following properties: The line $\Re(s) = -\varepsilon$ contains no pole of $D(s)$, and there is $\sigma_0 \in (0, \alpha)$ 
  such that every pole $p \neq \alpha$ of $D(s)$ satisfies $\Re(p) < \sigma_0$, and
    \begin{equation}
    \label{eq:Exponential bound}
    D(\sigma + it) \ll_{\sigma} e^{a |t|}, \qquad t \in \mathbb{R}, \ \sigma \in \{-\varepsilon, \sigma_0\}.
    \end{equation}
    Moreover, for every $c > \alpha$ the horizontal segments $\{\sigma \pm iT_j : -\varepsilon \leq \sigma \leq c\}$ contain no pole of $D(s)$, and
  \footnote{Since $|\Gamma(\sigma + it)| \asymp |t|^{\sigma - 1/2} e^{-\pi |t|/2}$ as $|t| \to \infty$, the condition $a < \pi/2$ makes the integrands of the Mellin integrals in \Cref{sec:Generating function} decay exponentially on the lines $\Re(z) = -\varepsilon$ and $\Re(z) = \sigma_0$. The heights $T_j$ are used to close the contour of integration between the poles of $D(s)$, of which there may be infinitely many.}
    \begin{equation}
    \label{eq:Horizontal bound}
    D(\sigma \pm iT_j) \ll_{c} e^{a T_j}, \qquad -\varepsilon \leq \sigma \leq c, \ j \geq 1.
    \end{equation}

    \item[(A3)] (Growth of weights) Let $\Lambda := \{n \in \mathbb{N} : f(n) \neq 0\}$ be the support of $f$ and assume that there exists a real number $M_0 > \alpha/2$, such that for all primes $p$ we have
    \[
    \sum_{n \in \Lambda \setminus p\mathbb{N}} f(n) \geq M_0.
\]
    Since $\Lambda \setminus p\mathbb{N} \subset \Lambda \setminus b\mathbb{N}$ for every prime $p \mid b$, the same inequality holds with $p$ replaced by any $b \geq 2$. If $f$ takes integer values, the condition follows from $|\Lambda \setminus p\mathbb{N}| \geq M_0$.\footnote{The condition (A3) is needed for the minor arcs. If $f$ were supported on the multiples of a prime $p$, then $P(\rho + 2\pi i/p) = P(\rho)$, and the integral \cref{eq:Ordinary Cauchy formula} would not be concentrated near $t = 0$. Quantitatively, \Cref{thm:Minor part} bounds the minor arcs by $\rho^{M_1 + 1} P(\rho)$ for every $M_1 < M_0$, while the main term has size $P(\rho)/\sqrt{\Psi''(\rho)}$, which is $\rho^{1 + \alpha/2} P(\rho)$ up to a power of $\log \rho^{-1}$. This is why $M_0 > \alpha/2$ is needed.}

\end{enumerate}

\subsection{The asymptotic formulas}
\label{sec:Main result}
Let us define the following residual sum: 
\begin{equation}
\label{eq:Delta}
\Delta(\rho) := \sum_{p} \operatorname{Res}_{z = p} \Gamma(z) \zeta(z + 1) D(z) \rho^{-z},
\end{equation}
where $p$ runs over the set 
$\{0\} \cup \big( (\mathcal{A} \cup \mathcal{B}) \setminus \{\alpha\} \big)$, so that $z = 0$ is counted once also when $0 \in \mathcal{A}$, and the sum over $\mathcal{B}$ is understood as in \Cref{thm:Infinite B}.
Define the following polynomials: 
\[
Q_0(s) := \sum_{r = 0}^{k - 1} \frac{F(\alpha, r)}{r!} s^r, 
\qquad Q_{\ell + 1} := (\alpha + \ell) Q_\ell + Q_\ell', \quad \ell \geq 0,
\]
where $F(p, r)$ denotes the coefficient at $(z - p)^{-r-1}$ in $\Gamma(z) \zeta(z + 1) D(z)$. Let us also define the following constant which is repeatedly used throughout the text: 
\[
C := \frac{\alpha (1 + \alpha) A(\alpha, k)\,  \Gamma(\alpha)\, \zeta(1 + \alpha)}{(k - 1)!}.
\]


One of our main results is then stated as follows. 

\begin{theorem}
\label{thm:Explicit saddle point}
Assume the conditions (A1)-(A3) hold, where (A2) holds with some $\sigma_0 < \alpha/2$.
\begin{enumerate}
    \item[(i)] If $k = 1$, then as $n \to \infty$ we have 
    \begin{equation}
    \label{eq:Simple pole formula}
    p_f(n) \sim 
      \frac{\exp\left\{ \left( 1 + \frac{1}{\alpha} \right) \left( \frac{C}{1 + \alpha} \right)^{\frac{1}{1 + \alpha}} n^{\frac{\alpha}{1 + \alpha}} 
        + \Delta\left( \left( \frac{C}{(1 + \alpha) n} \right)^{\frac{1}{1 + \alpha}} \right) 
        \right\}}
      {\sqrt{2 \pi C \left( \frac{(1 + \alpha) n}{C} \right)^{\frac{2 + \alpha}{1 + \alpha}}}}.
    \end{equation}
    If moreover $h_0 = 0$, the residue at $z = 0$ contributes $\frac{D(0)}{1 + \alpha} \log \frac{(1 + \alpha) n}{C} + D'(0)$ to the exponent.
    \item[(ii)] If $k = 2$, let $W$ be the Lambert function, the inverse of the map $w \mapsto w e^{w}$ on $[0, \infty)$, and put
    \[
    v_n := W\left( \frac{(1 + \alpha)^2 e^{(1 + \alpha) \kappa}}{C}\, n \right), \qquad 
  \kappa := \frac{F(\alpha, 0)}{F(\alpha, 1)} + \frac{1}{\alpha}.
    \]
    Then as $n \to \infty$ we have 
    \begin{equation}
    \label{eq:Double pole formula}
    p_f(n) \sim \frac{\exp\left\{ \frac{(1 + \alpha)^{\frac{\alpha - 1}{1 + \alpha}} C^{\frac{1}{1 + \alpha}}}{\alpha} \left( v_n - \frac{1}{\alpha} \right) \left( \frac{n}{v_n} \right)^{\frac{\alpha}{1 + \alpha}} + \Delta\left( \left( \frac{C v_n}{(1 + \alpha)^2 n} \right)^{\frac{1}{1 + \alpha}} \right) \right\}}{\sqrt{2 \pi \frac{C (v_n + 1)}{1 + \alpha} \left( \frac{(1 + \alpha)^2 n}{C v_n} \right)^{\frac{2 + \alpha}{1 + \alpha}}}}.
    \end{equation}
    If moreover $h_0 = 0$, the residue at $z = 0$ contributes $\frac{D(0)}{1 + \alpha} \log \frac{(1 + \alpha)^2 n}{C v_n} + D'(0)$ to the exponent.
    \item[(iii)] If $k \geq 3$, write
    \[
    Q_1\left( \frac{w}{1 + \alpha} \right) = \frac{C}{(1 + \alpha)^{k}} \prod_{i = 1}^{k - 1} (w - w_i), \qquad \mathcal{H}_k(w) := e^{w} \prod_{i = 1}^{k - 1} (w - w_i),
    \]
    with $w_1, \ldots, w_{k - 1} \in \mathbb{C}$, and let $w_n$ be the largest real solution of the equation
    \begin{equation}
    \label{eq:Generalized Lambert equation}
    \mathcal{H}_k(w) = \frac{(1 + \alpha)^{k}}{C}\, n.
    \end{equation}
    Then as $n \to \infty$ we have 
    \begin{equation}
    \label{eq:Higher pole formula}
    p_f(n) \sim \frac{\exp\left\{ e^{\frac{\alpha w_n}{1 + \alpha}} (Q_0 + Q_1)\left( \frac{w_n}{1 + \alpha} \right) + \Delta\left( e^{-\frac{w_n}{1 + \alpha}} \right) \right\}}{\sqrt{2 \pi e^{\frac{(2 + \alpha) w_n}{1 + \alpha}} Q_2\left( \frac{w_n}{1 + \alpha} \right)}}.
    \end{equation}
    If moreover $h_0 = 0$, the residue at $z = 0$ contributes $\frac{D(0)}{1 + \alpha} w_n + D'(0)$ to the exponent. The function $n \mapsto w_n$ is a generalized Lambert function\footnote{In the sense of \cite{MezoBaricz}.}, and with $y_n := (1 + \alpha)^{k} n/C$, it has the following asymptotics: 
    \[
    w_n = \log y_n - (k - 1) \log \log y_n + O\left( \frac{\log \log y_n}{\log y_n} \right), \quad n \to \infty.
    \]
\end{enumerate}
\end{theorem}

\begin{remark}
While this theorem has a better form for the case $\sigma_0 < \alpha/2$, we also show how similar asymptotic formulas are obtained when $\sigma_0 \ge \alpha/2$. 
  We present our results in full generality in \Cref{sec:Further results}, 
  where we also 
  show that for $k \geq 2$ no finite combination of powers of $n$, $\log n$, $\log \log n, \ldots$ approximates the saddle point well enough to give an asymptotic formula for $p_f(n)$. 
\end{remark}

\begin{remark}
The dominant pole $\alpha$ and its order $k$ alone determine the order of magnitude of the main term and in particular, by \Cref{thm:Log asymptotics}, assuming (A1)-(A3) we also obtain the following explicit weak asymptotics:
\[
\log p_f(n) \sim C_1 \cdot n^{\frac{\alpha}{1 + \alpha}} (\log n)^{\frac{k - 1}{1 + \alpha}}, \quad n \to \infty. 
\]
\end{remark}

\begin{remark}
\label{thm:Remark on poles}
In \Cref{thm:Explicit saddle point}, the poles other than $\alpha$ enter the formula through $\Delta$, and a pole $p$ with $\Re(p) > 0$ (real or complex), contributes to $\Delta(\rho)$ a term of size $\rho^{-\Re(p)}$ times a power of $\log \rho^{-1}$, which tends to infinity, so such residues cannot be omitted. For a complex pole $\beta$ the residue contains the factor $\Gamma(\beta)$, which decays like $e^{-\pi|\Im(\beta)|/2}$, so that these terms are often numerically small. A pole $p$ with $\Re(p) < 0$ contributes a term which tends to $0$; finitely many such terms may be omitted, but they can slow down the convergence. A real pole in $(-\varepsilon, 0)$ can also be removed from the region by decreasing $\varepsilon$.
\end{remark}


\subsection{Applications}

We give 
a number of applications of our results, presented in detail in \Cref{sec:Examples}. 
As noted, 
we
extend the saddle-point method of \cite{Tenenbaum2020} in three directions: arbitrary non-negative real weights, dominant poles of arbitrary order, and Dirichlet series with complex poles.

This covers many classical arithmetic weights, such as $d_r(n)$, $d(n)^2$, $2^{\omega(n)}$, $d(n^2)$, $\mu(n)^2$, $(\log n)^j$, $J_r(n)$, $\psi_r(n)$ and $\sigma^{*}(n)$, related to divisors, M\"obius, Jordan, Dedekind functions, as well as partitions into $r$-free parts (see \Cref{sec:Examples}); most of their Dirichlet series are products of shifted Riemann zeta functions and their reciprocals \cite{Ramanujan}, \cite[Chapter~XVII]{HardyWright}. 

It also applies to subgroup counts of groups, through their {\it subgroup zeta functions} \cite{LubotzkySegal}: for $\mathbb{Z}^d$ it gives the number of commuting $(d + 1)$-tuples of permutations \cite{BryanFulman}, whose asymptotics were found in \cite{BringmannFrankeHeim}, and for the Heisenberg group, whose zeta function was computed in \cite{GrunewaldSegalSmith}, it gives an example with a double pole, several real poles and infinitely many complex poles. Our results also apply to counting problems in which the dominant pole has higher order and which, as far as we know, are not covered by earlier results: the number of $n$-dimensional representations of $SU(2)^r$ and of $SU(3)^r$ for $r \geq 2$, and the number of multisets of subrings of $\mathbb{Z}^3$ with indices adding up to $n$ (for $SU(2)^r$, the asymptotics of $\log p_f(n)$ was obtained in \cite{AhmadiGomezAizaWard}); for all of them \Cref{thm:Dominant saddle point}(iii) gives asymptotic formulas in terms of the saddle point $\rho_\alpha$ of the dominant pole (see \Cref{sec:Divisor functions}, \Cref{sec:Representations} and \Cref{sec:Subrings}). 

An example of ordered factorizations in \Cref{sec:Ordered factorizations} show where the method stops: their Dirichlet series $1/(2 - \zeta(s))$ has complex poles accumulating at the line $\Re(s) = \alpha$, and only the asymptotics of $\log p_f(n)$ can be obtained.
 
Our results show how the analytic structure of $D(s)$ determines the growth of $p_f(n)$. The dominant pole $\alpha$ and its order $k$ alone determine the order of magnitude of the main term. 
A pole of higher order produces an additional power of $\log n$ in the exponent, a phenomenon that does not occur in the setting of Meinardus \cite{Meinardus}, where all poles are simple. Such poles arise naturally, already for some basic weights: $2^{\omega(n)}$, $d(n^2)$ and $(\log n)^j$ give poles of order $2$, $3$ and $j + 1$. Zeta functions of nilpotent groups are another natural source: by du Sautoy and Grunewald \cite{duSautoyGrunewald}, their only pole on the line of convergence is the abscissa itself, and its order is often larger than $1$, as for the Heisenberg group. The remaining poles, both real and complex, affect only lower-order terms, as long as they lie strictly to the left of $\alpha$. In particular, complex poles coming from the non-trivial zeros of $\zeta$, as for $J_r$, $\psi_r$, $2^{\omega(n)}$ and $r$-free parts, do not change the main term, and no hypothesis on the location of the zeros is needed in \Cref{thm:Main theorem}; such a hypothesis can enter only through the conditions on $\sigma_0$ in \Cref{thm:Dominant saddle point} (see \Cref{thm:Remark on explicit forms}). The non-negativity of $f$ plays a structural role here: by Landau's theorem it makes $\alpha$ the abscissa of convergence, and it makes $\Psi''(\rho)$ real and positive, so that oscillating contributions cannot destroy the saddle point.

\subsection{Plan of the proof}
The proof extends the saddle-point method of \cite{Tenenbaum2020}. 

In \Cref{sec:Generating function} we write $p_f(n)$ as a Cauchy integral and express $\Psi = \log P$ through the Mellin inversion formula. Shifting the line of integration and computing the residues, we obtain the expansion \cref{eq:Closed form of psi} of $\Psi(\rho)$ as $\rho \to 0^{+}$. In this expansion, every pole $p \neq 0$ of $D(s)$, whether real or complex, contributes $\rho^{-p}$ times a polynomial in $\log \rho^{-1}$, and $z = 0$ contributes a polynomial in $\log \rho^{-1}$.

\Cref{sec:Contributions} contains the main estimates. We first show, by Landau's theorem, that $\alpha$ is the abscissa of convergence of $D(s)$ and that $A(\alpha, k) > 0$. Using (A3) and Delange's Tauberian theorem, we then bound the integral over the minor arcs $\rho^{\chi} \leq |t| \leq \pi$. Next, by shifting the line of integration only to $\Re(z) = \sigma_0$, we prove that the pole $\alpha$ determines the size of every derivative of $\Psi$:
\[
(-1)^{\ell} \Psi^{(\ell)}(\rho) \asymp \rho^{-\alpha - \ell} (\log \rho^{-1})^{k - 1}.
\]
These bounds hold uniformly on vertical lines. On the major arcs $|t| \leq \rho^{\chi}$, they justify the Taylor expansion of $\Psi(\rho + it)$ and the evaluation of the resulting Gaussian integrals. The correction terms $\lambda_{2\mu}(\rho)$ have relative size $O(\rho^{\alpha} (\log \rho^{-1})^{1 - k})$. Hence $p_f(n) \sim e^{\rho n} P(\rho) / \sqrt{2 \pi \Psi''(\rho)}$, which is \Cref{thm:Main theorem}, and by \cref{eq:Closed form of psi} the exponent is given by the residues of all poles.

\Cref{thm:Explicit saddle point} and the other statements of \Cref{sec:Further results} are deduced from \Cref{thm:Main theorem} and the estimates of \Cref{sec:Contributions}; their proofs are given in \Cref{sec:Proofs}. The saddle point $\rho_n$ minimizes $n\rho + \Psi(\rho)$, whose second derivative has size $\rho^{-\alpha - 2} (\log \rho^{-1})^{k - 1}$, and $\Psi(\rho) = \rho^{-\alpha} Q_0(\log \rho^{-1}) + \Delta(\rho) + O(\rho^{\varepsilon})$, where the derivative of $\Delta(\rho)$ is $O(\rho^{-\sigma_0 - 1})$. A Taylor expansion around $\rho_\alpha$ shows that replacing $\rho_n$ by $\rho_\alpha$ changes the exponent by $O(\rho^{\alpha - 2\sigma_0})$, which tends to $0$ if $\sigma_0 < \alpha/2$; this gives \Cref{thm:Dominant saddle point}(i), and solving \cref{eq:Dominant saddle point equation} by the substitution $w = (1 + \alpha) \log \rho^{-1}$ gives \Cref{thm:Explicit saddle point}. When $\sigma_0 < 2\alpha/3$, the next term of this expansion is $-\rho_\alpha^{\alpha + 2} \Delta'(\rho_\alpha)^2/(2Q_2(\log \rho_\alpha^{-1}))$, and the remaining terms are $o(1)$, which gives \Cref{thm:Dominant saddle point}(ii). In the same way, a small change of $\log \rho^{-1}$ changes the exponent by $n\rho$ times the square of this change, up to a constant factor, which gives \Cref{thm:Accuracy of saddle point}(i). Newton's method for the derivative of the explicit function $n\rho + \rho^{-\alpha} Q_0(\log \rho^{-1}) + \Delta(\rho)$, started at $\rho_\alpha$, converges quadratically, so that after a fixed number of steps it reaches the precision required in \Cref{thm:Accuracy of saddle point}(i); this gives \Cref{thm:Dominant saddle point}(iii). Expanding instead the minimum of this function around $\rho_\alpha$ in powers of the derivatives of $\Delta$ gives \Cref{thm:Exponent expansion}. For $k \geq 2$, the saddle-point equation contains the term $(k - 1) \log \log \rho^{-1}$, so that the asymptotic expansion of $\log \rho_n^{-1}$ in powers of $1/\log n$ and $\log \log n/\log n$ does not terminate; this gives \Cref{thm:Accuracy of saddle point}(ii). Finally, \Cref{thm:Log asymptotics} follows from \Cref{thm:Main theorem}, \cref{eq:Derivative asymptotics} and \Cref{thm:Size of rho}.

In \Cref{sec:Examples} we apply these results to various weights mentioned above, grouped by the poles of $D(s)$: no poles other than $\alpha$ and $0$, complex poles coming from the zeros of $\zeta$, and poles $p \neq \alpha$ with $\Re(p) \geq \alpha/2$, which need \Cref{thm:Dominant saddle point}(ii), (iii) or \Cref{thm:Exponent expansion}. \Cref{tab:Examples} compares the resulting formulas with the exact values of $p_f(n)$. 

In \Cref{sec:Conclusion} we give some further remarks, discuss which features of $D(s)$ decide how explicit the formula can be, and limitations of the method.

\subsection{Overview of the asymptotic formulas}
\label{sec:Main result1}
Let us now explain where the above asymptotic formulas come from in some more detail. By Cauchy's formula \cref{eq:Ordinary Cauchy formula}, for every $\rho > 0$,
\[
p_f(n) = \frac{e^{\rho n}}{2 \pi} \int_{-\pi}^{\pi} P(\rho + it) e^{itn}\, dt.
\]
Since $P$ has non-negative coefficients, $|P(\rho + it)| \leq P(\rho)$, so $p_f(n) \leq e^{\rho n} P(\rho) = e^{\rho n + \Psi(\rho)}$ for every $\rho > 0$. This bound is smallest when the derivative of $\rho n + \Psi(\rho)$ vanishes, that is, when $\Psi'(\rho) + n = 0$. For this $\rho$ the linear term in $t$ disappears from the exponent, and near $t = 0$
\[
\Psi(\rho + it) + itn = \Psi(\rho) - \frac{1}{2} \Psi''(\rho) t^2 + O\big( |\Psi'''(\rho)| |t|^3 \big)
\]
so that the integral over a small arc around $t = 0$ is close to the Gaussian integral $P(\rho) \sqrt{2\pi/\Psi''(\rho)}$. This gives the denominator $\sqrt{2 \pi \Psi''(\rho)}$ in the formula below. The rest of the circle is negligible because of condition (A3). Under the conditions (A1)-(A3) these steps can be made rigorous, and we obtain the saddle-point formula
\[
p_f(n) \sim \frac{e^{\rho_n n} P(\rho_n)}{\sqrt{2 \pi \Psi''(\rho_n)}}
\]
where $\rho_n$ is the unique positive solution of $\Psi'(\rho) + n = 0$ (\Cref{thm:Main theorem}). This formula has the classical form, but $\rho_n$ is defined only implicitly, through all poles of $D(s)$. To make it explicit, we need $\Psi(\rho)$ near $\rho = 0$ and a good approximation of $\rho_n$. 
Expanding the logarithm in \cref{eq:Generating function} gives $\Psi(\rho) = \sum_{n, m \geq 1} f(n) e^{-mn\rho}/m$. Since $\int_0^\infty e^{-mn\rho} \rho^{z - 1}\, d\rho = \Gamma(z) (mn)^{-z}$ for $\Re(z) > 0$, the Mellin transform of $\Psi$ is $\Gamma(z) \zeta(z + 1) D(z)$. We invert it (\Cref{thm:Mellin inversion formula}) and move the line of integration to the left, past the poles of $D(z)$ and the pole of $\Gamma(z) \zeta(z + 1)$ at $z = 0$. By \Cref{thm:Closed form of psi}, $\Psi(\rho) = \log P(\rho)$ is then the sum of the residues of $\Gamma(z) \zeta(z + 1) D(z) \rho^{-z}$ at the points of $\{0\} \cup \mathcal{A} \cup \mathcal{B}$, up to an error $O(\rho^{\varepsilon})$.

These residues have a simple form. Near a point $p$ we write $\rho^{-z} = \rho^{-p} e^{(z - p) \log \rho^{-1}}$, so the residue at $p$ is $\rho^{-p}$ times a polynomial in $\log \rho^{-1}$, whose coefficients are the Laurent coefficients of $\Gamma(z) \zeta(z + 1) D(z)$ at $z = p$. To write them down, for a point $p$ with $\Re(p) > -1$, let $\phi(p, \nu)$ denote the coefficients of the Laurent expansion of $\Gamma(z) \zeta(z + 1)$ at $z = p$,
\[
\Gamma(p + s) \zeta(p + 1 + s) = \sum_{\nu} \phi(p, \nu) s^{\nu}.
\]
For $p \neq 0$ both factors are holomorphic at $s = 0$, so that $\phi(p, \nu) = 0$ for $\nu < 0$ and
\[
\phi(p,\nu) = \sum_{i = 0}^{\nu} \frac{\Gamma^{(i)}(p) \zeta^{(\nu - i)}(p + 1)}{i!(\nu - i)!}, \quad \nu \geq 0.
\]
For $p = 0$ the product has a double pole, and the expansion starts at $\nu = -2$, with $\phi(0, -2) = 1$ and $\phi(0, -1) = 0$ (see the proof of \Cref{thm:Res at 0}). If $p \neq 0$ is a pole of $D(s)$ of order $h$, with principal part \cref{eq:Principal part}, then for an integer $r \geq 0$ the coefficient of $(z - p)^{-r - 1}$ in $\Gamma(z) \zeta(z + 1) D(z)$ is
\[
F(p,r) := \sum_{r + 1 \leq j \leq h} \phi(p, j - r - 1) \cdot A(p,j)
\]
and the residue at $z = p$ is $\rho^{-p} \sum_{r = 0}^{h - 1} \frac{F(p, r)}{r!} (\log \rho^{-1})^{r}$ (\Cref{thm:Res at p}). Thus a pole $p \neq 0$ of order $h$ contributes $\rho^{-p}$ times a polynomial of degree at most $h - 1$ in $\log \rho^{-1}$. The point $z = 0$ contributes a polynomial of degree at most $h_0 + 1$ (exactly $h_0 + 1$ if $h_0 \geq 1$) in $\log \rho^{-1}$, which equals $D(0) \log \rho^{-1} + D'(0)$ if $h_0 = 0$; the corresponding numbers $F(0, r)$ are defined in \Cref{thm:Res at 0}.

The residue at the dominant pole $z = \alpha$ is $\rho^{-\alpha} Q_0(\log \rho^{-1})$, where
\[
Q_0(s) := \sum_{r = 0}^{k - 1} \frac{F(\alpha, r)}{r!} s^r
\]
The saddle-point equation involves $\Psi'$, and the denominator involves $\Psi''$, so we also need the derivatives of this term. Since $\frac{d}{d\rho} \big( \rho^{-\alpha - \ell} Q_\ell(\log \rho^{-1}) \big) = -\rho^{-\alpha - \ell - 1} \big( (\alpha + \ell) Q_\ell + Q_\ell' \big)(\log \rho^{-1})$, we define
\[
Q_{\ell + 1} := (\alpha + \ell) Q_\ell + Q_\ell', \quad \ell \geq 0, \quad \text{so that} \quad (-1)^{\ell} \frac{d^{\ell}}{d\rho^{\ell}} \big( \rho^{-\alpha} Q_0(\log \rho^{-1}) \big) = \rho^{-\alpha - \ell} Q_\ell(\log \rho^{-1}).
\]
Thus $Q_1 = \alpha Q_0 + Q_0'$ and $Q_2 = \alpha(1 + \alpha) Q_0 + (1 + 2\alpha) Q_0' + Q_0''$. The polynomial $Q_1$ appears in the saddle-point equation \cref{eq:Dominant saddle point equation} of the residue at $z = \alpha$ below, and $Q_2$ in the main term of $\Psi''$. Every $Q_\ell$ is a polynomial of degree $k - 1$, and its leading coefficient comes only from the highest term $A(\alpha, k) (s - \alpha)^{-k}$ of the principal part, since $F(\alpha, k - 1) = \Gamma(\alpha) \zeta(1 + \alpha) A(\alpha, k)$. We collect these factors in the constant
$
C := \frac{A(\alpha, k) \alpha (1 + \alpha) \Gamma(\alpha) \zeta(1 + \alpha)}{(k - 1)!},
$
which is normalized so that the leading coefficients of $Q_0$, $Q_1$ and $Q_2$ are $\frac{C}{\alpha(1 + \alpha)}$, $\frac{C}{1 + \alpha}$ and $C$, by \cref{eq:Leading coefficient}. For example, for ordinary partitions we have $D(s) = \zeta(s)$, $\alpha = k = 1$ and $C = \pi^2/3$, and the residue at $z = 1$ is $\rho^{-1} Q_0 = \pi^2/(6\rho)$.

We split off the residue at the dominant pole. Recall that 
\begin{equation*}
\Delta(\rho) = \sum_{p} \operatorname{Res}_{z = p} \Gamma(z) \zeta(z + 1) D(z) \rho^{-z}
\end{equation*}
where $p$ runs over the set $\{0\} \cup \big( (\mathcal{A} \cup \mathcal{B}) \setminus \{\alpha\} \big)$, so that $z = 0$ is counted once also when $0 \in \mathcal{A}$, and the sum over $\mathcal{B}$ is understood as in \Cref{thm:Infinite B}. By \Cref{thm:Closed form of psi}, $\Psi(\rho) = \rho^{-\alpha} Q_0(\log \rho^{-1}) + \Delta(\rho) + O(\rho^{\varepsilon})$.

The first term dominates. For large $n$, let $\rho_\alpha = \rho_\alpha(n)$ be the smallest positive solution of the saddle-point equation of this term alone,
\begin{equation}
\label{eq:Dominant saddle point equation}
\rho^{-\alpha - 1} Q_1(\log \rho^{-1}) = n
\end{equation}
so that $\rho_\alpha \to 0$ as $n \to \infty$. The poles in $\Delta$ move the true saddle point $\rho_n$ away from $\rho_\alpha$ by a relative amount $O(\rho_\alpha^{\alpha - \sigma_0})$, up to a power of $\log \rho_\alpha^{-1}$. On the other hand, by \Cref{thm:Accuracy of saddle point}, the saddle-point formula remains correct for an approximation of $\rho_n$ exactly when the relative error is $o\big( (n\rho_n)^{-1/2} \big)$, which is of order $\rho_n^{\alpha/2}$ up to a power of $\log \rho_n^{-1}$. Hence $\rho_\alpha$ itself may be used when $\sigma_0 < \alpha/2$. The left-hand side of \cref{eq:Dominant saddle point equation} is $-\frac{d}{d\rho} \rho^{-\alpha} Q_0(\log \rho^{-1})$, and the substitution $w = (1 + \alpha) \log \rho^{-1}$ turns $\rho^{-\alpha - 1}$ into $e^{w}$ and $Q_1(\log \rho^{-1})$ into a polynomial in $w$ of degree $k - 1$. For $k = 1$ and $k = 2$ the equation can then be solved explicitly, and for $k \geq 3$ its solution is a generalized Lambert function of $n$. This leads to our main theorem. 

The substitution $w = (1 + \alpha) \log \rho^{-1}$ in \Cref{thm:Explicit saddle point}(iii) turns $\rho^{-\alpha - 1}$ into $e^{w}$ and $Q_1(\log \rho^{-1})$ into a polynomial in $w$, so that \cref{eq:Dominant saddle point equation} becomes \cref{eq:Generalized Lambert equation}. The same substitution works for every $k$, and \cref{eq:Higher pole formula} also holds for $k = 1$ and $k = 2$, where it reduces to (i) and (ii). For $k = 1$ we have $\mathcal{H}_1(w) = e^{w}$, and for $k = 2$ we have $\mathcal{H}_2(w) = e^{w} (w - w_1)$ with $w_1 = -(1 + \alpha)\kappa$, so that $v_n = w_n - w_1$. The condition $\sigma_0 < \alpha/2$ is needed only to use $\rho_\alpha$ itself, and the cases $\sigma_0 \geq \alpha/2$ are covered as well. The identity $\rho_\alpha = e^{-w_n/(1 + \alpha)}$ does not use the condition on $\sigma_0$, and for every $\sigma_0 < \alpha$ a fixed number of Newton steps starting from $\rho_\alpha$, namely the least $J \geq 0$ with $2^{J}(\alpha - \sigma_0) > \alpha/2$, gives an asymptotic formula for $p_f(n)$ (\Cref{thm:Dominant saddle point}(iii)). For $\sigma_0 < 2\alpha/3$ a single correction term in the exponent suffices (\Cref{thm:Dominant saddle point}(ii)), and an expansion of the exponent around $\rho_\alpha$ to a fixed order gives a formula that is explicit in $\rho_\alpha$ (\Cref{thm:Exponent expansion}). In \Cref{sec:Examples} these results are applied when $\sigma_0 \geq \alpha/2$. We state them, together with the saddle-point formula, the asymptotics of $\log p_f(n)$ and the precision required of the saddle point, in \Cref{sec:Further results}.


\subsection{Comparison with previous works}
Our argument extends and differs from that of \cite{Tenenbaum2020} in the following respects.
First, 
  the weights $f(n)$ are arbitrary non-negative real numbers, and $D(s)$ may have several real poles of arbitrary orders as well as complex poles. 
Secondly, 
  since the remaining poles contribute to $\Psi$, possibly through oscillating terms $\rho^{-\beta}$, the domination of the pole $\alpha$ has to be proved. This is done by the partial shift of the line of integration described above. It requires the assumption \cref{eq:Complex poles to the left} in (A1); an example in \Cref{rem:ex} 
  shows that without it the pole $\alpha$ need not dominate $\Psi$. 
Thirdly, 
  since $\alpha$ is a pole of order $k$, every derivative of $\Psi$ carries a factor $(\log \rho^{-1})^{k - 1}$. All estimates are arranged so that these logarithms are absorbed, and they enter the formula through the residue $\rho^{-\alpha} Q_0(\log \rho^{-1})$ at $z = \alpha$. 
Fourthly, 
  the Gaussian integrals require $\Psi^{(\ell)}(\rho)$ to be real and $\Psi''(\rho)$ to be positive. When complex poles are present, this does not follow from the closed form of $\Psi$, and we deduce it from the non-negativity of $f$ instead. 
Finally, 
  we prove an explicit bound for the coefficients $\lambda_{2\mu}(\rho)$, whose contribution is stated to be of lower order without proof in \cite{Tenenbaum2020} (for simple real poles, a bound of this kind is given in \cite[Lemma~3.6]{Bridges2023}); this bound also gives the relative error on the major arcs.

Compared with \cite{Bridges2023}, where the poles of the Dirichlet series are real and simple, and an explicit formula is obtained when there are two positive poles \cite[Theorem~4.4]{Bridges2023}, we allow poles of higher order and complex poles. Our condition (A3) is a weighted form of the condition on $\Lambda$ used there. 
The results of \cite{GranovskyStark2012, GranovskyStark2015} also assume that all poles are simple and lie on the non-negative real axis; poles of higher order, such as those of $d(n)$, $2^{\omega(n)}$ or the Heisenberg group, are not covered by these works. 
In \cite{AhmadiGomezAizaWard} poles of higher order are allowed, but for a pole of order at least $2$ at $s = 1$ or $s = 2$, as for $\zeta(s)^{k}$ with $k \geq 2$, they obtain only the asymptotics of $\log p_f(n)$; 
we also note that their families with $i = k = 0$, in which the only pole of $D(s)$ is at $s = 0$, lie outside our setting, since we assume $\alpha > 0$.

For the application of commuting tuples of permutations, our work in a sense complements the approach from \cite{BringmannFrankeHeim}. 
The work \cite{BringmannFrankeHeim} treat the Dirichlet series $D(s) = \zeta(s) \zeta(s - 1) \cdots \zeta(s - \ell + 2)$, whose poles are all real and simple, and obtain an asymptotic expansion of the number of commuting $\ell$-tuples to all orders in $n^{-1/\ell}$, with the constants in the exponent given in closed form for $\ell = 3$ and through Lagrange inversion for $\ell \geq 4$, together with log-concavity for large $n$. While in these examples, we do not obtain such complete expansions, for $\ell = 3$ our \Cref{thm:Dominant saddle point}(ii) gives the same formula, and for every $\ell \geq 4$ \Cref{thm:Dominant saddle point}(iii) gives a closed asymptotic formula after $\lfloor \log_2 (\ell - 1) \rfloor$ Newton steps starting from $\rho_\alpha = (C/(\ell n))^{1/\ell}$. On the other hand, our results apply to Dirichlet series with poles of higher order and with complex poles, which lie outside the scope of \cite{Bridges2023, BringmannFrankeHeim, GranovskyStark2012, GranovskyStark2015}, such as those of $d(n)$, $2^{\omega(n)}$, $r$-free parts, $\varphi(n)$, the Heisenberg group, the representations of $SU(3)^r$ and the subrings of $\mathbb{Z}^3$. Moreover, \Cref{thm:Accuracy of saddle point}(ii) shows that for poles of higher order the saddle point cannot be replaced by finitely many terms of an expansion as in \cite{BringmannFrankeHeim}: for $k \geq 2$ no finite combination of powers of $n$, $\log n$, $\log \log n, \ldots$ gives an asymptotic formula, while \Cref{thm:Dominant saddle point}(iii) still gives one in terms of $\rho_\alpha$.

\subsection{Notation}
\label{sec:Notation}
The following 
notation is used throughout the paper.
\begin{itemize}
    \item $\mathbb{N} = \{1, 2, 3, \ldots\}$, and for $b \in \mathbb{N}$, $b\mathbb{N}$ is the set of positive multiples of $b$.
    \item For functions $F$ and $G \geq 0$, we write $F \ll G$, $G \gg F$ or $F = O(G)$ if there is a constant $c > 0$ such that $|F| \leq c\, G$ in the stated range of the variables. A subscript, as in $\ll_{M_0}$ or $O_N$, indicates that $c$ depends on the parameter in the subscript. Unless stated otherwise, $c$ may depend on $f$ and on the fixed parameters $\chi$, $N_0$, $\ell$, but not on $\rho$, $t$ or $n$.
    \item $F \asymp G$ means that both $F \ll G$ and $G \ll F$ hold.
    \item $F \sim G$ means that $F/G \to 1$, and $F = o(G)$ means that $F/G \to 0$, as $n \to \infty$ or $\rho \to 0^{+}$, depending on the context.
    \item For real $y$, $\lceil y \rceil$ is the smallest integer $\geq y$, and $\|y\|$ is the distance from $y$ to the nearest integer.
    \item $F^{(\ell)}$ is the $\ell$-th derivative of $F$.
    \item $\Gamma$ is Euler's Gamma function and $\zeta$ is the Riemann zeta function. $D(0)$ and $D'(0)$ are the values at $s = 0$ of the meromorphic continuation of $D(s)$.
    \item $\operatorname{Res}_{z = p} F(z)$ is the residue of $F$ at $z = p$. The integral $\int_{c - i\infty}^{c + i\infty}$ is taken upwards along the line $\Re(z) = c$.
    \item $\rho_n$ denotes the solution of the saddle-point equation $\Psi'(\rho) + n = 0$ from \Cref{thm:Rho existence}. In \Cref{sec:Contributions} we write $\rho = \rho_n$; elsewhere $\rho$ is a positive real variable.
\end{itemize}

\section{The saddle-point formula and further results}
\label{sec:Further results}
In this section we state our results in full generality. \Cref{thm:Explicit saddle point} is a special case of them, and the proofs are given in \Cref{sec:Contributions} and \Cref{sec:Proofs}.

\subsection{The saddle-point formula}
\label{sec:Saddle point formula}
Recall that we denote $\Psi(\rho) := \log P(\rho)$. 
The following theorem makes the heuristic argument of \Cref{sec:Main result1} rigorous. It holds for every $\sigma_0$ in (A2), and no explicit form of the saddle point is needed.

\begin{theorem}
\label{thm:Main theorem}
Assume the conditions (A1)-(A3) hold. Then, as $n \to \infty$,
\begin{equation}
\label{eq:Main formula}
p_f(n) \sim \frac{e^{\rho n} P(\rho)}{\sqrt{2 \pi \Psi''(\rho)}}, 
\end{equation}
where $\rho = \rho_n$ is the unique positive solution of the saddle-point equation $\Psi'(\rho) + n = 0$ (see \Cref{thm:Rho existence}).
\end{theorem}

\noindent
By \Cref{thm:Closed form of psi} below, the exponent $\rho n + \Psi(\rho)$ is the sum of $\rho n$ and the residues of all poles of $D(s)$, up to an error $O(\rho^{\varepsilon})$. The saddle point $\rho_n$, however, depends on all these poles, and the next subsection shows how to replace it by $\rho_\alpha$.

\subsection{The dominant saddle point}
\label{sec:Dominant saddle point}
\Cref{thm:Explicit saddle point} uses the solution $\rho_\alpha$ of \cref{eq:Dominant saddle point equation} in place of $\rho_n$. Once $\rho_\alpha$ has been computed, it remains to control the error caused by this replacement. The main term $\rho^{-\alpha} Q_0(\log \rho^{-1})$ of $\Psi$ has size $\asymp n\rho$ near the saddle point, and $\Delta$ is smaller by a factor $\rho^{\alpha - \sigma_0}$, up to a power of $\log \rho^{-1}$. The poles in $\Delta$ move the minimum point of $n\rho + \Psi(\rho)$ from $\rho_\alpha$ to $\rho_n$, by a relative amount of order $\rho^{\alpha - \sigma_0}$. Since $\rho_n$ is a minimum point, a point at relative distance $\delta$ from $\rho_n$ changes the value of $n\rho + \Psi(\rho)$ only by about $n\rho\, \delta^2$. At $\rho_\alpha$ the value is therefore off by about $n\rho \cdot \rho^{2(\alpha - \sigma_0)}$, which is of order $\rho^{\alpha - 2\sigma_0}$ up to a power of $\log \rho^{-1}$. This tends to $0$ exactly when $\sigma_0 < \alpha/2$, and then $\rho_\alpha$ may be used as it is; this is \Cref{thm:Explicit saddle point}. When $\sigma_0 \geq \alpha/2$, there are two ways to repair the error, and either one suffices.

The first way corrects the point. Newton's method for the derivative of the exponent, started at $\rho_\alpha$, squares the relative error at each step. After $J$ steps the relative error is about $\rho^{2^{J}(\alpha - \sigma_0)}$, so that the error in the value is about $n\rho \cdot \rho^{2^{J + 1}(\alpha - \sigma_0)}$. This tends to $0$ when $2^{J}(\alpha - \sigma_0) > \alpha/2$, and the exponent is then evaluated exactly at the last Newton point. This is \Cref{thm:Dominant saddle point}(iii).

The second way corrects the value. One keeps the point $\rho_\alpha$ and expands the minimum value of the exponent around $\rho_\alpha$ in powers of $\rho^{\alpha - \sigma_0}$. The correction $T_m$ has size about $n\rho \cdot \rho^{m(\alpha - \sigma_0)}$, so that $T_m$ is needed exactly when $\sigma_0 \geq \alpha (1 - 1/m)$: the term $T_2$ from $\sigma_0 = \alpha/2$ on, $T_3$ from $\sigma_0 = 2\alpha/3$ on, and so on. This is \Cref{thm:Exponent expansion}, and its first case, with $T_2$ alone, is \Cref{thm:Dominant saddle point}(ii).

Newton's method needs fewer steps, about $\log_2$ of the number of terms $T_m$, but its result is an iterated expression in $\rho_\alpha$, while the expansion gives a formula written directly in $\rho_\alpha$. In both cases the denominator only requires $\rho_\alpha$, since there a factor $1 + o(1)$ suffices. Ignoring the powers of $\log \rho^{-1}$, the least number $J$ of Newton steps and the corrections $T_m$ that are needed are as follows.
\begin{center}
\begin{tabular}{lccccc}
\toprule
$\sigma_0$ in & $[0, \frac{\alpha}{2})$ & $[\frac{\alpha}{2}, \frac{2\alpha}{3})$ & $[\frac{2\alpha}{3}, \frac{3\alpha}{4})$ & $[\frac{3\alpha}{4}, \frac{4\alpha}{5})$ & $[\frac{4\alpha}{5}, \frac{5\alpha}{6})$ \\
\midrule
Newton steps $J$ & $0$ & $1$ & $1$ & $2$ & $2$ \\
corrections & none & $T_2$ & $T_2, T_3$ & $T_2, \ldots, T_4$ & $T_2, \ldots, T_5$ \\
\bottomrule
\end{tabular}
\end{center}
In \Cref{thm:Exponent expansion}, $\sigma_0$ may be replaced by the exponent $\sigma_1$ defined there, which is often smaller. In borderline cases the powers of $\log \rho^{-1}$ can decide: for the subrings of $\mathbb{Z}^3$ in \Cref{sec:Subrings}, the largest other pole is $2\alpha/3$, on the boundary, and $N = 2$ is still admissible. The two ways are stated in the following two theorems. In \Cref{thm:Dominant saddle point}, part (i) is \Cref{thm:Explicit saddle point} before the substitution $w = (1 + \alpha) \log \rho^{-1}$.

\begin{theorem}
\label{thm:Dominant saddle point}
Assume the conditions (A1)-(A3) hold, and let $\rho_\alpha$ be the smallest positive solution of the equation 
  \begin{equation}
\label{eq:Dominant saddle point equation1}
\rho^{-\alpha - 1} Q_1(\log \rho^{-1}) = n.
\end{equation}

\begin{enumerate}
    \item[(i)] If (A2) holds with some $\sigma_0 < \alpha/2$, then, as $n \to \infty$,
    \begin{equation}
    \label{eq:Dominant formula}
    p_f(n) \sim \frac{\exp\left\{ \rho_\alpha n + \rho_\alpha^{-\alpha} Q_0(\log \rho_\alpha^{-1}) + \Delta(\rho_\alpha) \right\}}{\sqrt{2 \pi \rho_\alpha^{-\alpha - 2} Q_2(\log \rho_\alpha^{-1})}}.
    \end{equation}
    \item[(ii)] If (A2) holds with some $\sigma_0 < 2\alpha/3$, then, as $n \to \infty$,
    \begin{equation}
    \label{eq:Dominant formula corrected}
    p_f(n) \sim \frac{\exp\left\{ \rho_\alpha n + \rho_\alpha^{-\alpha} Q_0(\log \rho_\alpha^{-1}) + \Delta(\rho_\alpha) - \dfrac{\rho_\alpha^{\alpha + 2} \Delta'(\rho_\alpha)^2}{2 Q_2(\log \rho_\alpha^{-1})} \right\}}{\sqrt{2 \pi \rho_\alpha^{-\alpha - 2} Q_2(\log \rho_\alpha^{-1})}}.
    \end{equation}
    \item[(iii)] Let $\sigma_0$ be as in (A2), and let $J \geq 0$ be an integer with $\sigma_0 < \left(1 - \frac{1}{2^{J+1}}\right)\alpha$.
  Put
    \[
    \Phi_\Delta(x) := nx + x^{-\alpha} Q_0(\log x^{-1}) + \Delta(x),\quad x > 0.
    \]
    Let $x_0 := \rho_\alpha$, and $x_{j + 1} := x_j - \Phi_\Delta'(x_j)/\Phi_\Delta''(x_j)$ for $0 \leq j < J$. Then for large $n$ the points $x_j$ are well defined, and as $n \to \infty$,
    \begin{equation}
    \label{eq:Newton formula}
    p_f(n) \sim \frac{\exp\left\{ \Phi_\Delta(x_J) \right\}}{\sqrt{2 \pi \Phi_\Delta''(x_J)}} \sim \frac{\exp\left\{ \Phi_\Delta(x_J) \right\}}{\sqrt{2 \pi \rho_\alpha^{-\alpha - 2} Q_2(\log \rho_\alpha^{-1})}}.
    \end{equation}
\end{enumerate}
In all three cases, the residues at finitely many poles $p$ with $\Re(p) < 0$ may be omitted from $\Delta$.
\end{theorem}

\noindent
Since $\Phi_\Delta(\rho_\alpha)$ is the exponent in \cref{eq:Dominant formula}, part (i) is the case $J = 0$ of part (iii). One Newton step suffices for $\sigma_0 < 3\alpha/4$, and the Taylor expansion of $\Phi_\Delta(x_1)$ around $\rho_\alpha$ gives back the exponent in \cref{eq:Dominant formula corrected}. In general the least admissible $J$ is $\max\{0, \lfloor \log_2 (\alpha/(2(\alpha - \sigma_0))) \rfloor + 1\}$, which grows as $\sigma_0$ approaches $\alpha$ but is fixed for a given $D(s)$, and \cref{eq:Newton formula} holds for every larger $J$ as well. For $k = 1$, $\rho_\alpha = (C/((1 + \alpha) n))^{1/(1 + \alpha)}$ is an explicit power of $n$, so \cref{eq:Newton formula} is a closed formula for $p_f(n)$ whenever $\Delta$ is explicit, and for $k = 2$, $\rho_\alpha$ is given by the Lambert function as in \Cref{thm:Explicit saddle point}(ii). Numerically, the first form in \cref{eq:Newton formula} is more accurate than the second, and when $\alpha - \sigma_0$ is small, one or two additional steps improve the approximation considerably.

The exponents in \cref{eq:Dominant formula} and \cref{eq:Dominant formula corrected} are the first two cases of an expansion of the exponent around $\rho_\alpha$ that can be carried to any order. For every $\sigma_0$ it gives a formula that is explicit in $\rho_\alpha$.

\begin{theorem}
\label{thm:Exponent expansion}
Assume the conditions (A1)-(A3) hold, and let $\sigma_1 < \alpha$ and $h_1 \geq 0$ be such that, for every $\ell \geq 0$,
\begin{equation}
\label{eq:Delta bound}
\Delta^{(\ell)}(\rho) \ll_{\ell} \rho^{-\sigma_1 - \ell} (\log \rho^{-1})^{h_1}, \quad 0 < \rho \leq 1/2.
\end{equation}
Let $N \geq 1$ be an integer such that
\[
\rho^{(N + 1)(\alpha - \sigma_1) - \alpha} (\log \rho^{-1})^{(N + 1) h_1 - N(k - 1)} \to 0, \quad \rho \to 0^{+}.
\]
With $u_\alpha := \log \rho_\alpha^{-1}$, put
\[
q_j := Q_j(u_\alpha), \qquad b_j := (-\rho_\alpha)^{j} \Delta^{(j)}(\rho_\alpha), \quad j \geq 1.
\]
Then, as $n \to \infty$,
\begin{equation}
\label{eq:Expanded formula}
p_f(n) \sim \frac{\exp\left\{ n \rho_\alpha + \rho_\alpha^{-\alpha} Q_0(u_\alpha) + \Delta(\rho_\alpha) + T_2 + \cdots + T_N \right\}}{\sqrt{2 \pi \rho_\alpha^{-\alpha - 2} \left( q_2 + \rho_\alpha^{\alpha} (b_2 - b_1 q_3/q_2) \right)}},
\end{equation}
where
  $T_m$ is the coefficient of $t^m$ in $U(\eta(t), t)$ given by 
\[
U(\eta, t) := \rho_\alpha^{-\alpha} \sum_{j \geq 2} \frac{q_j}{j!} (-\eta)^{j} + t \sum_{j \geq 1} \frac{b_j}{j!} (-\eta)^{j}
\]
and $\eta(t)$ is the power series in $t$ without constant term with $\frac{\partial U}{\partial \eta}(\eta(t), t) = 0$; e.g we have 
\[
T_2 = -\frac{\rho_\alpha^{\alpha} b_1^2}{2 q_2}, \qquad T_3 = \frac{\rho_\alpha^{2\alpha} b_1^2 (3 q_2 b_2 - q_3 b_1)}{6 q_2^3}
\]
\[
T_4 = -\frac{\rho_\alpha^{3\alpha} b_1^2 \left( 4 q_2^2 b_1 b_3 + 12 q_2^2 b_2^2 - 12 q_2 q_3 b_1 b_2 - q_2 q_4 b_1^2 + 3 q_3^2 b_1^2 \right)}{24 q_2^5}.
\]  
  The term $\rho_\alpha^{\alpha} (b_2 - b_1 q_3/q_2)$ in the denominator may be omitted.
\end{theorem}

\noindent
The proof is given in \Cref{sec:Proofs}. By \cref{eq:Remainder bounds}, \cref{eq:Delta bound} holds with $\sigma_1 = \sigma_0$ and $h_1 = 0$. Then $N = 1$ is admissible if $\sigma_0 < \alpha/2$, and $N = 2$ if $\sigma_0 < 2\alpha/3$, and without the term $\rho_\alpha^{\alpha} (b_2 - b_1 q_3/q_2)$ in the denominator \cref{eq:Expanded formula} becomes \cref{eq:Dominant formula} and \cref{eq:Dominant formula corrected}, respectively. Often \cref{eq:Delta bound} holds with a smaller $\sigma_1$. Suppose that \cref{eq:Exponential bound} holds on a line $\Re(s) = \sigma'$ with $-\varepsilon < \sigma' < \sigma_0$ that contains no pole of $D(s)$, and that only finitely many poles $p \neq \alpha$ of $D(s)$ satisfy $\Re(p) > \sigma'$. Shifting the line of integration from $\Re(z) = \sigma_0$ to $\Re(z) = \sigma'$ as in the proof of \Cref{thm:Psi decomposition}, we see that $\Delta$ is the sum of the residues of $\Gamma(z) \zeta(z + 1) D(z) \rho^{-z}$ at these poles and of a function whose $\ell$-th derivative is $\ll_\ell \rho^{-\sigma' - \ell}$. By \Cref{thm:Res at p}, \cref{eq:Delta bound} then holds with $\sigma_1$ the largest of $\sigma'$ and the numbers $\Re(p)$, and with $h_1 + 1$ the largest order of a pole $p$ with $\Re(p) = \sigma_1$, or $h_1 = 0$ if there is no such pole. In the examples of \Cref{sec:Examples} in which \Cref{thm:Exponent expansion} is used, $\sigma_1$ is the second largest real part of a pole of $D(s)$. The least admissible $N$ is about $\alpha/(\alpha - \sigma_1)$, while the number of Newton steps in \Cref{thm:Dominant saddle point}(iii) grows only like its logarithm; on the other hand, \cref{eq:Expanded formula} is explicit in $\rho_\alpha$, and hence in $n$ when $\rho_\alpha$ is.

\subsection{Logarithmic asymptotics}
\label{sec:Precision1}
\noindent
\Cref{thm:Dominant saddle point} and \Cref{thm:Exponent expansion} use the residues collected in $\Delta$. The logarithm of $p_f(n)$ needs much less: without any condition on $\sigma_0$, it is determined by the size of $\rho_n$, that is, by $\alpha$, $k$ and $C$.

\begin{corollary}
\label{thm:Log asymptotics}
Assume the conditions (A1)-(A3) hold. Then, as $n \to \infty$,
\[
\log p_f(n) \sim \left( 1 + \frac{1}{\alpha} \right) \left( \frac{C}{(1 + \alpha)^{k}} \right)^{\frac{1}{1 + \alpha}} n^{\frac{\alpha}{1 + \alpha}} (\log n)^{\frac{k - 1}{1 + \alpha}}.
\]
\end{corollary}

\subsection{Precision of the saddle point}
\label{sec:Precision}
\noindent
On the other hand, \cref{eq:Main formula} is sensitive to the precision of the saddle point. The following result gives the precision which is needed, and shows that for $k \geq 2$ it cannot be reached by finitely many terms of an asymptotic expansion of $\log \rho_n^{-1}$.

\begin{theorem}
\label{thm:Accuracy of saddle point}
Assume the conditions (A1)-(A3) hold, and let $\tilde{\rho} = \tilde{\rho}(n) > 0$ with $\tilde{\rho} \to 0$ as $n \to \infty$.
\begin{enumerate}
    \item[(i)] We have
    \begin{equation}
    \label{eq:Perturbed saddle point}
    p_f(n) \sim \frac{e^{\tilde{\rho} n} P(\tilde{\rho})}{\sqrt{2 \pi \Psi''(\tilde{\rho})}}
    \end{equation}
    if and only if
    \[
    \log \tilde{\rho}^{-1} - \log \rho_n^{-1} = o\left( n^{-\frac{\alpha}{2(1 + \alpha)}} (\log n)^{-\frac{k - 1}{2(1 + \alpha)}} \right).
    \]
    \item[(ii)] Let $k \geq 2$, and let $\tilde{\rho} = e^{-L(n)}$ with $L(n) \to \infty$, where
    \[
    L(n) = \sum_{j \in \mathcal{J}} c_j\, n^{a_{j,0}} (\log n)^{a_{j,1}} (\log \log n)^{a_{j,2}} \cdots (\log_r n)^{a_{j,r}}
    \]
    with a finite index set $\mathcal{J}$, $r \in \mathbb{N}$, real numbers $c_j$ and $a_{j,i}$, and $\log_r$ the $r$-fold iterated logarithm. Then there is $c > 0$ such that
    \[
    \left| L(n) - \log \rho_n^{-1} \right| \geq  (\log n)^{-c}
    \]
    for all large $n$. In particular, \cref{eq:Perturbed saddle point} fails, and in fact
    \[
    p_f(n) = o\left( \frac{e^{\tilde{\rho} n} P(\tilde{\rho})}{\sqrt{2 \pi \Psi''(\tilde{\rho})}} \right), \qquad n \to \infty.
    \]
\end{enumerate}
\end{theorem}

\noindent
For $k = 1$ the situation is different: then $\rho_\alpha = (C/((1 + \alpha) n))^{1/(1 + \alpha)}$, and by the proof of \Cref{thm:Dominant saddle point}(iii), \cref{eq:Perturbed saddle point} holds for the point $\tilde{\rho} = x_J$ obtained from $\rho_\alpha$ by finitely many Newton steps; if $\sigma_0 < \alpha/2$, one may take $\tilde{\rho} = \rho_\alpha$. For $k \geq 2$, part (ii) excludes finite sums of monomials, but it says nothing about functions such as $\log(\log n - \log \log n)$, whose expansions have infinitely many terms. To formulate what we expect, we use the {\it logarithmico-exponential functions} of Hardy \cite{Hardy1912}, or {\it LE-functions}: the real functions, defined for all large $x$, which are obtained from $x$ and real constants by finitely many additions, multiplications, divisions, real powers, exponentials and logarithms. By \Cref{thm:Explicit saddle point}(iii), for $k \geq 2$ the saddle point is governed by the inverse function of $\mathcal{H}_k$, and we expect that this function cannot be approximated by LE-functions up to an error of a negative power of $n$. The motivation for the following conjecture is discussed in \Cref{sec:Motivation}.

\begin{conjecture}
\label{thm:No closed form}
Let $P$ be a real polynomial of degree $d \geq 1$ with positive leading coefficient, let $\lambda > 0$, and 
  let $w_n$ be the largest real solution of the equation $e^{w} P(w) = \lambda n$. Then there is no LE-function $L$ such that
\[
w_n - L(n) = O(n^{-c}), \qquad n \to \infty,
\]
  for some constant $c > 0$. 
\end{conjecture}

\noindent
For $P(w) = w$ and $\lambda = 1$ we have $w_n = W(n)$, 
and the conjecture says that the Lambert function cannot be approximated by an LE-function with an error $O(n^{-c})$. For $k \geq 2$, the generalized Lambert function $w_n$ of \Cref{thm:Explicit saddle point}(iii) is also of this form with $P(w) = \prod_{i = 1}^{k - 1} (w - w_i)$, $d = k - 1$ and $\lambda = (1 + \alpha)^{k}/C$; for $k = 2$ it reduces to the Lambert function, since $w_n = w_1 + W(\lambda e^{-w_1} n)$. In \Cref{sec:Motivation} we show that the conjecture has the following consequence: if $k \geq 2$, then for every LE-function $L$ and every $c > 0$ we have $\log \rho_n^{-1} - L(n) \neq O(n^{-c})$, so that, by \Cref{thm:Accuracy of saddle point}(i), \cref{eq:Perturbed saddle point} fails for $\tilde{\rho} = e^{-L(n)}$. In other words, no explicit formula of this kind gives the saddle point with the precision needed for an asymptotic formula for $p_f(n)$, while the generalized Lambert function followed by finitely many Newton steps does (\Cref{thm:Dominant saddle point}(iii)).

\begin{remark}
\label{thm:Remark on explicit forms}
\Cref{thm:No closed form} is not a consequence of the fact that $W$ is not an elementary (Liouville) function \cite{BCDJ}. That fact concerns exact identities, while the conjecture concerns approximation up to an error of a negative power of $n$, and a function which is not elementary may well be approximated by elementary functions to such a precision; for instance, $\log n + W(n)/n^2$ is not elementary, but it differs from $\log n$ by $O(n^{-3/2})$. Finite truncations of the classical expansion $W(n) = \log n - \log \log n + \frac{\log \log n}{\log n} + \cdots$, and also a fixed number of Newton steps for $w e^{w} = n$ started at $\log n$, are LE-functions whose error is a negative power of $\log n$; an error $O(n^{-c})$ is reached only after about $\log_2 \log n$ Newton steps (\Cref{sec:Motivation}). For $k \geq 3$, the generalized Lambert function $w_n$ is easily computed numerically by Newton's method. In \Cref{thm:Dominant saddle point}(i) and (ii), the bound \cref{eq:Exponential bound} is required on the line $\Re(s) = \sigma_0$, with $\sigma_0 < \alpha/2$ or $\sigma_0 < 2\alpha/3$. Part (iii) holds with every $\sigma_0$ from (A2), but a smaller $\sigma_0$ means fewer Newton steps. When $D(s)$ contains a factor $\zeta(as - b)^{-1}$, as for $\mu(n)^2$, $\varphi(n)$ and $2^{\omega(n)}$, this may require the Riemann hypothesis, under which $1/\zeta(\sigma + it)$ grows more slowly than every positive power of $|t|$ for every fixed $\sigma > 1/2$ \cite[Chapter~XIV]{Titchmarsh}.
\end{remark}

\begin{remark}
The algebraic theory of elementary functions, which shows that $W$ is not elementary \cite{BCDJ}, decides whether a function is equal to an LE-function, but not whether it can be approximated by one with an error $O(n^{-c})$. Questions about approximation by LE-functions go back to Hardy \cite{Hardy1912}, who suggested that the inverse function of $(\log x)(\log \log x)$ is not asymptotic to any LE-function. This was confirmed by van den Dries, Macintyre and Marker \cite{DriesMacintyreMarker}, and Shackell \cite{Shackell1993} proved the corresponding statement for $(\log \log x)(\log \log \log x)$. These results concern asymptotic equivalence, while for the Lambert function the relation $W(x) \sim \log x$ is trivial, and \Cref{thm:No closed form} concerns approximation with an error $O(x^{-c})$. The proofs of these results use expansions of LE-functions into generalized power series, and such expansions may also be useful for \Cref{thm:No closed form}. 
\end{remark}


\section{Analysis of the generating function}
\label{sec:Generating function}

In this section we prepare the two tools on which the saddle-point method rests: the Cauchy integral for $p_f(n)$ (\Cref{thm:Ordinary Cauchy formula}), and the expansion of $\Psi = \log P$ near $\rho = 0$ into the residues of a Mellin integral (\Cref{thm:Closed form of psi}). This expansion is the source of all terms in the formulas of \Cref{sec:Main result} and \Cref{sec:Further results}. These facts are mostly standard.

\begin{lemma}
\label{thm:Ordinary Cauchy formula}
Let $n \geq 0$, $\rho > 0$. Then $p_f(n)$ has the following integral formula:  
\begin{equation}
\label{eq:Ordinary Cauchy formula}
p_f(n) = \frac{e^{\rho n}}{2 \pi} \int_{-\pi}^{\pi} P(\rho + it) e^{itn} dt.
\end{equation}
\end{lemma}

\begin{proof}
By Cauchy's integral formula, $p_f(n)$ admits the representation
\[
p_f(n) = \frac{1}{2 \pi i} \int_{\rho - i \pi}^{\rho + i \pi} e^{\tau n} P(\tau)\, d\tau, \quad n \geq 0, \rho > 0.
\]
Writing $\tau = \rho + it$ and changing variables, we obtain \cref{eq:Ordinary Cauchy formula}. 
\end{proof}

\subsection{Mellin inversion formula}
Let 
\begin{equation}
\label{eq:Psi(s)}
\Psi(\tau) := \log P(\tau) = -\sum_{n = 1}^{\infty} f(n) \log(1 - e^{-\tau n}).
\end{equation}
With this notation, by using \Cref{thm:Ordinary Cauchy formula}, $p_f(n)$ takes the following form 
\begin{equation}
\label{eq:Exponential Cauchy formula}
p_f(n) = \frac{e^{\rho n}}{2 \pi} \int_{-\pi}^{\pi} \exp(\Psi(\rho + it)) e^{itn}dt.
\end{equation}

\begin{proposition}
\label{thm:Mellin inversion formula}
$\Psi(\tau)$ has the Mellin inversion formula: 
\begin{equation}
\label{eq:Mellin inversion formula}
\Psi(\tau) = \frac{1}{2 \pi i} \int_{c - i \infty}^{c + i \infty} \Gamma(z) \zeta(z + 1) D(z) \tau^{-z} dz, \quad \Re(\tau) > 0, c > \alpha.
\end{equation}
\end{proposition}
\begin{proof}
Using $\displaystyle -\log(1 - z) = \sum_{r = 1}^{\infty} \frac{z^r}{r}$, for $|z| < 1$, expand \cref{eq:Psi(s)}:
\[
\Psi(\tau) = \sum_{n = 1}^{\infty} \sum_{r = 1}^{\infty}f(n) \frac{e^{-nr\tau}}{r}
= \sum_{n = 1}^{\infty} \frac{e^{-n\tau}}{n} \sum_{m \mid n} m f(m)
\]
Then 
\begin{equation}
\label{eq:g(n) form of Psi}
\Psi(\tau) = \sum_{n = 1}^{\infty} \frac{g(n)}{n} e^{-n\tau}, \qquad g(n) := \sum_{m \mid n} m f(m)
\end{equation}
For $\Re(z) > \alpha$, using the Dirichlet convolution we have 
\begin{equation}
\label{eq:Dirichlet convolution}
\sum_{n = 1}^{\infty} \frac{g(n)}{n^{z + 1}} = \zeta(z + 1) \sum_{m = 1}^{\infty} \frac{f(m)}{m^z} = \zeta(z + 1) D(z)
\end{equation}
Using the Mellin inversion theorem, for $c > \alpha$ we have 
\[
e^{-\tau n} = \frac{1}{2 \pi i} \int_{c - i \infty}^{c + i \infty} \Gamma(z) (\tau n)^{-z} dz
\]
Inserting it into \cref{eq:g(n) form of Psi} gives 
\[
\Psi(\tau) = \sum_{n = 1}^{\infty} \frac{g(n)}{n} \cdot \frac{1}{2 \pi i} \int_{c - i \infty}^{c + i \infty} \Gamma(z) (\tau n)^{-z} dz
\]
Since $c > \alpha$, absolute convergence allows us to apply Fubini's theorem and swap the integral sign with the sum, together with \cref{eq:Dirichlet convolution}, we obtain \cref{eq:Mellin inversion formula}.
\end{proof}

\noindent
Using this form of $\Psi(\tau)$, we show that the saddle-point equation has a unique solution. 

\begin{lemma}
\label{thm:Rho existence}
For each $n \in \mathbb{N}$, the equation $\Psi'(\rho) + n = 0$ has a unique solution $\rho = \rho_n > 0$, and $\rho_n \to 0^{+}$ as $n \to \infty$.
\end{lemma}

\begin{proof}
Differentiating \cref{eq:g(n) form of Psi} we get
\[
\Psi'(\tau) = -\sum_{n = 1}^{\infty} g(n) e^{-n\tau}
\]
Since $f(n) \ge 0$, 
the function $g(n)$ is always non-negative. Moreover, if $f(m) > 0$, then $g(jm) \geq m f(m)$ for every $j \in \mathbb{N}$, so that $\sum_{n} g(n) = \infty$. Hence, $-\Psi'(\rho)$ is continuous and strictly decreasing from $\infty$ to $0$ on $(0,\infty)$, so $-\Psi'(\rho) = n$ has a unique solution $\rho_n$, and $\rho_n \to 0^{+}$ as $n \to \infty$.
\end{proof}

\subsection{Meromorphic continuation}
To expand $\Psi(\rho)$ as $\rho \to 0^{+}$, we move the line of integration in \cref{eq:Mellin inversion formula} from $\Re(z) = c$ to $\Re(z) = -\varepsilon$. On the way we pass the poles of the integrand, and their residues give the expansion. Let 
\begin{equation}
\label{eq:M function}
M(z,\tau) := \Gamma(z) \zeta(z + 1) D(z) \tau^{-z}
\end{equation}
The functions $\Gamma(z)$ and $\zeta(z+1)$ admit meromorphic continuation to the entire complex plane, while $D(z)$ admits a meromorphic continuation to the half-plane $\Re(z) \geq -\varepsilon$. It follows that, for each fixed $\tau$ with $\Re(\tau) > 0$, the function $M(z,\tau)$ admits a meromorphic continuation to the half-plane $\Re(z) \geq -\varepsilon$. For $c > \alpha, T > 0$ define 
\[
(c) := \{c + it : t \in [-T,T]\}, \quad (-\varepsilon) := \{-\varepsilon + it: t \in [-T,T]\}
\]
\[
L_1 := \{\sigma + iT : \sigma \in [-\varepsilon,c]\}, \quad L_2 := \{\sigma - iT : \sigma \in [-\varepsilon,c]\}
\]
We can connect them to make a closed contour $\mathcal{C}$:  
\[
\mathcal{C} :=  (c) \oplus L_1^- \oplus (-\varepsilon)^- \oplus L_2
\]
where $L_1^-, (-\varepsilon)^-$ are taken in the negative direction so that we have 
\[
\oint_{\mathcal{C}} M(z,\tau) dz = \left(\int_{(c)} - \int_{L_1} - \int_{(-\varepsilon)} + \int_{L_2}\right) M(z,\tau) dz  = 2 \pi i\sum_{\operatorname{Int}(\mathcal{C})} \operatorname{Res}
\]
where the sum of residues is taken in the vertical strip $[-\varepsilon, c]$. To let $T \to \infty$, we need the integrals over the horizontal segments to vanish; this is where the heights $T_j$ of (A2) are used.

\begin{lemma}
\label{thm:0 integrals}
Let $\tau > 0$ and $T = T_j$, where $(T_j)$ is the sequence of condition (A2). As $j \to \infty$, we have 
\[
\left| \int_{L_1} M(z,\tau) dz \right| \to 0, \quad \left| \int_{L_2} M(z,\tau) dz \right| \to 0
\]
\end{lemma}

\begin{proof}
On $L_1$, the Gamma function has the following asymptotics: 
\[
|\Gamma(\sigma + iT)| \sim \sqrt{2 \pi} |T|^{\sigma - 1/2} e^{-\frac{\pi}{2}|T|} \quad (|T| \to \infty, \sigma \text{ bounded})
\]
So uniformly for $\sigma \in [-\varepsilon, c]$ we have 
\[
|\Gamma(\sigma + iT)|\ll |T|^{c - 1/2} e^{-\frac{\pi}{2}|T|}
\]
For the Riemann zeta function we have 
\[
|\zeta(\sigma + 1 + iT)| \ll |T|^b
\]
for some constant $b > 0$. For the exponential factor, since $\tau > 0$,
\[
|\tau^{-z}| = \tau^{-\sigma}
\]
which is bounded for $\sigma \in [-\varepsilon,c]$. By our assumption \cref{eq:Horizontal bound},
\[
D(\sigma \pm iT_j) \ll e^{a T_j}
\]
uniformly for $\sigma \in [-\varepsilon, c]$, with $a < \pi/2$. Combining all we get  
\[\left| \int_{L_1} M(z,\tau) dz \right| \leq (c + \varepsilon) \sup_{\sigma \in [-\varepsilon,c]} |M(\sigma + iT,\tau)| \ll |T|^{b + c - 1/2} e^{\left(a -\frac{\pi}{2}\right)|T|} \to 0, \quad |T| \to \infty. \]
The same bound holds on $L_2$, so the integral over $L_2$ also tends to zero.
\end{proof}

\noindent
Letting $T = T_j \to \infty$, we obtain $\Psi(\rho)$ as a sum of residues.

\begin{proposition}    
\label{thm:Residues}
For $\rho > 0$, the function $\Psi(\rho)$ is given by the residues of $M(z, \rho)$ at the poles enclosed by the contour $\mathcal{C}$:
\begin{equation}
\label{eq:Residues}
\Psi(\rho) = \sum_{\operatorname{Int}(\mathcal{C})} \operatorname{Res} + O(\rho^{\varepsilon})
\end{equation}
\end{proposition}

\begin{proof}
Let $T = T_j$. Since the segments $L_1$ and $L_2$ contain no poles, letting $j \to \infty$ in \Cref{thm:0 integrals} gives
\[
\oint_{\mathcal{C}} M(z,\tau) dz  = \int_{(c)} M(z,\tau) dz - \int_{(-\varepsilon)} M(z,\tau) dz = 2 \pi i \sum_{\operatorname{Int}(\mathcal{C})} \operatorname{Res}
\]
But by \Cref{thm:Mellin inversion formula} we have 
\[
\Psi(\tau) = \frac{1}{2 \pi i} \int_{(c)} M(z,\tau) dz
\]
So, 
\[
\Psi(\tau) = \sum_{\operatorname{Int}(\mathcal{C})} \operatorname{Res} + \frac{1}{2 \pi i}\int_{(-\varepsilon)} M(z,\tau) dz 
\]
Here the integrals are taken over the whole lines $\Re(z) = c$ and $\Re(z) = -\varepsilon$, and the sum over the poles is the limit described in \Cref{thm:Infinite B}. Now let $\tau = \rho > 0$. On the line $\Re(z) = -\varepsilon$ we have $|\rho^{-z}| = \rho^{\varepsilon}$. Moreover, $|\Gamma(-\varepsilon + it)| \ll (1 + |t|)^{-\varepsilon - 1/2} e^{-\pi |t|/2}$ and $\zeta(1 - \varepsilon + it) \ll (1 + |t|)^{b}$, as in the proof of \Cref{thm:0 integrals}, and $D(-\varepsilon + it) \ll e^{a |t|}$ by \cref{eq:Exponential bound}, with $a < \pi/2$. Hence
\[
\int_{(-\varepsilon)} M(z,\rho)\, dz \ll \rho^{\varepsilon} \int_{-\infty}^{\infty} \left| \Gamma(-\varepsilon + it) \zeta(1 - \varepsilon + it) D(-\varepsilon + it) \right| dt \ll \rho^{\varepsilon}
\]
which gives \cref{eq:Residues}.
\end{proof}

\subsection{Residues}
By \Cref{thm:Residues}, $\Psi(\rho)$ is, up to $O(\rho^{\varepsilon})$, the sum of the residues of $M(z, \rho)$ inside the contour $\mathcal{C}$, which contains every pole of $D(z)$ in $\Re(z) > -\varepsilon$. Since $0 < \varepsilon < 1$, the functions $\Gamma(z)$ and $\zeta(z + 1)$ are holomorphic in this region except at $z = 0$, where $\Gamma(z) \zeta(z + 1)$ has a pole of order $2$, which combines with the pole of $D(z)$ of order $h_0$. So we treat the point $z = 0$ and the poles $p \neq 0$ differently, and we will use the notation $M(z,\rho)$ for calculations.
\begin{lemma}
\label{thm:Res at 0}
Let $d_j$ be the Laurent coefficients of $D$ at $z = 0$, so that, with $\phi(0, \nu)$ as in \Cref{sec:Main result1},
\[
\Gamma(z) \zeta(z + 1) = \sum_{\nu \geq -2} \phi(0, \nu) z^{\nu}, \qquad D(z) = \sum_{j \geq -h_0} d_j z^{j}
\]
Then the residue of $M(z,\rho)$ at $z = 0$ is
\begin{equation}
\label{eq:Res at 0}
\operatorname{Res}_{z = 0} M(z,\rho) = \sum_{r = 0}^{h_0 + 1} F(0, r) \frac{(\log \rho^{-1})^r}{r!}, \qquad F(0, r) := \sum_{\nu + j = -1 - r} \phi(0, \nu)\, d_j
\end{equation}
In particular, if $h_0 = 0$, then $\operatorname{Res}_{z = 0} M(z,\rho) = D(0) \log \rho^{-1} + D'(0)$.
\end{lemma}

\noindent
\begin{proof}
First, write the Laurent series expansions around $z = 0$, where $\gamma$ is Euler's constant
\[
\Gamma(z) = \frac{1}{z} - \gamma + O(|z|), \quad \zeta(z + 1) = \frac{1}{z} + \gamma  + O(|z|)
\]
which imply 
\[
\Gamma(z) \zeta(z + 1) = \frac{1}{z^2} + O(1)
\]
so that $\phi(0, -2) = 1$ and $\phi(0, -1) = 0$. Next, we expand 
\[
\rho^{-z}= \exp(-z \log \rho) = \exp(z \log \rho^{-1}) = \sum_{r = 0}^{\infty} \frac{(\log \rho^{-1})^r}{r!} z^r
\]
Multiplying the three expansions, we obtain
\[
M(z,\rho) = \left( \sum_{\nu \geq -2} \phi(0, \nu) z^{\nu} \right) \left( \sum_{r = 0}^{\infty} \frac{(\log \rho^{-1})^r}{r!} z^r \right) \left( \sum_{j \geq -h_0} d_j z^{j} \right)
\]
We need to extract the coefficient at $z^{-1}$, hence we require $\nu + r + j = -1$. Since $\nu \geq -2$ and $j \geq -h_0$, this gives $r \leq h_0 + 1$, and \cref{eq:Res at 0} follows. If $h_0 = 0$, then $d_0 = D(0)$ and $d_1 = D'(0)$, and
\[
F(0, 0) = \phi(0, -2) d_1 + \phi(0, -1) d_0 = D'(0), \qquad F(0, 1) = \phi(0, -2) d_0 = D(0)
\]
\end{proof}

\begin{lemma}
\label{thm:Res at p}
The residue of $M(z,\rho)$ at a pole $z = p \neq 0$ of $D(s)$, real or complex, with order $h$ and principal part 
\cref{eq:Principal part} is 
\begin{equation}
\label{eq:Res at p}
\operatorname{Res}_{z = p} M(z,\rho) = \rho^{-p} \sum_{r = 0}^{h - 1} F(p,r) \frac{(\log \rho^{-1})^r}{r!}
\end{equation}
where $F(p,r)$ is defined in \Cref{sec:Main result}. 
\end{lemma}

\begin{proof}
Since $p \neq 0$ and $\Re(p) > -1$, both $\Gamma(p + s)$ and $\zeta(p + 1 + s)$ are holomorphic at $s = 0$, hence they have Taylor expansions: 
\[
\Gamma(p + s) = \sum_{\nu = 0}^{\infty} \frac{\Gamma^{(\nu)}(p)}{\nu!}s^{\nu} \quad \zeta(p + 1 + s) = \sum_{\nu = 0}^{\infty} \frac{\zeta^{(\nu)}(p + 1)}{\nu!} s^{\nu}
\]
By the Cauchy product, with $\phi(p, \nu)$ as in \Cref{sec:Main result1},
\[
\Gamma(p + s) \zeta(p + 1 + s) = \sum_{\nu = 0}^{\infty} \phi(p,\nu) s^{\nu}
\]
Now expand $\rho^{-(p + s)}$
\[
\rho^{-(p + s)} = \rho^{-p} \sum_{r = 0}^{\infty} \frac{(\log \rho^{-1})^{r}}{r!} s^r
\]
By our assumption \cref{eq:Principal part}, if we multiply all the factors 
\[
M(p + s,\rho) = \rho^{-p} \left(\sum_{\nu = 0}^{\infty} \phi(p,\nu) s^{\nu}\right)\left(\sum_{r = 0}^{\infty} \frac{(\log \rho^{-1})^r}{r!} s^r\right)\left( \sum_{j = 1}^{h} A(p, j) s^{-j} + H_p(p + s)\right)
\]
We need to extract the coefficient at $s^{-1}$, hence we require 
$
s^{\nu} s^r s^{-j} = s^{-1} \iff \nu + r - j = -1. 
$
Therefore,
\[
\operatorname{Res}_{z = p} M(z,\rho) = \rho^{-p} \sum_{\nu + r - j = -1} \phi(p,\nu) \cdot A(p,j) \frac{(\log \rho^{-1})^r}{r!} 
\]
For fixed $r$, the sum over $\nu$ and $j$ is
\[
\sum_{\substack{j - \nu = r + 1\\ 1 \leq j \leq h,\; \nu\ge0}} \phi(p,\nu)\cdot A(p,j) = \sum_{r + 1 \leq j \leq h} \phi(p, j - r - 1) \cdot A(p,j) = F(p, r)
\]
which gives \cref{eq:Res at p}.
\end{proof}

\noindent
Now, we are ready to give the closed form for $\Psi(\rho)$.
\begin{corollary}    
\label{thm:Closed form of psi}
$\Psi(\rho)$ has the form 
\begin{equation}
\label{eq:Closed form of psi}
\begin{aligned}
\Psi(\rho) = {} & \sum_{\substack{1 \leq i \leq m \\ \alpha_i \neq 0}} \rho^{-\alpha_i} \sum_{r = 0}^{k_i - 1} F(\alpha_i,r) \frac{(\log \rho^{-1})^r}{r!} + \sum_{\beta \in \mathcal{B}} \rho^{-\beta} \sum_{r = 0}^{\gamma_\beta - 1} F(\beta,r) \frac{(\log \rho^{-1})^r}{r!} \\
& + \sum_{r = 0}^{h_0 + 1} F(0, r) \frac{(\log \rho^{-1})^r}{r!} + O(\rho^\varepsilon)
\end{aligned}
\end{equation}
where $\gamma_\beta$ is the order of the pole $\beta \in \mathcal{B}$.
\end{corollary}

\begin{proof}
By putting the values of residues by \Cref{thm:Res at 0} and \Cref{thm:Res at p} for every pole into \Cref{thm:Residues}, we have \cref{eq:Closed form of psi}.
\end{proof}



\begin{remark}
\label{thm:Infinite B}
    If $\mathcal{B}$ is infinite, then $|\Im(\beta)| \to \infty$ along $\mathcal{B}$. Indeed, $D(s)$ is meromorphic on a neighbourhood of the closed half-plane $\Re(s) \geq -\varepsilon$, so its poles are isolated, and they all satisfy $\Re(s) \leq \alpha$ by \Cref{thm:Landau}. Hence every compact rectangle $\{-\varepsilon \leq \Re(s) \leq \alpha, \ |\Im(s)| \leq T\}$ contains only finitely many poles. Accordingly, the sums over the poles in \cref{eq:Residues} and over $\mathcal{B}$ in \cref{eq:Closed form of psi} are understood as
\[
\sum_{\beta \in \mathcal{B}} = \lim_{j \to \infty} \sum_{\substack{\beta \in \mathcal{B} \\ |\Im(\beta)| < T_j}} 
\]
This limit exists: for $T = T_j$, the residue theorem expresses the finite sum over the poles inside $\mathcal{C}$ through the integrals over $(c)$, $(-\varepsilon)$, $L_1$ and $L_2$. As $j \to \infty$, the integrals over $(c)$ and $(-\varepsilon)$ converge absolutely, and those over $L_1$ and $L_2$ tend to $0$ by \Cref{thm:0 integrals}. Condition (A2) gives no bound on the coefficients $A(\beta, j)$, so we do not claim that the series over $\mathcal{B}$ converges absolutely.
\end{remark}
\section{Analysis of contributions}
\label{sec:Contributions}
Throughout this section we assume the conditions (A1)-(A3), and we prove \Cref{thm:Main theorem}. We estimate the integral \cref{eq:Exponential Cauchy formula} with $\rho = \rho_n$ from \Cref{thm:Rho existence}. Since $\rho_n \to 0^{+}$, we may assume $0 < \rho \leq 1/2$. The integral is split into the major arcs $|t| \leq \rho^{\chi}$, where $P(\rho + it)$ is close to a Gaussian, and the minor arcs $\rho^{\chi} \leq |t| \leq \pi$, which are negligible by (A3). 
The exponent of the major arcs is denoted by $\chi$, 
the order of a derivative by $\ell$, 
and the number of terms in the expansion of the exponential by $J_0$. 

Both parts use two consequences of the non-negativity of $f$, which we prove first: $\alpha$ is the abscissa of convergence of $D(s)$, and the leading coefficient $A(\alpha, k)$ is positive.

\begin{proposition}
\label{thm:Landau}
The abscissa of convergence of $D(s)$ is $\alpha$, and the series converges absolutely for $\Re(s) > \alpha$. Consequently, every complex pole $\beta \in \mathcal{B}$ satisfies $\Re(\beta) \leq \alpha$, and a pole of $D(s)$ on the line $\Re(s) = \alpha$ has order at most $k$.
\end{proposition}

\begin{proof}
Let $\sigma_c$ be the abscissa of convergence of $D(s)$. Since $f(n) \geq 0$, convergence and absolute convergence coincide. In $\Re(s) > \sigma_c$ the series defines a holomorphic function, which coincides with the meromorphic continuation of $D(s)$. Hence $D(s)$ has no poles in $\Re(s) > \sigma_c$, and since $\alpha$ is a pole, $\sigma_c \geq \alpha > 0$.

By Landau's theorem (see e.g. \cite[Theorem~1.7]{MontgomeryVaughan}), since $f(n) \geq 0$, the real point $s = \sigma_c$ is a singularity of $D(s)$. In the half-plane $\Re(s) \geq -\varepsilon$ the only singularities of $D(s)$ are its poles, and the real ones are $\alpha_1 > \cdots > \alpha_m$. Thus $\sigma_c \in \{\alpha_1, \ldots, \alpha_m\}$, and $\sigma_c \geq \alpha = \alpha_1$ gives $\sigma_c = \alpha$. Since $D(s)$ is holomorphic in $\Re(s) > \alpha$, every $\beta \in \mathcal{B}$ satisfies $\Re(\beta) \leq \alpha$.

Finally, let $\beta = \alpha + it_0$ be a pole of order $h$. By \cref{eq:Principal part}, as $\sigma \to \alpha^{+}$ we have $|D(\sigma + it_0)| \asymp (\sigma - \alpha)^{-h}$ and $D(\sigma) \asymp (\sigma - \alpha)^{-k}$. On the other hand,
\[
|D(\sigma + it_0)| \leq \sum_{n = 1}^{\infty} \frac{f(n)}{n^{\sigma}} = D(\sigma),
\]
and hence $h \leq k$.
\end{proof}

\begin{remark}\label{rem:ex}
Landau's theorem gives $\Re(\beta) \leq \alpha$, but it cannot exclude poles on the line $\Re(s) = \alpha$, which is why the assumption \cref{eq:Complex poles to the left} in (A1) is needed. For example, $f(n) = 1 + \cos(t_0 \log n) \geq 0$ with $t_0 \neq 0$ gives
\[
D(s) = \zeta(s) + \tfrac{1}{2} \zeta(s - it_0) + \tfrac{1}{2} \zeta(s + it_0)
\]
which satisfies (A3) and the growth bounds in (A2), but violates \cref{eq:Complex poles to the left} and hence also the condition on the poles in (A2), with $\alpha = 1$, $k = 1$, but has the poles $1 \pm it_0$. By \Cref{thm:Res at p}, these contribute $\Re\big( \Gamma(1 + it_0) \zeta(2 + it_0) \rho^{-1 - it_0} \big)$ to $\Psi(\rho)$, a term of the same order $\rho^{-1}$ as the main term, which oscillates as $\rho \to 0^{+}$. For such $f$, \Cref{thm:Main theorem} does not apply.
\end{remark}

\begin{proposition}
\label{thm:Positive leading coefficient}
$F(\alpha, k - 1) = \Gamma(\alpha) \zeta(\alpha + 1) A(\alpha, k) > 0$. In particular $C > 0$.
\end{proposition}

\begin{proof}
In the formula for $F(p, r)$ in \Cref{sec:Main result} with $p = \alpha$, $r = k - 1$, only $j = k$ contributes, and $\phi(\alpha, 0) = \Gamma(\alpha) \zeta(\alpha + 1) > 0$. For real $\sigma > \alpha$ the series converges by \Cref{thm:Landau}, so $D(\sigma) = \sum_{n} f(n) n^{-\sigma} \geq 0$. Hence by \cref{eq:Principal part}, $A(\alpha, k) = \lim_{\sigma \to \alpha^{+}} (\sigma - \alpha)^{k} D(\sigma) \geq 0$, and $A(\alpha, k) \neq 0$ since the pole has order $k$. Finally $C = \alpha(1 + \alpha) F(\alpha, k - 1)/(k - 1)! > 0$.
\end{proof}

\noindent
We now split the integral into the major and minor arcs: for any $\chi \geq 1$,
\begin{equation}
\label{eq:Splitted integral}
p_f(n) = \frac{e^{\rho n}}{2 \pi} \left( I + \int_{\rho^{\chi} \leq |t| \leq \pi} e^{itn}\exp(\Psi(\rho + it)) dt \right), \quad I := \int_{|t| \leq \rho^{\chi}} e^{itn}\exp(\Psi(\rho + it)) dt
\end{equation}
For the minor arcs we follow the argument of \cite[Lemma~3.1]{Tenenbaum2020}, see also \cite[Lemma~3.1]{Bridges2023}.

\begin{lemma}
\label{thm:Minor part}
Let $1 < \chi < 1 + \frac{\alpha}{2}$ and $\alpha/2 < M_1 < M_0$. Then the integral over the minor arcs satisfies
\begin{equation}
\label{eq:Minor part}
\int_{\rho^{\chi} \leq |t| \leq \pi} e^{itn}\exp(\Psi(\rho + it)) dt  \ll_{M_1} \rho^{M_1 + 1} P(\rho).
\end{equation}
\end{lemma}

\noindent
\begin{proof}
Since $|P(\rho - it)| = |P(\rho + it)|$, it suffices to show that
\[
\int_{\rho^{\chi}}^{\pi} |P(\rho + it)|\, dt \ll \rho^{M_1 + 1} P(\rho)
\]
As $f(m)$ is real, we have $|(1 - e^{-m\tau})^{-f(m)}| = |1 - e^{-m\tau}|^{-f(m)}$, and
\[
|1 - e^{-m(\rho + it)}|^2 = (1 - e^{-m\rho})^2 + 4 e^{-m\rho} \sin^2(mt/2)
\]
Hence, by \cref{eq:Generating function},
\begin{equation}
\label{eq:Ratio product}
\frac{|P(\rho + it)|}{P(\rho)} = \prod_{m = 1}^{\infty} \left( 1 + \frac{4 \sin^2(mt/2)}{e^{m\rho} (1 - e^{-m\rho})^2} \right)^{-f(m)/2} 
\end{equation}
Every factor in \cref{eq:Ratio product} is at most $1$, so any of them may be discarded in an upper bound. We put $x = t/2\pi$. Since $|\sin(\pi y)| \geq 2\|y\|$ for real $y$, and $e^{m\rho}(1 - e^{-m\rho})^2 \leq e\, (m\rho)^2$ for $m\rho \leq 1$, the $m$-th factor satisfies
\begin{equation}
\label{eq:Single factor}
\left( 1 + \frac{4 \sin^2(mt/2)}{e^{m\rho} (1 - e^{-m\rho})^2} \right)^{-f(m)/2} \leq \left( 1 + \frac{16 \|mx\|^2}{e\, m^2 \rho^2} \right)^{-f(m)/2}, \qquad m\rho \leq 1.
\end{equation}

We next bound the partial sums of $f$ from below. By \cref{eq:Complex poles to the left} and \Cref{thm:Landau}, the function $D(s) - \sum_{j = 1}^{k} A(\alpha, j)(s - \alpha)^{-j}$ is holomorphic in a neighbourhood of the closed half-plane $\Re(s) \geq \alpha$. Since $f(n) \geq 0$, Delange's Tauberian theorem \cite{Delange1954} yields
\[
\sum_{1 \leq m \leq x} f(m) \sim \frac{A(\alpha, k)}{\alpha (k - 1)!} x^{\alpha} (\log x)^{k - 1}, \qquad x \to \infty,
\]
where $A(\alpha, k) > 0$ by \Cref{thm:Positive leading coefficient}. Consequently, there exist $X_0 \geq 3$ and $c_0 > 0$ such that
\begin{equation}
\label{eq:Partial sums lower bound}
\sum_{1 \leq m \leq X} f(m) \geq c_0 X^{\alpha}, \qquad \sum_{X < m \leq 2X} f(m) \geq c_0 X^{\alpha}, \quad X \geq X_0.
\end{equation}
We fix $R_0 \geq X_0$ such that $3R_0$ is an integer and $c_0 R_0^{\alpha} \geq 2M_1 + 2$. By (A3), $\sum_{m \in \Lambda \setminus p\mathbb{N}} f(m) \geq M_0 > M_1$, so for every prime $p \leq 3R_0$ there is a finite set $S_p \subset \Lambda \setminus p\mathbb{N}$ with $\sum_{m \in S_p} f(m) \geq M_1$. Let $m^{*}$ denote the largest element of the union of these sets, and let $C_0 := \frac{9\sqrt{e}}{4} R_0 m^{*}$. Since $\rho = \rho_n \to 0^{+}$, we may assume that $C_0 \rho \leq 1$, $2\rho \leq 1/X_0$ and $\rho^{\chi} < 2\pi\rho$; note that $C_0 \rho \leq 1$ implies $\rho \leq 1/m^{*}$.

Consider first the range $\rho^{\chi} \leq t \leq 2\pi\rho$, in which $0 < x \leq \rho$. For $m \leq 1/(2\rho)$ we have $mx \leq 1/2$, hence $\|mx\| = mx = mt/2\pi$, and therefore
\[
\frac{16 \|mx\|^2}{e\, m^2 \rho^2} = \frac{4 t^2}{\pi^2 e\, \rho^2} \geq c_1 \rho^{2\chi - 2}, \qquad c_1 := \frac{4}{\pi^2 e} < 1 
\]
We retain in \cref{eq:Ratio product} only the factors with $m \leq 1/(2\rho)$. By \cref{eq:Single factor} and the last display, each of them is at most $(1 + y)^{-f(m)/2}$ with $y := c_1 \rho^{2\chi - 2}$, since the factor decreases as the quantity in the last display increases. As $\chi > 1$ and $c_1 < 1$, we have $0 < y \leq 1$, hence $\log(1 + y) \geq y/2$, and \cref{eq:Partial sums lower bound} with $X = 1/(2\rho)$ gives
\[
\frac{|P(\rho + it)|}{P(\rho)} \leq \exp\left( -\frac{c_1}{4} \rho^{2\chi - 2} \sum_{m \leq 1/(2\rho)} f(m) \right) \leq \exp\left( -c_2 \rho^{-(2 + \alpha - 2\chi)} \right), \qquad c_2 := \frac{c_0 c_1}{2^{\alpha + 2}} 
\]
Since $\chi < 1 + \frac{\alpha}{2}$, the exponent $2 + \alpha - 2\chi$ is positive, and the right-hand side is $\ll \rho^{M_1 + 1}$.

Consider now the range $2\pi\rho \leq t \leq \pi$, in which $\rho \leq x \leq \frac{1}{2}$. By Dirichlet's approximation theorem there are coprime integers $a$ and $q$ with $1 \leq q \leq 3R_0$ such that
\[
x = \frac{a}{q} + \theta r, \qquad \theta \in \{-1, 1\}, \qquad 0 \leq r \leq \frac{1}{3qR_0} 
\]
If $q = 1$, then $a = 0$, because $0 < x \leq \frac{1}{2}$ and $r < \frac{1}{2}$; in that case $r = x \geq \rho$.

Suppose first that $r \geq \frac{2\rho}{3q}$. We set $X = \frac{1}{3rq}$, so that $X \geq R_0$, and consider the integers $m$ with $X < m \leq 2X$. For these $m$ we have $mr \in \left( \frac{1}{3q}, \frac{2}{3q} \right]$, while $\frac{ma}{q}$ is a multiple of $\frac{1}{q}$; hence $\|mx\| \geq \frac{1}{3q}$. Moreover, $q m \rho \leq \frac{2\rho}{3r} \leq q \leq 3R_0$, and in particular $m\rho \leq 1$. By \cref{eq:Single factor}, the $m$-th factor in \cref{eq:Ratio product} is therefore at most
\[
\left( 1 + \frac{16}{9e\, (qm\rho)^2} \right)^{-f(m)/2} 
\]
If in addition $r \leq \sqrt{\rho}$, then $\frac{16}{9e (qm\rho)^2} \geq \frac{16}{81 e R_0^2} \eqqcolon c_3$ and $X \geq \frac{1}{3q\sqrt{\rho}} \geq \frac{1}{9R_0 \sqrt{\rho}}$. Retaining only the factors with $X < m \leq 2X$ and applying \cref{eq:Partial sums lower bound}, we find
\[
\frac{|P(\rho + it)|}{P(\rho)} \leq (1 + c_3)^{-\frac{c_0}{2} X^{\alpha}} \leq \exp\left( -c_4 \rho^{-\alpha/2} \right) \ll \rho^{M_1 + 1}
\]
for some $c_4 > 0$ depending only on $R_0$. If instead $r > \sqrt{\rho}$, then $(qm\rho)^2 \leq \frac{4\rho^2}{9r^2} < \frac{4\rho}{9}$, so that $\frac{16}{9e(qm\rho)^2} > \frac{4}{e\rho} \geq \frac{1}{\rho}$, and each of these factors is at most $(1 + \rho^{-1})^{-f(m)/2} \leq \rho^{f(m)/2}$. By \cref{eq:Partial sums lower bound} and the choice of $R_0$, it follows that
\[
\frac{|P(\rho + it)|}{P(\rho)} \leq \rho^{\frac{1}{2} \sum_{X < m \leq 2X} f(m)} \leq \rho^{\frac{c_0}{2} R_0^{\alpha}} \leq \rho^{M_1 + 1} 
\]

Suppose finally that $r < \frac{2\rho}{3q}$. Then $q \geq 2$, since $q = 1$ forces $r \geq \rho$. Let $p$ be a prime divisor of $q$, so that $p \leq 3R_0$, and let $m \in S_p$. As $p \nmid m$, we have $q \nmid m$, and since $(a, q) = 1$ this gives $\|ma/q\| \geq 1/q$. Since $m\rho \leq 1$, we also have $mr < \frac{2}{3q}$, and hence $\|mx\| \geq \frac{1}{q} - \frac{2}{3q} = \frac{1}{3q}$. Retaining only the factors with $m \in S_p$, we deduce from \cref{eq:Single factor} that
\[
\frac{|P(\rho + it)|}{P(\rho)} \leq \prod_{m \in S_p} \left( \frac{16}{9e\, q^2 m^2 \rho^2} \right)^{-f(m)/2} = \prod_{m \in S_p} \left( \frac{3\sqrt{e}\, q m \rho}{4} \right)^{f(m)} \leq (C_0 \rho)^{M_1} 
\]
The values of $t \in [2\pi\rho, \pi]$ for which this situation occurs lie in the union of the intervals $\left| \frac{t}{2\pi} - \frac{a}{q} \right| < \frac{2\rho}{3q}$ with $2 \leq q \leq 3R_0$ and $1 \leq a \leq q$, whose total measure is at most
\[
\sum_{2 \leq q \leq 3R_0} q \cdot 2\pi \cdot \frac{4\rho}{3q} \ll_{R_0} \rho 
\]
Collecting the above estimates, we conclude that
\[
\int_{\rho^{\chi}}^{\pi} \frac{|P(\rho + it)|}{P(\rho)}\, dt \leq 2\pi\rho\, e^{-c_2 \rho^{-(2 + \alpha - 2\chi)}} + \pi \left( e^{-c_4 \rho^{-\alpha/2}} + \rho^{M_1 + 1} \right) + O(\rho) \cdot (C_0 \rho)^{M_1} \ll \rho^{M_1 + 1}
\]
which proves \cref{eq:Minor part}.
\end{proof}

\subsection{Derivatives of \texorpdfstring{$\Psi$}{Psi}}
On the major arcs we expand $\Psi(\rho + it)$ in powers of $t$, so we need the sizes of the derivatives $\Psi^{(\ell)}$, uniformly on vertical lines. Their signs come from the non-negativity of $f$, and their size from the pole $\alpha$, by a shift of the line of integration to $\Re(z) = \sigma_0$.

\begin{lemma}
\label{thm:Reality}
For $\rho > 0$ and $\ell \geq 0$, $\Psi^{(\ell)}(\rho)$ is real and $(-1)^{\ell} \Psi^{(\ell)}(\rho) > 0$; in particular $\Psi''(\rho) > 0$. Moreover, for every $t \in \mathbb{R}$,
\begin{equation}
\label{eq:Vertical domination}
\left| \Psi^{(\ell)}(\rho + it) \right| \leq (-1)^{\ell} \Psi^{(\ell)}(\rho). 
\end{equation}
\end{lemma}

\begin{proof}
Differentiating \cref{eq:g(n) form of Psi} term by term, as in the proof of \Cref{thm:Rho existence},
\[
\Psi^{(\ell)}(\tau) = (-1)^{\ell} \sum_{n = 1}^{\infty} g(n)\, n^{\ell - 1} e^{-n\tau}.
\]
Since $g(n) \geq 0$ and $g$ is not identically zero, the first claim follows at $\tau = \rho$. For $\tau = \rho + it$, use $|e^{-n(\rho + it)}| = e^{-n\rho}$ and the triangle inequality.
\end{proof}

\noindent
By induction, the polynomial $Q_\ell$ defined in \Cref{sec:Main result} has degree $k - 1$ and the leading coefficient
\begin{equation}
\label{eq:Leading coefficient}
\alpha(\alpha + 1) \cdots (\alpha + \ell - 1) \frac{F(\alpha, k - 1)}{(k - 1)!} = \frac{\Gamma(\alpha + \ell)}{\Gamma(\alpha + 2)}\, C. 
\end{equation}

\begin{lemma}
\label{thm:Derivatives of Psi}
Let $\sigma_0$ be as in (A2). For every $\ell \geq 0$, as $\rho \to 0^{+}$,
\begin{equation}
\label{eq:Derivative expansion}
(-1)^{\ell} \Psi^{(\ell)}(\rho) = \rho^{-\alpha - \ell} Q_\ell(\log \rho^{-1}) + O(\rho^{-\sigma_0 - \ell}).
\end{equation}
In particular $\Psi''(\rho) = \rho^{-\alpha - 2} Q_2(\log \rho^{-1}) + O(\rho^{-\sigma_0 - 2})$, and
\begin{equation}
\label{eq:Derivative asymptotics}
(-1)^{\ell} \Psi^{(\ell)}(\rho) = \frac{\Gamma(\alpha + \ell)}{\Gamma(\alpha + 2)}\, C\, \rho^{-\alpha - \ell} (\log \rho^{-1})^{k - 1} \left\{ 1 + O\left( \frac{1}{\log \rho^{-1}} \right) \right\}.
\end{equation}
Consequently, for every $\ell \geq 1$, uniformly for $0 < \rho \leq 1/2$ and $t \in \mathbb{R}$,
\begin{equation}
\label{eq:Derivative upper bound}
|\Psi^{(\ell)}(\rho + it)| \ll \rho^{-\alpha - \ell} (\log \rho^{-1})^{k - 1}, \qquad \Psi''(\rho) \gg \rho^{-\alpha - 2} (\log \rho^{-1})^{k - 1}.
\end{equation}
\end{lemma}

\begin{proof}
By (A2), every pole $p \neq \alpha$ of $D(s)$ satisfies $\Re(p) < \sigma_0$, where $0 < \sigma_0 < \alpha$. Then $z = \alpha$ is the only pole of $M(z, \rho)$ in the strip $\sigma_0 \leq \Re(z) \leq c$. Repeating the proof of \Cref{thm:Residues} with the line $(-\varepsilon)$ replaced by $(\sigma_0)$, where \Cref{thm:0 integrals} applies verbatim, and using \Cref{thm:Res at p} with $p = \alpha$, $h = k$, we get for every $\rho > 0$, 
\[
\Psi(\rho) = \rho^{-\alpha} Q_0(\log \rho^{-1}) + \frac{1}{2 \pi i} \int_{\sigma_0 - i\infty}^{\sigma_0 + i\infty} M(z, \rho)\, dz. 
\]
We differentiate $\ell$ times with respect to $\rho$. In the integral,
\[
\frac{\partial^{\ell}}{\partial \rho^{\ell}} M(z, \rho) = (-1)^{\ell} z(z + 1) \cdots (z + \ell - 1)\, \Gamma(z) \zeta(z + 1) D(z)\, \rho^{-z - \ell},
\]
and on $\Re(z) = \sigma_0$ the bounds of \Cref{thm:0 integrals} (with $\zeta(z + 1)$ bounded since $\Re(z + 1) > 1$) and $a < \pi/2$ show that this is absolutely integrable, locally uniformly in $\rho$. So we may differentiate under the integral sign, and the integral contributes $O(\rho^{-\sigma_0 - \ell})$. For the main term, since
\[
\frac{d}{d\rho} \left[ \rho^{-\alpha - \ell} Q_\ell(\log \rho^{-1}) \right] = -\rho^{-\alpha - \ell - 1} \left( (\alpha + \ell) Q_\ell + Q_\ell' \right)(\log \rho^{-1}) = -\rho^{-\alpha - \ell - 1} Q_{\ell + 1}(\log \rho^{-1}),
\]
induction gives $\frac{d^{\ell}}{d\rho^{\ell}} [\rho^{-\alpha} Q_0(\log \rho^{-1})] = (-1)^{\ell} \rho^{-\alpha - \ell} Q_\ell(\log \rho^{-1})$. This proves \cref{eq:Derivative expansion}, and \cref{eq:Derivative asymptotics} follows from \cref{eq:Leading coefficient} and $\rho^{\alpha - \sigma_0} \ll (\log \rho^{-1})^{-1}$.

By \cref{eq:Derivative asymptotics} and \Cref{thm:Positive leading coefficient}, $(-1)^{\ell} \Psi^{(\ell)}(\rho) \asymp \rho^{-\alpha - \ell} (\log \rho^{-1})^{k - 1}$ for $0 < \rho \leq \rho^{\circ}_\ell$ with some $\rho^{\circ}_\ell > 0$. On $[\rho^{\circ}_\ell, 1/2]$ both sides are continuous and positive (\Cref{thm:Reality}), so the relation holds on $(0, 1/2]$. Together with \cref{eq:Vertical domination} this gives \cref{eq:Derivative upper bound}.
\end{proof}

\noindent
With $\ell = 1$ at $\rho = \rho_n$, where $-\Psi'(\rho_n) = n$, the lemma gives the size of the saddle point.

\begin{corollary}
\label{thm:Size of rho}
As $n \to \infty$,
\[
\rho^{\alpha} (\log \rho^{-1})^{1 - k} \asymp \frac{1}{n \rho} \asymp n^{-\frac{\alpha}{1 + \alpha}} (\log n)^{-\frac{k - 1}{1 + \alpha}}. 
\]
\end{corollary}

\begin{proof}
By \cref{eq:Derivative asymptotics} with $\ell = 1$, $n = -\Psi'(\rho) \asymp \rho^{-\alpha - 1} (\log \rho^{-1})^{k - 1}$, which gives the first relation. Taking logarithms, $\log \rho^{-1} \sim \frac{\log n}{1 + \alpha}$, hence $\rho \asymp n^{-\frac{1}{1 + \alpha}} (\log n)^{\frac{k - 1}{1 + \alpha}}$ and the second relation follows.
\end{proof}

\subsection{The major arcs}

\noindent
On the major arcs, the integrand is close to the Gaussian 
$P(\rho) e^{-t^2 \Psi''(\rho)/2}$, whose factor $\Psi''(\rho)^{-1/2}$ is $\rho^{1 + \alpha/2}$ up to a power of $\log \rho^{-1}$, by \cref{eq:Derivative upper bound}. We choose the exponent $\chi$ so that the arc $|t| \leq \rho^{\chi}$ is much wider than this Gaussian, which needs $\chi < 1 + \frac{\alpha}{2}$, while on this arc the cubic term $t^3 \Psi'''(\rho)$, which is $\rho^{3\chi - 3 - \alpha}$ up to a power of $\log \rho^{-1}$, still tends to $0$, which needs $\chi > 1 + \frac{\alpha}{3}$.

Fix $1 + \frac{\alpha}{3} < \chi < 1 + \frac{\alpha}{2}$ and $N_0 > 0$, so that $0 < 3\chi - 3 - \alpha < \frac{\alpha}{2}$. Put
\[
\ell_0 := \left\lceil \frac{3(N_0 + \alpha)}{\alpha} \right\rceil, \qquad
J_0 := \left\lceil \frac{N_0}{3\chi - 3 - \alpha} \right\rceil, \qquad
E(t) := \sum_{3 \leq \ell \leq \ell_0 - 1} \frac{(it)^{\ell} \Psi^{(\ell)}(\rho)}{\ell!}
\]
and for $\ell \geq 3$ let
\begin{equation}
\label{eq:Lambda}
\lambda_{\ell}(\rho) := i^{\ell} \sum_{1 \leq q \leq J_0} \frac{1}{q!} \sum_{\substack{3 \leq m_1, \ldots, m_q \leq \ell_0 - 1 \\ m_1 + \cdots + m_q = \ell}} \prod_{1 \leq j \leq q} \frac{\Psi^{(m_j)}(\rho)}{m_j!}, 
\end{equation}
which vanishes for $\ell > J_0(\ell_0 - 1)$.

\begin{lemma}
\label{thm:Local expansion}
Uniformly for $|t| \leq \rho^{\chi}$,
\begin{equation}
\label{eq:Local expansion}
\Psi(\rho + it) + itn = \Psi(\rho) - \frac{t^2 \Psi''(\rho)}{2} + E(t) + O(\rho^{N_0}), \quad E(t) \ll \rho^{3\chi - 3 - \alpha} (\log \rho^{-1})^{k - 1} \ll 1. 
\end{equation}
\end{lemma}

\begin{proof}
Since $|t| \leq \rho^{\chi} < \rho$, Taylor's formula with integral remainder along the segment $[\rho, \rho + it]$ gives
\[
\Psi(\rho + it) = \sum_{\ell = 0}^{\ell_0 - 1} \frac{(it)^{\ell} \Psi^{(\ell)}(\rho)}{\ell!} + \frac{(it)^{\ell_0}}{(\ell_0 - 1)!} \int_0^1 (1 - x)^{\ell_0 - 1} \Psi^{(\ell_0)}(\rho + ixt)\, dx. 
\]
By \cref{eq:Derivative upper bound}, which is uniform in the imaginary part, $|t|^{\ell} |\Psi^{(\ell)}(\rho + ixt)| \ll \rho^{\ell(\chi - 1) - \alpha} (\log \rho^{-1})^{k - 1}$. For $\ell = \ell_0$ we have $\ell_0(\chi - 1) - \alpha > \frac{\ell_0 \alpha}{3} - \alpha \geq N_0$, so the remainder is $O(\rho^{N_0})$; for $3 \leq \ell < \ell_0$ the exponent is at least $3\chi - 3 - \alpha > 0$, which bounds $E(t)$. Finally, the term $\ell = 1$ is $it\Psi'(\rho) = -itn$.
\end{proof}

\noindent
Expanding $e^{E(t)}$ produces the coefficients $\lambda_\ell(\rho)$ of \cref{eq:Lambda}. In the integral over the major arcs only the even ones survive, and they are small.

\begin{lemma}
\label{thm:Size of lambda}
For every integer $\mu \geq 2$, $\lambda_{2\mu}(\rho)$ is real and, uniformly for $0 < \rho \leq 1/2$,
\begin{equation}
\label{eq:Size of lambda}
\frac{\lambda_{2\mu}(\rho)}{\Psi''(\rho)^{\mu}} \ll \rho^{\alpha} (\log \rho^{-1})^{1 - k}. 
\end{equation}
\end{lemma}

\begin{proof}
  Since $i^{2\mu} = (-1)^{\mu}$ and by \Cref{thm:Reality}, $\lambda_{2\mu}(\rho)$ is real. In a term of \cref{eq:Lambda} with $\ell = 2\mu$, every $m_j \geq 3$, so $3q \leq 2\mu$ and hence $\mu - q \geq 1$. By \cref{eq:Derivative upper bound}, such a term divided by $\Psi''(\rho)^{\mu}$ is
\[
\ll \frac{\rho^{-q\alpha - 2\mu} (\log \rho^{-1})^{q(k - 1)}}{\rho^{-\mu\alpha - 2\mu} (\log \rho^{-1})^{\mu(k - 1)}} = \left( \rho^{\alpha} (\log \rho^{-1})^{1 - k} \right)^{\mu - q} \ll \rho^{\alpha} (\log \rho^{-1})^{1 - k}
\]
as $\rho^{\alpha} (\log \rho^{-1})^{1 - k}$ is bounded on $(0, 1/2]$. There are finitely many terms.
\end{proof}

\noindent
We can now evaluate the integral $I$ over the major arcs.

\begin{lemma}
\label{thm:Major part}
We have
\begin{equation}
\label{eq:Major part}
I = \frac{\sqrt{2 \pi}\, P(\rho)}{\sqrt{\Psi''(\rho)}} \left\{ 1 + \sum_{2 \leq \mu \leq J_0(\ell_0 - 1)/2} \frac{(2\mu)!\, \lambda_{2\mu}(\rho)}{2^{\mu} \mu!\, \Psi''(\rho)^{\mu}} + O(\rho^{N_0}) \right\}.
\end{equation}
\end{lemma}

\begin{proof}
Using \Cref{thm:Local expansion} and $\exp(\Psi(\rho)) = P(\rho)$, for $|t| \leq \rho^{\chi}$, we have 
\[
e^{itn} \exp(\Psi(\rho + it)) = P(\rho)\, e^{-t^2 \Psi''(\rho)/2}\, e^{E(t)} \left( 1 + O(\rho^{N_0}) \right).
\]
Since $E(t) \ll 1$, we have $e^{E(t)} = \sum_{0 \leq q \leq J_0} E(t)^q/q! + O(|E(t)|^{J_0 + 1})$, and $|E(t)|^{J_0 + 1} \ll \rho^{N_0}$ because $(J_0 + 1)(3\chi - 3 - \alpha) > N_0$. Expanding the powers $E(t)^q$ and grouping by $\ell = m_1 + \cdots + m_q$,
\[
I = P(\rho) \int_{|t| \leq \rho^{\chi}} e^{-t^2 \Psi''(\rho)/2} \left\{ 1 + \sum_{3 \leq \ell \leq J_0(\ell_0 - 1)} \lambda_{\ell}(\rho)\, t^{\ell} + O(\rho^{N_0}) \right\} dt. 
\]
Since $\Psi''(\rho) > 0$ is real (\Cref{thm:Reality}), these are genuine Gaussian integrals. Odd $\ell$ contribute $0$ by symmetry. For $\ell = 2\mu$, using $\Gamma(\mu + \frac{1}{2}) = \frac{(2\mu)!\sqrt{\pi}}{4^{\mu} \mu!}$,
\[
\int_{-\infty}^{\infty} e^{-t^2 \Psi''(\rho)/2} t^{2\mu} dt = \Gamma\left(\mu + \tfrac{1}{2}\right) \left( \frac{2}{\Psi''(\rho)} \right)^{\mu + 1/2} = \frac{\sqrt{2 \pi} (2\mu)!}{2^{\mu} \mu!\, \Psi''(\rho)^{\mu + 1/2}},
\]
while for the tail, since $e^{-t^2\Psi''(\rho)/2} \leq e^{-\rho^{2\chi}\Psi''(\rho)/4} e^{-t^2\Psi''(\rho)/4}$ when $|t| \geq \rho^{\chi}$, we have 
\[
\int_{|t| \geq \rho^{\chi}} e^{-t^2 \Psi''(\rho)/2} t^{2\mu} dt \leq 2^{\mu + 1/2} e^{-\rho^{2\chi}\Psi''(\rho)/4} \frac{\sqrt{2 \pi} (2\mu)!}{2^{\mu} \mu!\, \Psi''(\rho)^{\mu + 1/2}} 
\]
By \cref{eq:Derivative upper bound}, $\rho^{2\chi} \Psi''(\rho) \gg \rho^{-(2 + \alpha - 2\chi)}$ with $2 + \alpha - 2\chi > 0$ as $\chi < 1 + \frac{\alpha}{2}$, so $e^{-\rho^{2\chi}\Psi''(\rho)/4} \ll \rho^{N_0}$. The $O(\rho^{N_0})$ term contributes $\ll \rho^{N_0} \sqrt{2\pi/\Psi''(\rho)}$. Factoring out $\sqrt{2\pi} P(\rho)/\sqrt{\Psi''(\rho)}$ and using that each term of the sum is $O(1)$ by \Cref{thm:Size of lambda}, we obtain \cref{eq:Major part}.
\end{proof}

\noindent
By \Cref{thm:Size of lambda}, all terms of the sum are of lower order.

\begin{corollary}
\label{thm:Dropping lambda}
Assume $N_0 > \alpha$. Then
\begin{equation}
\label{eq:Major part main term}
I = \frac{\sqrt{2 \pi}\, P(\rho)}{\sqrt{\Psi''(\rho)}} \left\{ 1 + O\left( \rho^{\alpha} (\log \rho^{-1})^{1 - k} \right) \right\}
= \frac{\sqrt{2 \pi}\, P(\rho)}{\sqrt{\Psi''(\rho)}} \left\{ 1 + O\left( n^{-\frac{\alpha}{1 + \alpha}} (\log n)^{-\frac{k - 1}{1 + \alpha}} \right) \right\}.
\end{equation}
\end{corollary}

\begin{proof}
The sum in \cref{eq:Major part} has a bounded number of terms, each $O(\rho^{\alpha} (\log \rho^{-1})^{1 - k})$ by \Cref{thm:Size of lambda}, and $\rho^{N_0} \ll \rho^{\alpha} (\log \rho^{-1})^{1 - k}$ since $N_0 > \alpha$. The second form follows from \Cref{thm:Size of rho}.
\end{proof}

\noindent
Adding the minor arcs, we obtain the saddle-point formula with an explicit relative error.

\begin{corollary}
\label{thm:Saddle point conclusion}
Let $\alpha/2 < M_1 < M_0$. Then, as $n \to \infty$,
\[
p_f(n) = \frac{e^{\rho n} P(\rho)}{\sqrt{2 \pi \Psi''(\rho)}} \left\{ 1 + O\left( \rho^{\alpha} (\log \rho^{-1})^{1 - k} + \rho^{M_1 - \frac{\alpha}{2}} (\log \rho^{-1})^{\frac{k - 1}{2}} \right) \right\} \sim \frac{e^{\rho n} P(\rho)}{\sqrt{2 \pi \Psi''(\rho)}}. 
\]
\end{corollary}

\begin{proof}
Insert \Cref{thm:Minor part} and \Cref{thm:Dropping lambda} into \cref{eq:Splitted integral}. Relative to the main term, the minor arcs contribute $\ll \rho^{M_1 + 1} \sqrt{\Psi''(\rho)} \ll \rho^{M_1 - \frac{\alpha}{2}} (\log \rho^{-1})^{\frac{k - 1}{2}}$ by \cref{eq:Derivative upper bound}, which tends to $0$ as $M_1 > \alpha/2$.
\end{proof}

\subsection{Proof of \texorpdfstring{\Cref{thm:Main theorem}}{Theorem 2.1}}

\noindent
\Cref{thm:Saddle point conclusion} gives $p_f(n) \sim e^{\rho n} P(\rho)/\sqrt{2 \pi \Psi''(\rho)}$, where $\rho = \rho_n$ is the solution of $\Psi'(\rho) + n = 0$ from \Cref{thm:Rho existence}. This is \cref{eq:Main formula}, which completes the proof of \Cref{thm:Main theorem}. \qed

\section{Proofs of the results}
\label{sec:Proofs}
In this section we prove \Cref{thm:Explicit saddle point} and the results of \Cref{sec:Dominant saddle point} and \Cref{sec:Precision}. We first prove \Cref{thm:Dominant saddle point} and \Cref{thm:Exponent expansion}, and then deduce \Cref{thm:Explicit saddle point} from \Cref{thm:Dominant saddle point}(i). The proof of \Cref{thm:Dominant saddle point}(iii) uses \Cref{thm:Accuracy of saddle point}(i), whose proof in \Cref{sec:Proof of accuracy} uses \Cref{thm:Main theorem}, the estimates of \Cref{sec:Contributions} and the lemmas of \Cref{sec:Preliminaries}. Throughout we assume the conditions (A1)-(A3), and $\rho_\alpha$ and $\Delta$ are as in \Cref{sec:Main result}. We write
\[
\Psi_\alpha(\rho) := \rho^{-\alpha} Q_0(\log \rho^{-1})
\]
for the residue of $M(z, \rho)$ at $z = \alpha$, so that the left-hand side of \cref{eq:Dominant saddle point equation} is $-\Psi_\alpha'(\rho)$. For $c \in \{-\varepsilon, \sigma_0\}$ we put
\[
K_c(\rho) := \frac{1}{2 \pi i} \int_{c - i\infty}^{c + i\infty} M(z, \rho)\, dz, \quad \rho > 0.
\]
By (A2), these lines contain no pole of $M(z, \rho)$.

By \Cref{thm:Saddle point conclusion}, the asymptotics of $p_f(n)$ is determined by $e^{n\rho_n} P(\rho_n)$ and $\Psi''(\rho_n)$, where the saddle point $\rho_n$ is defined only implicitly, by $-\Psi'(\rho_n) = n$. 
The idea 
is to replace $\Psi$ in this equation by its main part $\Psi_\alpha$, which leads to \cref{eq:Dominant saddle point equation} and its solution $\rho_\alpha$, and to treat $\Psi - \Psi_\alpha$ as a perturbation. By \Cref{thm:Psi decomposition} and \Cref{thm:Dominant saddle point location} below, this perturbation and its derivatives are smaller than those of $\Psi_\alpha$ by a factor of order at most $\rho^{\alpha - \sigma_0} (\log \rho^{-1})^{1 - k}$. So $\rho_n$ is close to $\rho_\alpha$, and the main task is to control how the exponent changes when $\rho_n$ is replaced by $\rho_\alpha$, or by a point closer to $\rho_n$.

\subsection{Preliminaries}
\label{sec:Preliminaries}

\noindent
We first split $\Psi$ into its main part $\Psi_\alpha$, the residues collected in $\Delta$, and a remainder $K_{-\varepsilon}$, which tends to $0$.

\begin{lemma}
\label{thm:Psi decomposition}
For $\rho > 0$,
\begin{equation}
\label{eq:Psi decomposition}
\Psi(\rho) = \Psi_\alpha(\rho) + K_{\sigma_0}(\rho), \qquad K_{\sigma_0}(\rho) = \Delta(\rho) + K_{-\varepsilon}(\rho).
\end{equation}
The functions $K_{\sigma_0}$, $K_{-\varepsilon}$ and $\Delta$ are real-valued. Moreover, $(-1)^{\ell} \Psi_\alpha^{(\ell)}(\rho) = \rho^{-\alpha - \ell} Q_\ell(\log \rho^{-1})$ for $\ell \geq 0$, and, uniformly for $0 < \rho \leq 1/2$,
\begin{equation}
\label{eq:Remainder bounds}
K_{\sigma_0}^{(\ell)}(\rho) \ll_\ell \rho^{-\sigma_0 - \ell}, \qquad K_{-\varepsilon}^{(\ell)}(\rho) \ll_\ell \rho^{\varepsilon - \ell}, \qquad \Delta^{(\ell)}(\rho) \ll_\ell \rho^{-\sigma_0 - \ell}.
\end{equation}
\end{lemma}

\begin{proof}
The proof of \Cref{thm:Derivatives of Psi} gives the first identity in \cref{eq:Psi decomposition}, the formula for $\Psi_\alpha^{(\ell)}$ and the bound for $K_{\sigma_0}^{(\ell)}$. The bound for $K_{-\varepsilon}^{(\ell)}$ follows in the same way on the line $\Re(z) = -\varepsilon$, where $|\rho^{-z - \ell}| = \rho^{\varepsilon - \ell}$ and the integrand is estimated by \cref{eq:Exponential bound} and the bounds for $\Gamma$ and $\zeta$ in the proof of \Cref{thm:0 integrals}; for $\ell = 0$ this is the bound $O(\rho^{\varepsilon})$ in the proof of \Cref{thm:Residues}. By the proof of \Cref{thm:Residues}, $\Psi(\rho)$ is the sum of the residues of $M(z, \rho)$ at the points of $\{0\} \cup \mathcal{A} \cup \mathcal{B}$, plus $K_{-\varepsilon}(\rho)$, where the sum over $\mathcal{B}$ is understood as in \Cref{thm:Infinite B}. This sum of residues is $\Psi_\alpha(\rho) + \Delta(\rho)$, and comparing with the first identity gives the second one. Since $f$ is real, $\overline{M(\overline{z}, \rho)} = M(z, \rho)$ for $\rho > 0$, so $K_{\sigma_0}(\rho)$ and $K_{-\varepsilon}(\rho)$ are real, and hence so is $\Delta(\rho)$. Finally, $\Delta^{(\ell)} = K_{\sigma_0}^{(\ell)} - K_{-\varepsilon}^{(\ell)}$, and $\rho^{\varepsilon - \ell} \leq \rho^{-\sigma_0 - \ell}$ for $0 < \rho \leq 1/2$.
\end{proof}

\noindent
For the rest of this section we write $u_\alpha := \log \rho_\alpha^{-1}$. The next lemma shows that $\rho_\alpha$ is well defined and close to $\rho_n$. The reason is that $\Psi' + n$ is increasing, with derivative $\Psi''$ of a fixed size near $\rho_\alpha$, while at $\rho_\alpha$ it equals the small quantity $K_{\sigma_0}'(\rho_\alpha)$, because $\Psi_\alpha' + n$ vanishes there.

\begin{lemma}
\label{thm:Dominant saddle point location}
There is $\rho^{\circ} \in (0, 1/2]$ such that $Q_\ell(\log \rho^{-1}) \asymp (\log \rho^{-1})^{k - 1}$ for $\ell \in \{1, 2\}$ and $0 < \rho \leq \rho^{\circ}$, where the implied constants are positive. For large $n$, the smallest positive solution $\rho_\alpha$ of \cref{eq:Dominant saddle point equation} exists, it is the only solution of \cref{eq:Dominant saddle point equation} in $(0, \rho^{\circ}]$, and $\rho_\alpha \to 0$ as $n \to \infty$. Moreover, for large $n$,
\begin{equation}
\label{eq:Local size}
\Psi''(x) \asymp \rho_\alpha^{-\alpha - 2} u_\alpha^{k - 1}, \qquad |\Psi'''(x)| \ll \rho_\alpha^{-\alpha - 3} u_\alpha^{k - 1}, \qquad \frac{\rho_\alpha}{2} \leq x \leq 2 \rho_\alpha,
\end{equation}
and $\rho_\alpha/2 < \rho_n < 2 \rho_\alpha$, with
\begin{equation}
\label{eq:Saddle points close}
\rho_n - \rho_\alpha = -\frac{K_{\sigma_0}'(\rho_\alpha)}{\Psi''(\xi)} \ll \rho_\alpha^{1 + \alpha - \sigma_0} u_\alpha^{1 - k}
\end{equation}
for some $\xi$ between $\rho_n$ and $\rho_\alpha$. In particular, $\rho_n \sim \rho_\alpha$.
\end{lemma}

\begin{proof}
By \cref{eq:Leading coefficient} and \Cref{thm:Positive leading coefficient}, $Q_\ell$ has degree $k - 1$ and positive leading coefficient $\Gamma(\alpha + \ell) C/\Gamma(\alpha + 2)$. So $Q_\ell(u)/u^{k - 1}$ has a positive limit as $u \to \infty$, and we can choose $u^{\circ} \geq \log 2$ with $Q_\ell(u) \asymp u^{k - 1}$ for $u \geq u^{\circ}$ and $\ell \in \{1, 2\}$. We take $\rho^{\circ} = e^{-u^{\circ}}$.

By \Cref{thm:Psi decomposition}, the left-hand side of \cref{eq:Dominant saddle point equation} is
\[
-\Psi_\alpha'(\rho) = \rho^{-\alpha - 1} Q_1(\log \rho^{-1})
\]
On $(0, \rho^{\circ}]$ this function is positive and continuous, it tends to $+\infty$ as $\rho \to 0^{+}$, and its derivative $-\rho^{-\alpha - 2} Q_2(\log \rho^{-1})$ is negative. For $n > -\Psi_\alpha'(\rho^{\circ})$ the equation $-\Psi_\alpha'(\rho) = n$ therefore has exactly one solution in $(0, \rho^{\circ}]$, and all other positive solutions exceed $\rho^{\circ}$, so this solution is $\rho_\alpha$. If $0 < r < \rho^{\circ}$ and $n > -\Psi_\alpha'(r)$, then $\rho_\alpha < r$, because $-\Psi_\alpha'$ is decreasing on $(0, \rho^{\circ}]$ and $-\Psi_\alpha'(\rho_\alpha) = n > -\Psi_\alpha'(r)$. Thus $\rho_\alpha \to 0$ and $u_\alpha \to \infty$.

Now let $n$ be so large that $\rho_\alpha \leq 1/4$, and let $\rho_\alpha/2 \leq x \leq 2\rho_\alpha$. Then $x \leq 1/2$ and $x/\rho_\alpha \in [\frac{1}{2}, 2]$. Moreover $|\log x^{-1} - u_\alpha| \leq \log 2 \leq u_\alpha/2$, as $u_\alpha \geq \log 4$, so that $\log x^{-1}$ lies between $u_\alpha/2$ and $3u_\alpha/2$. Hence
\[
x^{-\alpha - 2} (\log x^{-1})^{k - 1} \asymp \rho_\alpha^{-\alpha - 2} u_\alpha^{k - 1}, \qquad x^{-\alpha - 3} (\log x^{-1})^{k - 1} \asymp \rho_\alpha^{-\alpha - 3} u_\alpha^{k - 1},
\]
and \cref{eq:Local size} follows from \cref{eq:Derivative upper bound} with $t = 0$ and $\ell = 2, 3$, which gives the upper bounds for $\Psi''$ and $|\Psi'''|$ and the lower bound for $\Psi''$. We fix $c_5 > 0$ such that $\Psi''(x) \geq c_5 \rho_\alpha^{-\alpha - 2} u_\alpha^{k - 1}$ for $\rho_\alpha/2 \leq x \leq 2\rho_\alpha$, which is possible by the lower bound in \cref{eq:Local size}.

Since $\Psi_\alpha'(\rho_\alpha) = -n$, we get from \cref{eq:Psi decomposition} and \cref{eq:Remainder bounds}
\[
\Psi'(\rho_\alpha) + n = K_{\sigma_0}'(\rho_\alpha) \ll \rho_\alpha^{-\sigma_0 - 1} 
  = \rho_\alpha^{\alpha - \sigma_0} u_\alpha^{1 - k} \cdot \rho_\alpha^{-\alpha - 1} u_\alpha^{k - 1},
\]
and hence $K_{\sigma_0}'(\rho_\alpha) = o(\rho_\alpha^{-\alpha - 1} u_\alpha^{k - 1})$, because $\sigma_0 < \alpha$. Integrating $\Psi'' \geq c_5 \rho_\alpha^{-\alpha - 2} u_\alpha^{k - 1}$ over $[\rho_\alpha, 2\rho_\alpha]$ and over $[\rho_\alpha/2, \rho_\alpha]$, which have lengths $\rho_\alpha$ and $\rho_\alpha/2$, we find
\[
\Psi'(2\rho_\alpha) - \Psi'(\rho_\alpha) \geq c_5 \rho_\alpha^{-\alpha - 1} u_\alpha^{k - 1}, \qquad \Psi'(\rho_\alpha) - \Psi'(\rho_\alpha/2) \geq \frac{c_5}{2} \rho_\alpha^{-\alpha - 1} u_\alpha^{k - 1}.
\]
Since $\Psi'(\rho_\alpha) + n = K_{\sigma_0}'(\rho_\alpha)$, this gives
\[
\Psi'(2\rho_\alpha) + n \geq K_{\sigma_0}'(\rho_\alpha) + c_5 \rho_\alpha^{-\alpha - 1} u_\alpha^{k - 1}, \qquad \Psi'(\rho_\alpha/2) + n \leq K_{\sigma_0}'(\rho_\alpha) - \frac{c_5}{2} \rho_\alpha^{-\alpha - 1} u_\alpha^{k - 1}
\]
For large $n$, $|K_{\sigma_0}'(\rho_\alpha)| < \frac{c_5}{2} \rho_\alpha^{-\alpha - 1} u_\alpha^{k - 1}$, so $\Psi' + n$ is positive at $2\rho_\alpha$ and negative at $\rho_\alpha/2$. As $\Psi'' > 0$ by \Cref{thm:Reality}, $\Psi' + n$ is strictly increasing, and its zero $\rho_n$ lies between $\rho_\alpha/2$ and $2\rho_\alpha$.

Finally, by the mean value theorem, $\Psi'(\rho_n) - \Psi'(\rho_\alpha) = \Psi''(\xi)(\rho_n - \rho_\alpha)$ for some $\xi$ 
  between $\rho_n$ and $\rho_\alpha$. 
  The left-hand side is $-n - \Psi'(\rho_\alpha) = -K_{\sigma_0}'(\rho_\alpha)$, which proves the identity in \cref{eq:Saddle points close}. Since $\xi \in (\rho_\alpha/2, 2\rho_\alpha)$, we have $\Psi''(\xi) \geq c_5 \rho_\alpha^{-\alpha - 2} u_\alpha^{k - 1}$, and \cref{eq:Remainder bounds} gives
\[
|\rho_n - \rho_\alpha| \leq \frac{\rho_\alpha^{\alpha + 2} u_\alpha^{1 - k}}{c_5} |K_{\sigma_0}'(\rho_\alpha)| \ll \rho_\alpha^{\alpha + 2} u_\alpha^{1 - k} \rho_\alpha^{-\sigma_0 - 1} = \rho_\alpha^{1 + \alpha - \sigma_0} u_\alpha^{1 - k}.
\]
Dividing by $\rho_\alpha$, we see that $\rho_n \sim \rho_\alpha$.
\end{proof}

\begin{remark}
Note that \Cref{thm:Dominant saddle point location} uses only (A1)-(A3), and not the conditions on $\sigma_0$ in \Cref{thm:Dominant saddle point}.
\end{remark}

\subsection{Proof of \texorpdfstring{\Cref{thm:Dominant saddle point}}{Theorem 2.2}}

\noindent
Let $\Phi(x) := nx + \Psi(x)$, so that $\Phi'(\rho_n) = 0$, $\Phi'' = \Psi''$ and $e^{\Phi(x)} = e^{nx} P(x)$. By \Cref{thm:Saddle point conclusion},
\begin{equation}
\label{eq:Saddle form}
p_f(n) \sim \frac{e^{\Phi(\rho_n)}}{\sqrt{2 \pi \Psi''(\rho_n)}}.
\end{equation}
Since $\rho_n$ is not explicit, we expand the exponent and the denominator around $\rho_\alpha$. As $\Phi'(\rho_n) = 0$ and $\Phi'(\rho_\alpha)$ is small, the exponent changes only by a second-order quantity, approximately $-\Phi'(\rho_\alpha)^2/(2\Psi''(\rho_\alpha))$. Up to small errors, this is the correction term in \cref{eq:Dominant formula corrected}, and it is negligible when $\sigma_0 < \alpha/2$. We treat the denominator in \Cref{thm:Prefactor} and the exponent in \Cref{thm:Exponent at rho alpha}; together they give (i) and (ii), and \Cref{thm:Negative poles} shows that finitely many residues with negative real part may be omitted. For (iii), \Cref{thm:Critical point} locates the zero $\hat{\rho}$ of $\Phi_\Delta'$, and \Cref{thm:Newton steps} shows that the Newton points approach it quadratically.

To pass from $\rho_n$ to $\rho_\alpha$ we put $d_n := \rho_n - \rho_\alpha$. Since $-\Psi_\alpha'(\rho_\alpha) = n$, \Cref{thm:Psi decomposition} gives $\Phi'(\rho_\alpha) = n + \Psi_\alpha'(\rho_\alpha) + K_{\sigma_0}'(\rho_\alpha) = K_{\sigma_0}'(\rho_\alpha)$, and together with \Cref{thm:Dominant saddle point location} we obtain
\begin{equation}
\label{eq:Sizes at rho alpha}
\Phi'(\rho_\alpha) = K_{\sigma_0}'(\rho_\alpha) \ll \rho_\alpha^{-\sigma_0 - 1}, \qquad \Psi''(\rho_\alpha) \asymp \rho_\alpha^{-\alpha - 2} u_\alpha^{k - 1}, \qquad d_n \ll \rho_\alpha^{\alpha + 2} u_\alpha^{1 - k} |\Phi'(\rho_\alpha)|.
\end{equation}

\begin{lemma}
\label{thm:Prefactor}
As $n \to \infty$,
\begin{equation}
\label{eq:Prefactor}
\Psi''(\rho_\alpha) = \rho_\alpha^{-\alpha - 2} Q_2(u_\alpha) \left( 1 + O\left( \rho_\alpha^{\alpha - \sigma_0} u_\alpha^{1 - k} \right) \right), \qquad \Psi''(\rho_n) \sim \rho_\alpha^{-\alpha - 2} Q_2(u_\alpha).
\end{equation}
\end{lemma}

\begin{proof}
By \Cref{thm:Derivatives of Psi}, we have $\Psi''(\rho_\alpha) = \rho_\alpha^{-\alpha - 2} Q_2(u_\alpha) + O(\rho_\alpha^{-\sigma_0 - 2})$, and the mean value theorem, \cref{eq:Local size} and \cref{eq:Sizes at rho alpha} give $|\Psi''(\rho_n) - \Psi''(\rho_\alpha)| \ll \rho_\alpha^{-\alpha - 3} u_\alpha^{k - 1} |d_n| \ll \rho_\alpha^{-\sigma_0 - 2}$. Since $Q_2(u_\alpha) \asymp u_\alpha^{k - 1}$ by \Cref{thm:Dominant saddle point location}, the main term $\rho_\alpha^{-\alpha - 2} Q_2(u_\alpha)$ is of order $\rho_\alpha^{-\alpha - 2} u_\alpha^{k - 1}$, so both errors $O(\rho_\alpha^{-\sigma_0 - 2})$ are smaller than the main term by a factor $O(\rho_\alpha^{\alpha - \sigma_0} u_\alpha^{1 - k}) = o(1)$.
\end{proof}

\noindent
In the exponent, the correction term appears through a Taylor expansion at $\rho_\alpha$. The following lemma uses only $\sigma_0 < \alpha$.

\begin{lemma}
\label{thm:Exponent at rho alpha}
As $n \to \infty$,
\begin{equation}
\label{eq:Exponent expansion}
\Phi(\rho_n) = n\rho_\alpha + \rho_\alpha^{-\alpha} Q_0(u_\alpha) + \Delta(\rho_\alpha) - \frac{\rho_\alpha^{\alpha + 2} \Delta'(\rho_\alpha)^2}{2 Q_2(u_\alpha)} + O(E_\alpha),
\end{equation}
where
\[
E_\alpha := \rho_\alpha^{2\alpha - 3\sigma_0} u_\alpha^{2(1 - k)} + \rho_\alpha^{\alpha - \sigma_0 + \varepsilon} u_\alpha^{1 - k} + \rho_\alpha^{\varepsilon}.
\]
\end{lemma}

\begin{proof}
In the exponent we start from \cref{eq:Psi decomposition} and \cref{eq:Remainder bounds}, which give
\begin{equation}
\label{eq:Exponent at rho alpha}
\Phi(\rho_\alpha) = n\rho_\alpha + \rho_\alpha^{-\alpha} Q_0(u_\alpha) + \Delta(\rho_\alpha) + K_{-\varepsilon}(\rho_\alpha) = n\rho_\alpha + \rho_\alpha^{-\alpha} Q_0(u_\alpha) + \Delta(\rho_\alpha) + O(\rho_\alpha^{\varepsilon}).
\end{equation}
By Taylor's formula at $\rho_\alpha$, there are $\xi_1, \xi_2$ between $\rho_n$ and $\rho_\alpha$ such that
\[
\begin{aligned}
0 &= \Phi'(\rho_n) = \Phi'(\rho_\alpha) + \Psi''(\rho_\alpha) d_n + \frac{1}{2} \Psi'''(\xi_1) d_n^2, \\
\Phi(\rho_n) &= \Phi(\rho_\alpha) + \Phi'(\rho_\alpha) d_n + \frac{1}{2} \Psi''(\rho_\alpha) d_n^2 + \frac{1}{6} \Psi'''(\xi_2) d_n^3.
\end{aligned}
\]
Without the term with $\Psi'''$, the first relation would give $d_n = -\Phi'(\rho_\alpha)/\Psi''(\rho_\alpha)$, so we put $y := d_n + \Phi'(\rho_\alpha)/\Psi''(\rho_\alpha)$. Since $|\Psi'''(\xi_1)|/\Psi''(\rho_\alpha) \ll \rho_\alpha^{-1}$ by \cref{eq:Local size}, the first relation and \cref{eq:Sizes at rho alpha} give
\[
y = -\frac{\Psi'''(\xi_1) d_n^2}{2\Psi''(\rho_\alpha)} \ll \frac{d_n^2}{\rho_\alpha} \ll \rho_\alpha^{2\alpha + 3} u_\alpha^{2(1 - k)} \Phi'(\rho_\alpha)^2.
\]
Substituting $d_n = -\Phi'(\rho_\alpha)/\Psi''(\rho_\alpha) + y$ into the second relation, the terms linear in $y$ cancel:
\[
\Phi'(\rho_\alpha) d_n + \frac{1}{2} \Psi''(\rho_\alpha) d_n^2 = -\frac{\Phi'(\rho_\alpha)^2}{2\Psi''(\rho_\alpha)} + \frac{1}{2} \Psi''(\rho_\alpha) y^2.
\]
Hence, bounding the cubic term by \cref{eq:Local size} and then inserting the bounds for $y$ and $d_n$,
\[
\begin{aligned}
\Phi(\rho_n) &= \Phi(\rho_\alpha) - \frac{\Phi'(\rho_\alpha)^2}{2\Psi''(\rho_\alpha)} + O\left( \Psi''(\rho_\alpha) y^2 + \rho_\alpha^{-\alpha - 3} u_\alpha^{k - 1} |d_n|^3 \right) \\
&= \Phi(\rho_\alpha) - \frac{\Phi'(\rho_\alpha)^2}{2\Psi''(\rho_\alpha)} + O\left( \rho_\alpha^{3\alpha + 4} u_\alpha^{3(1 - k)} \Phi'(\rho_\alpha)^4 + \rho_\alpha^{2\alpha + 3} u_\alpha^{2(1 - k)} |\Phi'(\rho_\alpha)|^3 \right).
\end{aligned}
\]
As $\Phi'(\rho_\alpha) \ll \rho_\alpha^{-\sigma_0 - 1}$, the error term is $\ll \rho_\alpha^{3\alpha - 4\sigma_0} u_\alpha^{3(1 - k)} + \rho_\alpha^{2\alpha - 3\sigma_0} u_\alpha^{2(1 - k)} \ll \rho_\alpha^{2\alpha - 3\sigma_0} u_\alpha^{2(1 - k)}$, because the first term is the second one multiplied by $\rho_\alpha^{\alpha - \sigma_0} u_\alpha^{1 - k} = o(1)$. The term $\Phi'(\rho_\alpha)^2/(2\Psi''(\rho_\alpha))$ is close to the one in \cref{eq:Dominant formula corrected}. Indeed, $\Phi'(\rho_\alpha) = K_{\sigma_0}'(\rho_\alpha) = \Delta'(\rho_\alpha) + K_{-\varepsilon}'(\rho_\alpha)$ by \cref{eq:Psi decomposition}, so that $\Phi'(\rho_\alpha)^2 - \Delta'(\rho_\alpha)^2 = K_{-\varepsilon}'(\rho_\alpha) \big( 2\Delta'(\rho_\alpha) + K_{-\varepsilon}'(\rho_\alpha) \big)$, and \cref{eq:Remainder bounds} and \cref{eq:Sizes at rho alpha} give
\[
\frac{|\Phi'(\rho_\alpha)^2 - \Delta'(\rho_\alpha)^2|}{\Psi''(\rho_\alpha)} \ll \rho_\alpha^{\alpha + 2} u_\alpha^{1 - k} \rho_\alpha^{-\sigma_0 - 1} \rho_\alpha^{\varepsilon - 1} = \rho_\alpha^{\alpha - \sigma_0 + \varepsilon} u_\alpha^{1 - k}.
\]
Since $\Delta'(\rho_\alpha)^2 \ll \rho_\alpha^{-2\sigma_0 - 2}$ and $1/Q_2(u_\alpha) \ll u_\alpha^{1 - k}$, the first relation in \cref{eq:Prefactor} gives
\[
\begin{aligned}
\frac{\Delta'(\rho_\alpha)^2}{2\Psi''(\rho_\alpha)} &= \frac{\rho_\alpha^{\alpha + 2} \Delta'(\rho_\alpha)^2}{2 Q_2(u_\alpha)} + O\left( \rho_\alpha^{-2\sigma_0 - 2} \rho_\alpha^{\alpha + 2} u_\alpha^{1 - k} \rho_\alpha^{\alpha - \sigma_0} u_\alpha^{1 - k} \right) \\
&= \frac{\rho_\alpha^{\alpha + 2} \Delta'(\rho_\alpha)^2}{2 Q_2(u_\alpha)} + O\left( \rho_\alpha^{2\alpha - 3\sigma_0} u_\alpha^{2(1 - k)} \right).
\end{aligned}
\]
Combining these estimates with \cref{eq:Exponent at rho alpha}, we obtain \cref{eq:Exponent expansion}. The first term of $E_\alpha$ comes from the Taylor expansion and from \cref{eq:Prefactor}, the second from replacing $\Phi'(\rho_\alpha)$ by $\Delta'(\rho_\alpha)$, and the third from $K_{-\varepsilon}(\rho_\alpha)$.
\end{proof}

\noindent
We can now prove (i) and (ii). If $\sigma_0 < 2\alpha/3$, then all exponents of $\rho_\alpha$ in $E_\alpha$ are positive, so $E_\alpha = o(1)$, and inserting \cref{eq:Exponent expansion} and \cref{eq:Prefactor} into \cref{eq:Saddle form} gives \cref{eq:Dominant formula corrected}. This proves (ii). If even $\sigma_0 < \alpha/2$, then (ii) holds, since $\alpha/2 < 2\alpha/3$. Moreover $\Delta'(\rho_\alpha) \ll \rho_\alpha^{-\sigma_0 - 1}$ by \cref{eq:Remainder bounds} and $Q_2(u_\alpha) \asymp u_\alpha^{k - 1}$ by \Cref{thm:Dominant saddle point location}, so that
\[
\frac{\rho_\alpha^{\alpha + 2} \Delta'(\rho_\alpha)^2}{2 Q_2(u_\alpha)} \ll \rho_\alpha^{\alpha - 2\sigma_0} u_\alpha^{1 - k} = o(1)
\]
so the last term of the exponent in \cref{eq:Dominant formula corrected} can be dropped. This gives \cref{eq:Dominant formula} and proves (i).

\begin{lemma}
\label{thm:Negative poles}
Let $\Delta_{<0}$ be the sum of the residues of $M(z, \rho)$ at finitely many poles $p$ of $D(s)$ with $\Re(p) < 0$. There are $\theta_0 > 0$ and $c_6 \geq 0$ such that $\Delta_{<0}^{(j)}(\rho) \ll_j \rho^{\theta_0 - j} (\log \rho^{-1})^{c_6}$ for $j \geq 0$ and $0 < \rho \leq 1/2$, and \cref{eq:Dominant formula} and \cref{eq:Dominant formula corrected} remain valid when $\Delta$ is replaced by $\Delta - \Delta_{<0}$.
\end{lemma}

\begin{proof}
By \Cref{thm:Res at p}, each of these residues is $\rho^{-p}$ times a polynomial in $\log \rho^{-1}$, and $|\rho^{-p}| = \rho^{-\Re(p)}$. This gives the bounds, with $\theta_0$ the least of the numbers $-\Re(p)$, for instance. Replacing $\Delta$ by $\Delta - \Delta_{<0}$ changes $\Delta(\rho_\alpha)$ by $\Delta_{<0}(\rho_\alpha) = o(1)$. Since $\Delta'^2 - (\Delta' - \Delta_{<0}')^2 = \Delta_{<0}' (2\Delta' - \Delta_{<0}')$, it changes the last term of the exponent in \cref{eq:Dominant formula corrected} by
\[
\begin{aligned}
\frac{\rho_\alpha^{\alpha + 2} \left| \Delta'(\rho_\alpha)^2 - \left( \Delta'(\rho_\alpha) - \Delta_{<0}'(\rho_\alpha) \right)^2 \right|}{2 Q_2(u_\alpha)} &\ll \rho_\alpha^{\alpha + 2} \left( \rho_\alpha^{-\sigma_0 - 1} + \rho_\alpha^{\theta_0 - 1} u_\alpha^{c_6} \right) \rho_\alpha^{\theta_0 - 1} u_\alpha^{c_6} \\
&\ll \rho_\alpha^{\alpha - \sigma_0 + \theta_0} u_\alpha^{2c_6} = o(1),
\end{aligned}
\]
so \cref{eq:Dominant formula} and \cref{eq:Dominant formula corrected} remain valid.
\end{proof}

\noindent
It remains to prove (iii). If $\sigma_0 \geq 2\alpha/3$, the error $E_\alpha$ need not tend to $0$, so the point $\rho_\alpha$ is not accurate enough. We therefore approach $\rho_n$ by Newton's method. The remainder $K_{-\varepsilon}$ is negligible, so it suffices to approach the zero $\hat{\rho}$ of $\Phi_\Delta'$, which is very close to $\rho_n$, and each Newton step squares the relative distance to $\hat{\rho}$. Put $I := [\rho_\alpha/2, 2\rho_\alpha]$.

\begin{lemma}
\label{thm:Critical point}
For large $n$ we have
\[
\Phi_\Delta''(x) \geq \frac{c_5}{2} \rho_\alpha^{-\alpha - 2} u_\alpha^{k - 1}, \qquad \Phi_\Delta^{(\ell)}(x) \ll_\ell \rho_\alpha^{-\alpha - \ell} u_\alpha^{k - 1}, \quad \ell \geq 2,\quad x \in I.
\]
The derivative $\Phi_\Delta'$ has exactly one zero $\hat{\rho}$ in $I$, and
\begin{equation}
\label{eq:Critical point}
|\hat{\rho} - \rho_\alpha| \leq \frac{2}{c_5} \rho_\alpha^{\alpha + 2} u_\alpha^{1 - k} |\Delta'(\rho_\alpha)| \ll \rho_\alpha^{1 + \alpha - \sigma_0} u_\alpha^{1 - k}, \qquad |\rho_n - \hat{\rho}| \ll \rho_\alpha^{1 + \alpha + \varepsilon} u_\alpha^{1 - k}.
\end{equation}
\end{lemma}

\begin{proof}
By \cref{eq:Psi decomposition}, $\Phi = \Phi_\Delta + K_{-\varepsilon}$, and $K_{-\varepsilon}^{(\ell)}(x) \ll \rho_\alpha^{\varepsilon - \ell} = o(\rho_\alpha^{-\alpha - \ell} u_\alpha^{k - 1})$ on $I$ by \cref{eq:Remainder bounds}. Since $\Phi_\Delta^{(\ell)} = \Psi^{(\ell)} - K_{-\varepsilon}^{(\ell)}$ for $\ell \geq 2$, the bounds for $\Phi_\Delta^{(\ell)}$ follow from \cref{eq:Derivative upper bound}, \cref{eq:Local size} and the choice of $c_5$ in the proof of \Cref{thm:Dominant saddle point location}. Moreover, 
  $$
  \Phi_\Delta'(\rho_\alpha) = n + \Psi_\alpha'(\rho_\alpha) + \Delta'(\rho_\alpha) = \Delta'(\rho_\alpha) \ll \rho_\alpha^{-\sigma_0 - 1} = \rho_\alpha^{\alpha - \sigma_0} u_\alpha^{1 - k} \cdot \rho_\alpha^{-\alpha - 1} u_\alpha^{k - 1}.
  $$ 
  Arguing as in the proof of \Cref{thm:Dominant saddle point location}, with $\Phi_\Delta$ in place of $\Phi$ and $\Delta'$ in place of $K_{\sigma_0}'$, we see that $\Phi_\Delta'$ is positive at $2\rho_\alpha$, negative at $\rho_\alpha/2$, and increasing on $I$. So it has exactly one zero $\hat{\rho}$ in $I$, and the mean value theorem gives the first bound in \cref{eq:Critical point}. Since $\rho_n \in I$ and $\Phi_\Delta'(\rho_n) = \Phi'(\rho_n) - K_{-\varepsilon}'(\rho_n) = -K_{-\varepsilon}'(\rho_n) \ll \rho_\alpha^{\varepsilon - 1}$, the mean value theorem for $\Phi_\Delta'$ between $\hat{\rho}$ and $\rho_n$ gives the second bound.
\end{proof}

\begin{lemma}
\label{thm:Newton steps}
Let $J \geq 0$ be fixed, and let $x_0 := \rho_\alpha$ and $x_{j + 1} := x_j - \Phi_\Delta'(x_j)/\Phi_\Delta''(x_j)$. For large $n$, the points $x_0, \ldots, x_J$ are well defined and lie in $I$, and
\begin{equation}
\label{eq:Newton distance}
\frac{|x_J - \rho_n|}{\rho_\alpha} \ll \rho_\alpha^{2^J (\alpha - \sigma_0)} u_\alpha^{2^J (1 - k)} + \rho_\alpha^{\alpha + \varepsilon} u_\alpha^{1 - k}.
\end{equation}
\end{lemma}

\begin{proof}
For $x \in I$, Taylor's formula gives $0 = \Phi_\Delta'(\hat{\rho}) = \Phi_\Delta'(x) + \Phi_\Delta''(x)(\hat{\rho} - x) + \frac{1}{2} \Phi_\Delta'''(\xi)(\hat{\rho} - x)^2$ for some $\xi \in I$, so that, with a constant $c_7 \geq 1$ that exists by the bounds for $\Phi_\Delta''$ and $\Phi_\Delta'''$ on $I$,
\[
\left| x - \frac{\Phi_\Delta'(x)}{\Phi_\Delta''(x)} - \hat{\rho} \right| = \frac{|\Phi_\Delta'''(\xi)|}{2 \Phi_\Delta''(x)} (x - \hat{\rho})^2 \leq \frac{c_7}{\rho_\alpha} (x - \hat{\rho})^2.
\]
The left-hand side is the distance from the Newton step at $x$ to $\hat{\rho}$, so each step squares this distance, up to the factor $c_7/\rho_\alpha$. Put $\delta_j := c_7 |x_j - \hat{\rho}|/\rho_\alpha$. Since $x_0 = \rho_\alpha$, \cref{eq:Critical point} gives $\delta_0 \ll \rho_\alpha^{\alpha - \sigma_0} u_\alpha^{1 - k}$, so $\delta_0 \leq \frac{1}{4}$ for large $n$. If $x_j \in I$ and $\delta_j \leq \delta_0$, then $x_{j + 1}$ is defined, as $\Phi_\Delta'' > 0$ on $I$, and the last display with $x = x_j$, multiplied by $c_7/\rho_\alpha$, gives $\delta_{j + 1} \leq \delta_j^2 \leq \delta_0$. So $|x_{j + 1} - \hat{\rho}| \leq \rho_\alpha \delta_0 \leq \rho_\alpha/4$, and $x_{j + 1} \in I$ because $\hat{\rho} = \rho_\alpha (1 + o(1))$. By induction all $x_j$ are well defined, lie in $I$, and satisfy $\delta_j \leq \delta_0^{2^j}$. Hence \cref{eq:Newton distance} follows from the triangle inequality $|x_J - \rho_n| \leq |x_J - \hat{\rho}| + |\hat{\rho} - \rho_n|$, since $J$ is fixed.
\end{proof}

\noindent
We can now prove (iii). By \cref{eq:Size of n rho} and $\rho_n \sim \rho_\alpha$, which gives $\log \rho_n^{-1} = u_\alpha + o(1)$, we have $(n\rho_n)^{-1/2} \asymp \rho_\alpha^{\alpha/2} u_\alpha^{(1 - k)/2}$, and
\[
\frac{\rho_\alpha^{2^J (\alpha - \sigma_0)} u_\alpha^{2^J (1 - k)}}{\rho_\alpha^{\alpha/2} u_\alpha^{(1 - k)/2}} = \rho_\alpha^{2^J (\alpha - \sigma_0) - \alpha/2} u_\alpha^{(1 - k)(2^J - 1/2)} \to 0, \qquad \frac{\rho_\alpha^{\alpha + \varepsilon} u_\alpha^{1 - k}}{\rho_\alpha^{\alpha/2} u_\alpha^{(1 - k)/2}} \to 0
\]
by the choice of $J$. Since $x_J$ and $\rho_n$ lie in $I$, the mean value theorem gives $|\log x_J^{-1} - \log \rho_n^{-1}| \leq 2|x_J - \rho_n|/\rho_\alpha$, and \cref{eq:Newton distance} shows that $\log x_J^{-1} - \log \rho_n^{-1} = o\big( (n\rho_n)^{-1/2} \big)$, and \Cref{thm:Accuracy of saddle point}(i), applied with $\tilde{\rho} = x_J$, gives $p_f(n) \sim e^{n x_J} P(x_J)/\sqrt{2\pi \Psi''(x_J)}$. Here 
$$
n x_J + \log P(x_J) = \Phi(x_J) = \Phi_\Delta(x_J) + K_{-\varepsilon}(x_J) = \Phi_\Delta(x_J) + O(\rho_\alpha^{\varepsilon}),
$$ 
and $\Psi''(x_J) = \Phi_\Delta''(x_J) + K_{-\varepsilon}''(x_J) = \Phi_\Delta''(x_J) (1 + o(1))$, as $K_{-\varepsilon}''(x_J) \ll \rho_\alpha^{\varepsilon - 2}$ and $\Phi_\Delta''(x_J) \gg \rho_\alpha^{-\alpha - 2} u_\alpha^{k - 1}$. This gives the first form of \cref{eq:Newton formula}. Moreover $|x_J - \rho_\alpha| \leq |x_J - \hat{\rho}| + |\hat{\rho} - \rho_\alpha| \ll \rho_\alpha^{1 + \alpha - \sigma_0} u_\alpha^{1 - k}$, so $|\Psi''(x_J) - \Psi''(\rho_\alpha)| \ll \rho_\alpha^{-\sigma_0 - 2}$ by \cref{eq:Local size}, and $\Psi''(x_J) \sim \rho_\alpha^{-\alpha - 2} Q_2(u_\alpha)$ by \cref{eq:Prefactor}. This gives the second form of \cref{eq:Newton formula}. If $\Delta$ is replaced by $\Delta - \Delta_{<0}$, the same argument applies with $K_{-\varepsilon}$ replaced by $K_{-\varepsilon} + \Delta_{<0}$, whose derivatives of order at most $3$ satisfy the bounds used above with $\varepsilon$ replaced by $\min(\varepsilon, \theta_0)$, up to a factor $u_\alpha^{c_6}$, by \Cref{thm:Negative poles}. \qed

\begin{remark}
The proof gives explicit error terms. The exponent in \cref{eq:Dominant formula corrected} has error $O(E_\alpha)$, the exponent in \cref{eq:Dominant formula} has error $O(\rho_\alpha^{\alpha - 2\sigma_0} u_\alpha^{1 - k} + \rho_\alpha^{\varepsilon})$, and the denominator has relative error $O(\rho_\alpha^{\alpha - \sigma_0} u_\alpha^{1 - k})$. These are in addition to the relative error of \Cref{thm:Saddle point conclusion}. Even when $\sigma_0 < \alpha/2$, the last term of the exponent in \cref{eq:Dominant formula corrected} can be numerically noticeable, and \cref{eq:Dominant formula corrected} is then more accurate than \cref{eq:Dominant formula}. For example, if $z = 0$ is a pole of $D(s)$ of order $h_0 \geq 1$ and no pole $p \neq \alpha$ has $\Re(p) > 0$, this term has size $\rho_\alpha^{\alpha} (\log \rho_\alpha^{-1})^{2h_0 + 1 - k}$.
\end{remark}

\subsection{Proof of \texorpdfstring{\Cref{thm:Exponent expansion}}{Theorem 2.3}}

\noindent
The Newton points in \Cref{thm:Dominant saddle point}(iii) are not explicit in $\rho_\alpha$. Instead, we now expand the minimum value $\Phi_\Delta(\hat{\rho})$ around $\rho_\alpha$. After the substitution $x = \rho_\alpha (1 + \eta)$ and division by $\rho_\alpha^{-\alpha} u_\alpha^{k - 1}$, the function $\Phi_\Delta$ becomes a series in $\eta$ with bounded coefficients plus a perturbation, coming from $\Delta$, of relative size $\rho_\alpha^{\alpha - \sigma_1} u_\alpha^{h_1 + 1 - k}$. The terms $T_m$ are the terms of the minimum value of this function, ordered by powers of this small quantity.

Let $u_\alpha$ and $c_5$ be as above, and let $I = [\rho_\alpha/2, 2\rho_\alpha]$. Since $\sigma_1 < \alpha$, we have $\rho_\alpha^{\alpha - \sigma_1} u_\alpha^{h_1 + 1 - k} \to 0$, and \cref{eq:Delta bound} gives $\Delta'(\rho_\alpha) \ll \rho_\alpha^{-\sigma_1 - 1} u_\alpha^{h_1}$. By \Cref{thm:Critical point}, $\Phi_\Delta'' \geq \frac{c_5}{2} \rho_\alpha^{-\alpha - 2} u_\alpha^{k - 1}$ on $I$, $\Phi_\Delta'$ has exactly one zero $\hat{\rho} \in I$, and, by the first inequality in \cref{eq:Critical point},
\[
\hat{\rho} - \rho_\alpha \ll \rho_\alpha^{1 + \alpha - \sigma_1} u_\alpha^{h_1 + 1 - k}, \qquad \rho_n - \hat{\rho} \ll \rho_\alpha^{1 + \alpha + \varepsilon} u_\alpha^{1 - k}.
\]
We first replace $\rho_n$ by $\hat{\rho}$, which does not depend on the remainder $K_{-\varepsilon}$.

\begin{lemma}
\label{thm:Critical form}
As $n \to \infty$,
\begin{equation}
\label{eq:Critical form}
p_f(n) \sim \frac{e^{\Phi_\Delta(\hat{\rho})}}{\sqrt{2 \pi \Phi_\Delta''(\hat{\rho})}}.
\end{equation}
\end{lemma}

\begin{proof}
Since $\Phi = \Phi_\Delta + K_{-\varepsilon}$ and $\Phi_\Delta'(\hat{\rho}) = 0$, Taylor's formula for $\Phi_\Delta$ at $\hat{\rho}$, the bound $\Phi_\Delta'' \ll \rho_\alpha^{-\alpha - 2} u_\alpha^{k - 1}$ on $I$ and \cref{eq:Remainder bounds} give
\[
\Phi(\rho_n) = \Phi_\Delta(\hat{\rho}) + O\left( \rho_\alpha^{-\alpha - 2} u_\alpha^{k - 1} (\rho_n - \hat{\rho})^2 + \rho_\alpha^{\varepsilon} \right) = \Phi_\Delta(\hat{\rho}) + o(1),
\]
where $\rho_\alpha^{-\alpha - 2} u_\alpha^{k - 1} (\rho_n - \hat{\rho})^2 \ll \rho_\alpha^{\alpha + 2\varepsilon} u_\alpha^{1 - k}$. In the same way, 
  $$
  \Psi''(\rho_n) = \Phi_\Delta''(\rho_n) + O(\rho_\alpha^{\varepsilon - 2}) = \Phi_\Delta''(\hat{\rho}) (1 + o(1)),$$ since $\Phi_\Delta''(\rho_n) - \Phi_\Delta''(\hat{\rho}) \ll \rho_\alpha^{-\alpha - 3} u_\alpha^{k - 1} |\rho_n - \hat{\rho}|,
  $
  and $\Phi_\Delta''(\hat{\rho}) \gg \rho_\alpha^{-\alpha - 2} u_\alpha^{k - 1}$. Therefore \cref{eq:Saddle form} becomes \cref{eq:Critical form}.
\end{proof}

\noindent
It remains to expand $\Phi_\Delta(\hat{\rho})$ and $\Phi_\Delta''(\hat{\rho})$ around $\rho_\alpha$. We first explain how the series $U(\eta, t)$ and the terms $T_m$ arise. Write $x = \rho_\alpha (1 + \eta)$. By \Cref{thm:Formal expansion} below, the Taylor series of $\Phi_\Delta(\rho_\alpha (1 + \eta)) - \Phi_\Delta(\rho_\alpha)$ in $\eta$ is
\[
\rho_\alpha^{-\alpha} \sum_{j \geq 2} \frac{q_j}{j!} (-\eta)^{j} + \sum_{j \geq 1} \frac{b_j}{j!} (-\eta)^{j} = U(\eta, 1).
\]
The first sum comes from $nx + x^{-\alpha} Q_0(\log x^{-1})$. It has no linear term, because $\rho_\alpha$ solves \cref{eq:Dominant saddle point equation}, and since $q_2 > 0$ it has a minimum at $\eta = 0$. The second sum comes from $\Delta$. By \cref{eq:Delta bound}, $b_j \ll \rho_\alpha^{-\sigma_1} u_\alpha^{h_1}$, while $q_j \ll u_\alpha^{k - 1}$ and $q_2 \asymp u_\alpha^{k - 1}$ by \cref{eq:Leading coefficient} and \Cref{thm:Dominant saddle point location}. So the second sum is smaller than the first by a factor of order $\rho_\alpha^{\alpha - \sigma_1} u_\alpha^{h_1 + 1 - k}$, which tends to $0$.

The variable $t$ in $U(\eta, t)$ measures this perturbation. For $t = 0$ the critical point of $U(\cdot, t)$ is $\eta = 0$, and as $t$ grows it moves to $\eta(t)$. For $t = 1$ we recover $\Phi_\Delta$, whose minimum on $I$ is at $\hat{\rho}$, so the quantity we need, $\Phi_\Delta(\hat{\rho}) - \Phi_\Delta(\rho_\alpha)$, is formally $U(\eta(1), 1) = \sum_{m \geq 2} T_m$. Since $U(\eta, t)$ depends on the $b_j$ only through the products $t b_j$, the coefficient $T_m$ of $t^m$ in $U(\eta(t), t)$ is a homogeneous polynomial of degree $m$ in the $b_j$. Thus $T_m$ collects the contributions of order $m$ in $\Delta$, and by \Cref{thm:Value expansion} we have $T_m \ll \rho_\alpha^{-\alpha} u_\alpha^{k - 1} (\rho_\alpha^{\alpha - \sigma_1} u_\alpha^{h_1 + 1 - k})^{m}$. By the choice of $N$, the terms with $m > N$ are then negligible.

\begin{lemma}
\label{thm:Formal expansion}
Let $q_j$, $b_j$, $U(\eta, t)$, $\eta(t)$ and $T_m$ be as in \Cref{thm:Exponent expansion}. Then
\begin{equation}
\label{eq:Phi Delta coefficients}
\rho_\alpha \Phi_\Delta'(\rho_\alpha) = -b_1, \qquad \rho_\alpha^{j} \Phi_\Delta^{(j)}(\rho_\alpha) = (-1)^{j} \big( \rho_\alpha^{-\alpha} q_j + b_j \big), \quad j \geq 2,
\end{equation}
so that the Taylor series of $\Phi_\Delta(\rho_\alpha (1 + \eta)) - \Phi_\Delta(\rho_\alpha)$ in $\eta$ is $U(\eta, 1)$. Moreover, $U(\eta(t), t)$ has no terms of degree less than $2$ in $t$, $m T_m$ is the coefficient of $t^{m - 1}$ in $\sum_{j \geq 1} \frac{b_j}{j!} (-\eta(t))^{j}$, and $T_2$, $T_3$ and $T_4$ are given by the formulas in \Cref{thm:Exponent expansion}.
\end{lemma}

\begin{proof}
Since $(-1)^{j} \rho_\alpha^{j} \Psi_\alpha^{(j)}(\rho_\alpha) = \rho_\alpha^{-\alpha} q_j$ and $\Psi_\alpha'(\rho_\alpha) = -n$, we have
\[
\rho_\alpha \Phi_\Delta'(\rho_\alpha) = \rho_\alpha \Delta'(\rho_\alpha) = -b_1, \qquad \rho_\alpha^{j} \Phi_\Delta^{(j)}(\rho_\alpha) = (-1)^{j} \big( \rho_\alpha^{-\alpha} q_j + b_j \big), \quad j \geq 2.
\]
Hence the Taylor series of $\Phi_\Delta(\rho_\alpha (1 + \eta)) - \Phi_\Delta(\rho_\alpha)$ in $\eta$ is $U(\eta, 1)$. The terms $T_m$ are computed most easily as follows. Since $\frac{\partial U}{\partial \eta}(\eta(t), t) = 0$, the chain rule gives
\begin{equation}
\label{eq:Envelope}
\frac{d}{dt} U(\eta(t), t) = \frac{\partial U}{\partial t}(\eta(t), t) = \sum_{j \geq 1} \frac{b_j}{j!} \big( -\eta(t) \big)^{j}.
\end{equation}
As $\eta(0) = 0$, the right-hand side has no constant term. Hence $U(\eta(t), t)$ starts with $t^2$, and $m T_m$ is the coefficient of $t^{m - 1}$ on the right-hand side, so that $T_m$ needs $\eta(t)$ only up to $t^{m - 1}$. To find $\eta(t)$, we write
\[
\frac{\partial U}{\partial \eta}(\eta, t) = \rho_\alpha^{-\alpha} \left( q_2 \eta - \frac{q_3}{2} \eta^2 + \cdots \right) + t \left( -b_1 + b_2 \eta - \frac{b_3}{2} \eta^2 + \cdots \right).
\]
Setting $\eta = \eta(t)$ and comparing the coefficients of $t$ and of $t^2$ with $0$, we find
\[
\eta(t) = \frac{\rho_\alpha^{\alpha} b_1}{q_2}\, t + \frac{\rho_\alpha^{2\alpha} b_1 (q_3 b_1 - 2 q_2 b_2)}{2 q_2^3}\, t^2 + \cdots
\]
Inserting this into the right-hand side of \cref{eq:Envelope}, which begins with $-b_1 \eta(t) + \frac{1}{2} b_2 \eta(t)^2$, we obtain
\[
2 T_2 = -\frac{\rho_\alpha^{\alpha} b_1^2}{q_2}, \qquad 3 T_3 = -\frac{\rho_\alpha^{2\alpha} b_1^2 (q_3 b_1 - 2 q_2 b_2)}{2 q_2^3} + \frac{b_2}{2} \cdot \frac{\rho_\alpha^{2\alpha} b_1^2}{q_2^2} = \frac{\rho_\alpha^{2\alpha} b_1^2 (3 q_2 b_2 - q_3 b_1)}{2 q_2^3}
\]
which are the formulas for $T_2$ and $T_3$ in \Cref{thm:Exponent expansion}. The formula for $T_4$ is obtained in the same way from the coefficient of $t^3$ in $\eta(t)$.
\end{proof}

\noindent
The following lemma makes the formal computation rigorous, up to the order $N$.

\begin{lemma}
\label{thm:Value expansion}
For $2 \leq m \leq N$ we have $T_m \ll \rho_\alpha^{-\alpha} u_\alpha^{k - 1} (\rho_\alpha^{\alpha - \sigma_1} u_\alpha^{h_1 + 1 - k})^m$, and
\begin{equation}
\label{eq:Value expansion}
\Phi_\Delta(\hat{\rho}) = \Phi_\Delta(\rho_\alpha) + T_2 + \cdots + T_N + O\left( \rho_\alpha^{(N + 1)(\alpha - \sigma_1) - \alpha} u_\alpha^{(N + 1) h_1 - N(k - 1)} \right).
\end{equation}
Moreover,
\begin{equation}
\label{eq:Critical point expansion}
\hat{\rho} - \rho_\alpha = \frac{\rho_\alpha^{1 + \alpha} b_1}{q_2} + O\left( \rho_\alpha^{1 + 2(\alpha - \sigma_1)} u_\alpha^{2(h_1 + 1 - k)} \right).
\end{equation}
\end{lemma}

\begin{proof}
We normalize the coefficients so that they are bounded. Put $a_j := q_j u_\alpha^{1 - k}$ and $\beta_j := b_j \rho_\alpha^{\sigma_1} u_\alpha^{-h_1}$, so that $\rho_\alpha^{\alpha} u_\alpha^{1 - k} b_j = \rho_\alpha^{\alpha - \sigma_1} u_\alpha^{h_1 + 1 - k} \beta_j$. By \Cref{thm:Dominant saddle point location}, \Cref{thm:Critical point} and \cref{eq:Delta bound}, the numbers $a_j$ and $\beta_j$ are bounded, and $a_2$ is bounded below by a positive constant. By \cref{eq:Phi Delta coefficients} and Taylor's formula at $\rho_\alpha$, up to order $N + 1$, we have, for $|\eta| \leq \frac{1}{2}$,
\[
\frac{\Phi_\Delta(\rho_\alpha (1 + \eta)) - \Phi_\Delta(\rho_\alpha)}{\rho_\alpha^{-\alpha} u_\alpha^{k - 1}} = \overline{U}\big( \eta, \rho_\alpha^{\alpha - \sigma_1} u_\alpha^{h_1 + 1 - k} \big) + O(|\eta|^{N + 2}),
\]
where
\[
\overline{U}(\eta, t) := \sum_{j = 2}^{N + 1} \frac{a_j}{j!} (-\eta)^{j} + t \sum_{j = 1}^{N + 1} \frac{\beta_j}{j!} (-\eta)^{j},
\]
and the derivative of the remainder with respect to $\eta$ is $O(|\eta|^{N + 1})$, because $\rho_\alpha^{N + 2} \Phi_\Delta^{(N + 2)} \ll \rho_\alpha^{-\alpha} u_\alpha^{k - 1}$ on $I$. The series $U(\eta, t)$ and $\rho_\alpha^{-\alpha} u_\alpha^{k - 1} \overline{U}(\eta, \rho_\alpha^{\alpha - \sigma_1} u_\alpha^{h_1 + 1 - k} t)$ differ only in terms of degree at least $N + 2$ in $\eta$. Since $\eta(t)$ has no constant term, such terms change $\eta(t)$ only in terms of order $t^{N + 1}$, and, as $\frac{\partial U}{\partial \eta}(\eta(t), t) = 0$, they change $U(\eta(t), t)$ only in terms of order $t^{N + 2}$. In particular they do not affect the coefficients of $t^m$ with $m \leq N$ in $U(\eta(t), t)$. Let $\bar{\eta}(t) = \sum_{m \geq 1} e_m t^m$ be the formal solution of $\frac{\partial \overline{U}}{\partial \eta}(\bar{\eta}(t), t) = 0$. It exists and is unique, because the coefficient of $t^m$ in $\frac{\partial \overline{U}}{\partial \eta}(\bar{\eta}(t), t)$ is $a_2 e_m$ plus a polynomial in $e_1, \ldots, e_{m - 1}$, the $a_j$ and the $\beta_j$. Write $\overline{U}(\bar{\eta}(t), t) = \sum_{m \geq 2} \overline{T}_m t^m$. Replacing $t$ by $\rho_\alpha^{\alpha - \sigma_1} u_\alpha^{h_1 + 1 - k} t$, we get $T_m = \rho_\alpha^{-\alpha} u_\alpha^{k - 1} \overline{T}_m (\rho_\alpha^{\alpha - \sigma_1} u_\alpha^{h_1 + 1 - k})^m$ for $2 \leq m \leq N$, and $e_m$ and $\overline{T}_m$ are polynomials in the $a_j$, the $\beta_j$ and $1/a_2$, hence bounded. In particular $T_m \ll \rho_\alpha^{-\alpha} u_\alpha^{k - 1} (\rho_\alpha^{\alpha - \sigma_1} u_\alpha^{h_1 + 1 - k})^m$ for $2 \leq m \leq N$.

Let $\eta_N(t) := \sum_{m = 1}^{N} e_m t^m$ and $\hat{\eta} := \hat{\rho}/\rho_\alpha - 1$, so that $\hat{\eta} \ll \rho_\alpha^{\alpha - \sigma_1} u_\alpha^{h_1 + 1 - k}$. The derivative of the left-hand side above with respect to $\eta$ vanishes at $\hat{\eta}$, as $\Phi_\Delta'(\hat{\rho}) = 0$. It is $\ll (\rho_\alpha^{\alpha - \sigma_1} u_\alpha^{h_1 + 1 - k})^{N + 1}$ at $\eta_N(\rho_\alpha^{\alpha - \sigma_1} u_\alpha^{h_1 + 1 - k})$, because $\frac{\partial \overline{U}}{\partial \eta}(\eta_N(t), t)$ is a polynomial in $t$ with bounded coefficients and without terms of degree at most $N$, and $\eta_N(t) \ll t$. The second derivative is $\rho_\alpha^{\alpha + 2} u_\alpha^{1 - k} \Phi_\Delta''(\rho_\alpha (1 + \eta)) \geq \frac{c_5}{2}$, and it is bounded above. Hence, by the mean value theorem,
\[
\hat{\eta} - \eta_N\big( \rho_\alpha^{\alpha - \sigma_1} u_\alpha^{h_1 + 1 - k} \big) \ll \big( \rho_\alpha^{\alpha - \sigma_1} u_\alpha^{h_1 + 1 - k} \big)^{N + 1},
\]
and by Taylor's formula at $\hat{\eta}$, the values of the left-hand side at these two points differ by $O\big( (\rho_\alpha^{\alpha - \sigma_1} u_\alpha^{h_1 + 1 - k})^{2N + 2} \big)$. Moreover, $\overline{U}(\eta_N(t), t) = \sum_{m = 2}^{N} \overline{T}_m t^m + O(t^{N + 1})$ for $0 \leq t \leq 1$, because $\overline{U}(\eta_N(t), t)$ is a polynomial in $t$ with bounded coefficients which agrees with $\overline{U}(\bar{\eta}(t), t)$ up to terms of order $t^{2N + 2}$, as $\frac{\partial \overline{U}}{\partial \eta}$ vanishes at $\bar{\eta}(t)$ and $\bar{\eta}(t) - \eta_N(t) = O(t^{N + 1})$. Taking $t = \rho_\alpha^{\alpha - \sigma_1} u_\alpha^{h_1 + 1 - k}$, adding the remainder $O(|\eta|^{N + 2})$ at $\eta = \eta_N(t)$, and multiplying by $\rho_\alpha^{-\alpha} u_\alpha^{k - 1}$, we obtain \cref{eq:Value expansion}.

Finally, $e_1 = \beta_1/a_2$ and $e_1 \rho_\alpha^{\alpha - \sigma_1} u_\alpha^{h_1 + 1 - k} = \rho_\alpha^{\alpha} b_1/q_2$. Since $\hat{\rho} - \rho_\alpha = \rho_\alpha \hat{\eta}$, and $\hat{\eta}$ differs from $\eta_N(\rho_\alpha^{\alpha - \sigma_1} u_\alpha^{h_1 + 1 - k})$, and hence from $e_1 \rho_\alpha^{\alpha - \sigma_1} u_\alpha^{h_1 + 1 - k}$, by $O\big( \rho_\alpha^{2(\alpha - \sigma_1)} u_\alpha^{2(h_1 + 1 - k)} \big)$, we get $\hat{\rho} - \rho_\alpha = \rho_\alpha^{1 + \alpha} b_1/q_2 + O\big( \rho_\alpha^{1 + 2(\alpha - \sigma_1)} u_\alpha^{2(h_1 + 1 - k)} \big)$, which is \cref{eq:Critical point expansion}.
\end{proof}

\noindent
We can now prove \Cref{thm:Exponent expansion}. By the choice of $N$, the error term in \cref{eq:Value expansion} tends to $0$. Since $\Phi_\Delta''(\rho_\alpha) = \rho_\alpha^{-\alpha - 2} (q_2 + \rho_\alpha^{\alpha} b_2)$ and $\Phi_\Delta'''(\rho_\alpha) = -\rho_\alpha^{-\alpha - 3} (q_3 + \rho_\alpha^{\alpha} b_3)$ by \cref{eq:Phi Delta coefficients}, Taylor's formula and \cref{eq:Critical point expansion} give
\[
\begin{aligned}
\Phi_\Delta''(\hat{\rho}) &= \Phi_\Delta''(\rho_\alpha) + \Phi_\Delta'''(\rho_\alpha) (\hat{\rho} - \rho_\alpha) + O\big( \rho_\alpha^{\alpha - 2\sigma_1 - 2} u_\alpha^{2h_1 + 1 - k} \big) \\
&= \rho_\alpha^{-\alpha - 2} \left( q_2 + \rho_\alpha^{\alpha} (b_2 - b_1 q_3/q_2) \right) + O\big( \rho_\alpha^{\alpha - 2\sigma_1 - 2} u_\alpha^{2h_1 + 1 - k} \big),
\end{aligned}
\]
where we used $\Phi_\Delta^{(4)} \ll \rho_\alpha^{-\alpha - 4} u_\alpha^{k - 1}$ on $I$ and $\rho_\alpha^{2\alpha} b_1 b_3/q_2 \ll \rho_\alpha^{2(\alpha - \sigma_1)} u_\alpha^{2h_1 + 1 - k}$. As $q_2 \asymp u_\alpha^{k - 1}$ and $\rho_\alpha^{\alpha} b_j \ll \rho_\alpha^{\alpha - \sigma_1} u_\alpha^{h_1} = o(u_\alpha^{k - 1})$, both $\rho_\alpha^{-\alpha - 2} (q_2 + \rho_\alpha^{\alpha} (b_2 - b_1 q_3/q_2))$ and $\rho_\alpha^{-\alpha - 2} q_2$ are equal to $\Phi_\Delta''(\hat{\rho}) (1 + o(1))$. Inserting the expansions of $\Phi_\Delta(\hat{\rho})$ and $\Phi_\Delta''(\hat{\rho})$ into \cref{eq:Critical form}, and using $\Phi_\Delta(\rho_\alpha) = n\rho_\alpha + \rho_\alpha^{-\alpha} Q_0(u_\alpha) + \Delta(\rho_\alpha)$, we obtain \cref{eq:Expanded formula}. \qed

\subsection{Proof of \texorpdfstring{\Cref{thm:Explicit saddle point}}{Theorem 1.1}}

\noindent
\Cref{thm:Explicit saddle point} is \Cref{thm:Dominant saddle point}(i) with $\rho_\alpha$ computed explicitly. For $k \leq 2$ the polynomials $Q_\ell$ have degree at most $1$, so \cref{eq:Dominant saddle point equation} can be solved in closed form, by a power of $n$ for $k = 1$ and by the Lambert function for $k = 2$. For $k \geq 3$, the same substitution as for $k = 2$ turns it into the equation of a generalized Lambert function. It then remains to substitute the solution into \cref{eq:Dominant formula}.

In view of \Cref{thm:Dominant saddle point}(i), we only have to compute $\rho_\alpha$, $n\rho_\alpha + \rho_\alpha^{-\alpha} Q_0(u_\alpha)$ and $\rho_\alpha^{-\alpha - 2} Q_2(u_\alpha)$, where $u_\alpha = \log \rho_\alpha^{-1}$. By \cref{eq:Leading coefficient} with $\ell = 0$, the leading coefficient of $Q_0$ is $\frac{F(\alpha, k - 1)}{(k - 1)!} = \frac{C}{\alpha(1 + \alpha)}$, where $C > 0$ by \Cref{thm:Positive leading coefficient}. If $h_0 = 0$, then by \Cref{thm:Res at 0} the residue of $M(z, \rho_\alpha)$ at $z = 0$ is $D(0) u_\alpha + D'(0)$.

Let $k = 1$. Then the polynomial $Q_0(s) = F(\alpha, 0) = \frac{C}{\alpha(1 + \alpha)}$ is constant, so that the relation $Q_{\ell + 1} = (\alpha + \ell) Q_\ell + Q_\ell'$ gives
\[
Q_1 = \frac{C}{1 + \alpha}, \qquad Q_2 = \alpha(1 + \alpha) Q_0 = C.
\]
Thus \cref{eq:Dominant saddle point equation} reads $\frac{C}{1 + \alpha} \rho^{-\alpha - 1} = n$. Its left-hand side is strictly decreasing on $(0, \infty)$, so it has the unique solution
\[
\rho_\alpha = \left( \frac{C}{(1 + \alpha) n} \right)^{\frac{1}{1 + \alpha}}, \qquad u_\alpha = \frac{1}{1 + \alpha} \log \frac{(1 + \alpha) n}{C}.
\]
Using $\rho_\alpha^{-\alpha - 1} = \frac{(1 + \alpha) n}{C}$, we get
\[
\rho_\alpha^{-\alpha} Q_0(u_\alpha) = \rho_\alpha \cdot \frac{(1 + \alpha) n}{C} \cdot \frac{C}{\alpha(1 + \alpha)} = \frac{n\rho_\alpha}{\alpha},
\]
and therefore,
\[
n\rho_\alpha + \rho_\alpha^{-\alpha} Q_0(u_\alpha) = \left( 1 + \frac{1}{\alpha} \right) n\rho_\alpha = \left( 1 + \frac{1}{\alpha} \right) \left( \frac{C}{1 + \alpha} \right)^{\frac{1}{1 + \alpha}} n^{\frac{\alpha}{1 + \alpha}}.
\]
Moreover,
\[
\rho_\alpha^{-\alpha - 2} Q_2(u_\alpha) = C \rho_\alpha^{-\alpha - 2} = C \left( \frac{(1 + \alpha) n}{C} \right)^{\frac{2 + \alpha}{1 + \alpha}}.
\]
Inserting these expressions into \cref{eq:Dominant formula} gives \cref{eq:Simple pole formula}. If $h_0 = 0$, the residue at $z = 0$ is $D(0) u_\alpha + D'(0) = \frac{D(0)}{1 + \alpha} \log \frac{(1 + \alpha) n}{C} + D'(0)$. This proves (i).

Now let $k = 2$. The leading coefficient of $Q_0$ is $F(\alpha, 1) = \frac{C}{\alpha(1 + \alpha)} > 0$, so $\kappa$ is well defined. By the definition of $\kappa$, $F(\alpha, 0) = F(\alpha, 1)\left( \kappa - \frac{1}{\alpha} \right)$, and hence
\[
Q_0(s) = \frac{C}{\alpha(1 + \alpha)} \left( s + \kappa - \frac{1}{\alpha} \right), \qquad Q_0' = \frac{C}{\alpha(1 + \alpha)}, \qquad Q_0'' = 0,
\]
and therefore, as $Q_1 = \alpha Q_0 + Q_0'$ and $Q_2 = \alpha(1 + \alpha) Q_0 + (1 + 2\alpha) Q_0' + Q_0''$,
\[
Q_1(s) = \frac{C}{1 + \alpha} (s + \kappa), \qquad Q_2(s) = C (s + \kappa) + \frac{(1 + 2\alpha) C}{\alpha(1 + \alpha)} - \frac{C}{\alpha} = C \left( s + \kappa + \frac{1}{1 + \alpha} \right).
\]
With $v := (1 + \alpha)(\log \rho^{-1} + \kappa)$, we have $\rho^{-\alpha - 1} = e^{v - (1 + \alpha)\kappa}$ and $Q_1(\log \rho^{-1}) = \frac{C v}{(1 + \alpha)^2}$. So \cref{eq:Dominant saddle point equation} reads $\frac{C v}{(1 + \alpha)^2} e^{v - (1 + \alpha)\kappa} = n$, that is
\begin{equation}
\label{eq:Lambert equation}
v e^{v} = \frac{(1 + \alpha)^2 e^{(1 + \alpha)\kappa}}{C}\, n.
\end{equation}
Note that the map $\rho \mapsto v$ is a decreasing bijection from $(0, \infty)$ onto $\mathbb{R}$. The right-hand side of \cref{eq:Lambert equation} is positive, and $v \mapsto v e^{v}$ is negative on $(-\infty, 0)$ and strictly increasing from $0$ to $\infty$ on $[0, \infty)$. Hence \cref{eq:Lambert equation} has exactly one real solution, namely $v = v_n$. So \cref{eq:Dominant saddle point equation} has a unique positive solution $\rho_\alpha$, and
\[
u_\alpha + \kappa = \frac{v_n}{1 + \alpha}, \qquad \rho_\alpha^{-\alpha - 1} = \frac{(1 + \alpha) n}{C (u_\alpha + \kappa)} = \frac{(1 + \alpha)^2 n}{C v_n}.
\]
In particular, $\rho_\alpha = \left( \frac{C v_n}{(1 + \alpha)^2 n} \right)^{\frac{1}{1 + \alpha}}$, which is the argument of $\Delta$ in \cref{eq:Double pole formula}. Next,
\[
n\rho_\alpha^{1 + \alpha} = \frac{C v_n}{(1 + \alpha)^2}, \qquad Q_0(u_\alpha) = \frac{C}{\alpha(1 + \alpha)} \left( \frac{v_n}{1 + \alpha} - \frac{1}{\alpha} \right),
\]
so that, as $n\rho_\alpha + \rho_\alpha^{-\alpha} Q_0(u_\alpha) = \rho_\alpha^{-\alpha} \big( n\rho_\alpha^{1 + \alpha} + Q_0(u_\alpha) \big)$, we have 
\[
n\rho_\alpha + \rho_\alpha^{-\alpha} Q_0(u_\alpha) = \rho_\alpha^{-\alpha} \left( \frac{C v_n}{(1 + \alpha)^2} \left( 1 + \frac{1}{\alpha} \right) - \frac{C}{\alpha^2 (1 + \alpha)} \right) = \rho_\alpha^{-\alpha} \frac{C}{\alpha(1 + \alpha)} \left( v_n - \frac{1}{\alpha} \right).
\]
Since $\rho_\alpha^{-\alpha} = \left( \frac{(1 + \alpha)^2 n}{C v_n} \right)^{\frac{\alpha}{1 + \alpha}}$ and $\frac{2\alpha}{1 + \alpha} - 1 = \frac{\alpha - 1}{1 + \alpha}$, it gives 
\[
\frac{(1 + \alpha)^{\frac{\alpha - 1}{1 + \alpha}} C^{\frac{1}{1 + \alpha}}}{\alpha} \left( v_n - \frac{1}{\alpha} \right) \left( \frac{n}{v_n} \right)^{\frac{\alpha}{1 + \alpha}}.
\]
Finally, $Q_2(u_\alpha) = C \left( \frac{v_n}{1 + \alpha} + \frac{1}{1 + \alpha} \right) = \frac{C (v_n + 1)}{1 + \alpha}$, so that
\[
\rho_\alpha^{-\alpha - 2} Q_2(u_\alpha) = \frac{C (v_n + 1)}{1 + \alpha} \left( \frac{(1 + \alpha)^2 n}{C v_n} \right)^{\frac{2 + \alpha}{1 + \alpha}}.
\]
Substituting everything into \cref{eq:Dominant formula}, we obtain \cref{eq:Double pole formula}; for $h_0 = 0$ the residue at $z = 0$ is $D(0) u_\alpha + D'(0) = \frac{D(0)}{1 + \alpha} \log \frac{(1 + \alpha)^2 n}{C v_n} + D'(0)$. This proves (ii).

Finally let $k \geq 3$. By \cref{eq:Leading coefficient} with $\ell = 1$, $Q_1$ is a polynomial of degree $k - 1$ with real coefficients and leading coefficient $\frac{C}{1 + \alpha}$. So $Q_1(w/(1 + \alpha))$ has leading coefficient $C/(1 + \alpha)^{k}$, which gives the factorization in (iii), and the non-real $w_i$ occur in complex conjugate pairs. Hence $\mathcal{H}_k$ is real on $\mathbb{R}$. With $w := (1 + \alpha) \log \rho^{-1}$ we have $\rho^{-\alpha - 1} = e^{w}$ and $Q_1(\log \rho^{-1}) = Q_1(w/(1 + \alpha))$, so \cref{eq:Dominant saddle point equation} reads
\[
\frac{C}{(1 + \alpha)^{k}}\, e^{w} \prod_{i = 1}^{k - 1} (w - w_i) = n,
\]
which is \cref{eq:Generalized Lambert equation}. The map $\rho \mapsto w$ is a decreasing bijection from $(0, \infty)$ onto $\mathbb{R}$, so the smallest positive solution $\rho_\alpha$ of \cref{eq:Dominant saddle point equation} corresponds to the largest real solution of \cref{eq:Generalized Lambert equation}. This solution exists, because $\mathcal{H}_k$ is continuous on $\mathbb{R}$ and tends to $0$ as $w \to -\infty$ and to $\infty$ as $w \to \infty$, and because $\mathcal{H}_k - y_n$ is a non-constant real analytic function which does not vanish outside a bounded interval, so that it has only finitely many real zeros. Thus $\rho_\alpha = e^{-w_n/(1 + \alpha)}$ and $u_\alpha = w_n/(1 + \alpha)$.

Since $n\rho_\alpha^{1 + \alpha} = Q_1(u_\alpha)$ by \cref{eq:Dominant saddle point equation}, we have
\[
n\rho_\alpha + \rho_\alpha^{-\alpha} Q_0(u_\alpha) = \rho_\alpha^{-\alpha} \big( Q_1(u_\alpha) + Q_0(u_\alpha) \big) = e^{\frac{\alpha w_n}{1 + \alpha}} (Q_0 + Q_1)\left( \frac{w_n}{1 + \alpha} \right)
\]
and $\rho_\alpha^{-\alpha - 2} Q_2(u_\alpha) = e^{\frac{(2 + \alpha) w_n}{1 + \alpha}} Q_2\big( \frac{w_n}{1 + \alpha} \big)$. Inserting these expressions into \cref{eq:Dominant formula} gives \cref{eq:Higher pole formula}, and for $h_0 = 0$ the residue at $z = 0$ is $D(0) u_\alpha + D'(0) = \frac{D(0)}{1 + \alpha} w_n + D'(0)$.

It remains to expand $w_n$. Since $\mathcal{H}_k(w_n) = y_n \to \infty$ and $\mathcal{H}_k$ is bounded on every interval $(-\infty, A]$, we have $w_n \to \infty$. For large real $w$, $\prod_{i = 1}^{k - 1} (w - w_i) = w^{k - 1} (1 + O(1/w))$, so taking logarithms in \cref{eq:Generalized Lambert equation} gives
\[
w_n + (k - 1) \log w_n + O\left( \frac{1}{w_n} \right) = \log y_n.
\]
Hence $w_n \sim \log y_n$, so that $\log w_n = \log \log y_n + O(\log \log y_n/\log y_n)$, and inserting this into the last relation gives the expansion of $w_n$ in (iii). \qed

\subsection{Proof of \texorpdfstring{\Cref{thm:Log asymptotics}}{Corollary 2.4}}

\noindent
Here only the leading term of $\log p_f(n)$ is needed. So we work with $\rho_n$ itself, keep only the leading terms of $\Psi$ and $\Psi'$, and finally express $n\rho_n$ through $n$.

Let $\rho = \rho_n$ and $u := \log \rho^{-1}$. By \Cref{thm:Main theorem}, we have 
\[
\log p_f(n) = n\rho + \Psi(\rho) - \frac{1}{2} \log\big( 2\pi \Psi''(\rho) \big) + o(1).
\]
By \cref{eq:Derivative asymptotics} with $\ell = 0$ and $\ell = 1$, since $\Gamma(\alpha)/\Gamma(\alpha + 2) = 1/(\alpha(1 + \alpha))$ and $\Gamma(\alpha + 1)/\Gamma(\alpha + 2) = 1/(1 + \alpha)$,
\begin{equation}
\label{eq:Psi and n}
\begin{aligned}
\Psi(\rho) &= \frac{C}{\alpha(1 + \alpha)}\, \rho^{-\alpha} u^{k - 1} \left( 1 + O\left( \frac{1}{u} \right) \right), \\
n = -\Psi'(\rho) &= \frac{C}{1 + \alpha}\, \rho^{-\alpha - 1} u^{k - 1} \left( 1 + O\left( \frac{1}{u} \right) \right).
\end{aligned}
\end{equation}
Dividing the first relation by the second, we get $\Psi(\rho) = \frac{n\rho}{\alpha} \left( 1 + O(1/u) \right)$, and $n\rho + \Psi(\rho) = \left( 1 + \frac{1}{\alpha} \right) n\rho \left( 1 + O(1/u) \right)$. By \cref{eq:Derivative upper bound}, $\Psi''(\rho) \asymp \rho^{-\alpha - 2} u^{k - 1}$, so $\log \Psi''(\rho) = (\alpha + 2) u + O(\log u) = O(u)$, while $n\rho \asymp \rho^{-\alpha} u^{k - 1} = e^{\alpha u} u^{k - 1}$ by \cref{eq:Psi and n}. Hence $u^2 = o(n\rho)$, so the term $\frac{1}{2} \log\big( 2\pi \Psi''(\rho) \big)$, which is $O(u)$, is also $O(n\rho/u)$. Therefore
\[
\log p_f(n) = \left( 1 + \frac{1}{\alpha} \right) n\rho \left( 1 + O\left( \frac{1}{u} \right) \right)
\]
To express $n\rho$ in terms of $n$, recall from the proof of \Cref{thm:Size of rho} that $u \sim \frac{\log n}{1 + \alpha}$. The second relation in \cref{eq:Psi and n} gives $\rho^{1 + \alpha} = \frac{C}{1 + \alpha} \cdot \frac{u^{k - 1}}{n} \left( 1 + o(1) \right)$, and therefore
\[
n\rho = \left( \frac{C}{1 + \alpha} \right)^{\frac{1}{1 + \alpha}} n^{\frac{\alpha}{1 + \alpha}} u^{\frac{k - 1}{1 + \alpha}} \left( 1 + o(1) \right) = \left( \frac{C}{(1 + \alpha)^k} \right)^{\frac{1}{1 + \alpha}} n^{\frac{\alpha}{1 + \alpha}} (\log n)^{\frac{k - 1}{1 + \alpha}} \left( 1 + o(1) \right)
\]
where in the last step we used $u^{\frac{k - 1}{1 + \alpha}} = \left( \frac{\log n}{1 + \alpha} \right)^{\frac{k - 1}{1 + \alpha}} (1 + o(1))$. This proves \Cref{thm:Log asymptotics}. \qed

\noindent
Taking logarithms in the second relation of \cref{eq:Psi and n}, we also get $(1 + \alpha) u = \log n - (k - 1) \log u + O(1)$, so $u = \frac{\log n}{1 + \alpha} \left( 1 + O\left( \frac{\log \log n}{\log n} \right) \right)$, and the relative error in \Cref{thm:Log asymptotics} is $O(\log \log n/\log n)$.

\subsection{Proof of \texorpdfstring{\Cref{thm:Accuracy of saddle point}}{Proposition 2.5}}
\label{sec:Proof of accuracy}

\noindent
In the variable $u = \log \rho^{-1}$, the exponent $\Phi(e^{-u})$ has a minimum at $\log \rho_n^{-1}$, with second derivative of order $n\rho_n$ (\Cref{thm:Convexity} below). So moving this point by $o\big( (n\rho_n)^{-1/2} \big)$ changes the exponent by $o(1)$, while a larger displacement changes it by an amount that does not tend to $0$; this gives (i). For (ii) we show that, for $k \geq 2$, $\log \rho_n^{-1}$ has an expansion in $1/\log n$ and $\log \log n/\log n$ with infinitely many nonzero terms, so that a finite sum as in (ii) differs from $\log \rho_n^{-1}$ by at least a power of $1/\log n$.

Let $\Phi(x) = nx + \Psi(x)$ be as in the proof of \Cref{thm:Dominant saddle point}, and put
\[
G(u) := \Phi(e^{-u}) = n e^{-u} + \Psi(e^{-u}), \qquad u_n := \log \rho_n^{-1}, \qquad \tilde{u} := \log \tilde{\rho}^{-1}, \qquad \eta_n := \tilde{u} - u_n.
\]
Then $G(u_n) = \Phi(\rho_n)$, $G(\tilde{u}) = \Phi(\tilde{\rho})$, and $G'(u_n) = -\rho_n \Phi'(\rho_n) = 0$. Since $\tilde{\rho} \to 0$, we have $\tilde{u} \to \infty$. By \Cref{thm:Size of rho},
\begin{equation}
\label{eq:Size of n rho}
n\rho_n \asymp \rho_n^{-\alpha} u_n^{k - 1} \asymp n^{\frac{\alpha}{1 + \alpha}} (\log n)^{\frac{k - 1}{1 + \alpha}}
\end{equation}
so the condition in (i) reads $\eta_n = o\big( (n\rho_n)^{-1/2} \big)$.

We first explain where this condition comes from. Replacing $\rho_n$ by $\tilde{\rho}$ in the formula of \Cref{thm:Main theorem} changes two things: the numerator $e^{\Phi(\rho_n)} = e^{G(u_n)}$ becomes $e^{G(\tilde{u})}$, and the denominator changes by the factor $\big( \Psi''(\tilde{\rho})/\Psi''(\rho_n) \big)^{1/2}$. Since $G'(u_n) = 0$, the first change is governed by $G''(u_n)$. With $x = e^{-u}$ we have $G''(u) = x\Phi'(x) + x^2 \Psi''(x)$, and since $\Phi'(\rho_n) = 0$, \cref{eq:Derivative upper bound} and \cref{eq:Size of n rho} give
\[
G''(u_n) = \rho_n^2 \Psi''(\rho_n) \asymp \rho_n^{-\alpha} u_n^{k - 1} \asymp n\rho_n.
\]
So, if $\eta_n$ is small, $G(\tilde{u}) - G(u_n) \approx \frac{1}{2} G''(u_n) \eta_n^2$ is of order $n\rho_n \eta_n^2$, while the denominator changes only by a factor $1 + O(|\eta_n|)$, which tends to $1$. Hence the formula remains correct exactly when $n\rho_n \eta_n^2 \to 0$: the number $(n\rho_n)^{-1/2}$ is the distance from $u_n$ over which $G$ changes by a bounded amount. For example, for ordinary partitions we have $\alpha = k = 1$ and $n\rho_n \asymp n^{1/2}$, so $\log \tilde{\rho}^{-1}$ must approximate $\log \rho_n^{-1}$ with error $o(n^{-1/4})$. Since (A2) holds for every $\sigma_0 \in (0, 1)$ (\Cref{sec:Examples}), \cref{eq:Saddle points close} shows that $\rho_\alpha = \pi/\sqrt{6n}$ does this, with error $O(n^{-\sigma})$ for every $\sigma < 1/2$. \Cref{thm:Convexity}, \Cref{thm:Denominator ratio} and \Cref{thm:Error bounds} make this argument rigorous, also when $\eta_n$ is not small; in that case the convexity of $G$ gives a lower bound which grows linearly in $|\eta_n|$.

\begin{lemma}
\label{thm:Convexity}
There are $u_0 > 0$ and $c_8, c_9 > 0$ such that $G''(u) > 0$ for $u \geq u_0$, and, for large $n$,
\begin{equation}
\label{eq:Xi second derivative}
c_8\, n\rho_n \leq G''(u) \leq c_9\, n\rho_n, \qquad |u - u_n| \leq 1.
\end{equation}
Consequently, for large $n$ and $\tilde{u} \geq u_0$,
\begin{equation}
\label{eq:Exponent gap}
\Phi(\tilde{\rho}) - \Phi(\rho_n) \geq \frac{c_8}{2}\, n\rho_n \min\left( \eta_n^2, |\eta_n| \right), \qquad \Phi(\tilde{\rho}) - \Phi(\rho_n) \leq \frac{c_9}{2}\, n\rho_n \eta_n^2 \quad \text{if } |\eta_n| \leq 1.
\end{equation}
\end{lemma}

\begin{proof}
Differentiating $G(u) = n e^{-u} + \Psi(e^{-u})$ twice, we get, with $x = e^{-u}$, $G''(u) = nx + x\Psi'(x) + x^2 \Psi''(x)$. By \cref{eq:Derivative expansion} with $\ell = 1, 2$,
\[
x\Psi'(x) + x^2 \Psi''(x) = x^{-\alpha} (Q_2 - Q_1)(u) + O(x^{-\sigma_0})
\]
where $Q_2 - Q_1 = \alpha Q_1 + Q_1'$ has degree $k - 1$ and, by \cref{eq:Leading coefficient}, the positive leading coefficient $\alpha C/(1 + \alpha)$. Since $x^{-\alpha} = e^{\alpha u}$ and the error term is $O(e^{\sigma_0 u})$, there is $u_0 \geq \log 2$ such that $x\Psi'(x) + x^2\Psi''(x) \asymp e^{\alpha u} u^{k - 1}$ for $u \geq u_0$. In particular $G''(u) > 0$ there, as also $nx > 0$. For $|u - u_n| \leq 1$ and large $n$ we have $u \geq u_0$ and $x/\rho_n \in [e^{-1}, e]$, hence $nx \asymp n\rho_n$, and $e^{\alpha u} u^{k - 1} \asymp \rho_n^{-\alpha} u_n^{k - 1} \asymp n\rho_n$ by \cref{eq:Size of n rho}, which gives \cref{eq:Xi second derivative}.

If $|\eta_n| \leq 1$, Taylor's formula at $u_n$, where $G'(u_n) = 0$, gives $G(\tilde{u}) - G(u_n) = \frac{1}{2} G''(\xi) \eta_n^2$ with $|\xi - u_n| \leq 1$, and \cref{eq:Exponent gap} follows from \cref{eq:Xi second derivative}. If $\eta_n \geq 1$, then, since $G$ is convex on $[u_0, \infty)$ and hence lies above its tangent at $u_n + 1$,
\[
G(\tilde{u}) - G(u_n) \geq G(u_n + 1) - G(u_n) + G'(u_n + 1)(\eta_n - 1) \geq \frac{c_8}{2} n\rho_n + c_8\, n\rho_n (\eta_n - 1) \geq \frac{c_8}{2} n\rho_n \eta_n
\]
since $G(u_n + 1) - G(u_n) \geq \frac{c_8}{2} n\rho_n$ by the first case and $G'(u_n + 1) = \int_{u_n}^{u_n + 1} G'' \geq c_8 n\rho_n$. The case $\eta_n \leq -1$ is the same, using $u_n - 1$ and $u_0 \leq \tilde{u}$.
\end{proof}

\noindent
We also need to compare $\Psi''(\tilde{\rho})$ with $\Psi''(\rho_n)$.

\begin{lemma}
\label{thm:Denominator ratio}
For large $n$ we have
\begin{equation}
\label{eq:Denominator ratio small}
\left| \log \Psi''(\tilde{\rho}) - \log \Psi''(\rho_n) \right| \ll |\eta_n|, \qquad |\eta_n| \leq 1.
\end{equation}
and, without any condition on $\eta_n$,
\begin{equation}
\label{eq:Denominator ratio large}
\log \Psi''(\tilde{\rho}) - \log \Psi''(\rho_n) \leq (\alpha + 2) |\eta_n| + (k - 1) \log(1 + |\eta_n|) + O(1).
\end{equation}
\end{lemma}

\begin{proof}
Let $n$ be so large that $\tilde{\rho} \leq 1/2$ and $\rho_n \leq 1/(2e)$; this is possible since $\tilde{\rho} \to 0$ and $\rho_n \to 0$. If $|\eta_n| \leq 1$, then $\tilde{\rho}$ and $\rho_n$ lie in $[\rho_n/e, e\rho_n]$, where $|\Psi'''|/\Psi'' \ll \rho_n^{-1}$ by \cref{eq:Derivative upper bound}. Hence, by the mean value theorem applied to $\log \Psi''$, and since $\tilde{\rho} = \rho_n e^{-\eta_n}$,
\[
\left| \log \Psi''(\tilde{\rho}) - \log \Psi''(\rho_n) \right| \ll \rho_n^{-1} |\tilde{\rho} - \rho_n| = |e^{-\eta_n} - 1| \ll |\eta_n|
\]
In general, \cref{eq:Derivative upper bound} gives $\Psi''(x) \asymp x^{-\alpha - 2} (\log x^{-1})^{k - 1}$ for $0 < x \leq 1/2$, so that 
  $$
  \log \Psi''(\tilde{\rho}) - \log \Psi''(\rho_n) = (\alpha + 2) \eta_n + (k - 1) \log(\tilde{u}/u_n) + O(1).
  $$ 
  Since $u_n \geq 1$ and $\tilde{u} \leq u_n + |\eta_n| \leq u_n (1 + |\eta_n|)$, this gives \cref{eq:Denominator ratio large}.
\end{proof}

\noindent
By \Cref{thm:Main theorem}, $p_f(n) \sim e^{\Phi(\rho_n)}/\sqrt{2\pi \Psi''(\rho_n)}$, and $e^{\tilde{\rho} n} P(\tilde{\rho}) = e^{\Phi(\tilde{\rho})}$. Hence
\begin{equation}
\label{eq:Ratio to saddle form}
\frac{e^{\tilde{\rho} n} P(\tilde{\rho})}{p_f(n) \sqrt{2 \pi \Psi''(\tilde{\rho})}} = \left( 1 + o(1) \right) e^{\mathcal{E}_n},\quad \mathcal{E}_n := \Phi(\tilde{\rho}) - \Phi(\rho_n) - \frac{1}{2} \left( \log \Psi''(\tilde{\rho}) - \log \Psi''(\rho_n) \right).
\end{equation}
Thus, up to a factor $1 + o(1)$, $e^{\mathcal{E}_n}$ is the factor by which the right-hand side of \cref{eq:Perturbed saddle point} differs from $p_f(n)$, and the two terms of $\mathcal{E}_n$ are the two changes described above. So \cref{eq:Perturbed saddle point} holds if and only if $\mathcal{E}_n \to 0$, and if $\mathcal{E}_n \to \infty$, the right-hand side is too large by an unbounded factor. The following lemma bounds $\mathcal{E}_n$ in terms of $\eta_n$.

\begin{lemma}
\label{thm:Error bounds}
There is $c_{10} > 0$ such that, for large $n$,
\begin{equation}
\label{eq:E large chi}
\mathcal{E}_n \geq \frac{c_8}{4} n\rho_n - O(1), \qquad |\eta_n| \geq 1.
\end{equation}
and
\begin{equation}
\label{eq:E small chi}
\frac{c_8}{2} n\rho_n \eta_n^2 - c_{10} |\eta_n| \leq \mathcal{E}_n \leq \frac{c_9}{2} n\rho_n \eta_n^2 + c_{10} |\eta_n|, \qquad |\eta_n| \leq 1.
\end{equation}
\end{lemma}

\begin{proof}
Since $\tilde{\rho} \to 0$, we have $\tilde{u} \geq u_0$ for large $n$, so that \cref{eq:Exponent gap} applies. If $|\eta_n| \geq 1$, then \cref{eq:Exponent gap} and \cref{eq:Denominator ratio large} give, since $\log(1 + |\eta_n|) \leq |\eta_n|$,
\[
\mathcal{E}_n \geq \left( \frac{c_8}{2} n\rho_n - \frac{\alpha + k + 1}{2} \right) |\eta_n| - O(1) \geq \frac{c_8}{4} n\rho_n - O(1)
\]
since $|\eta_n| \geq 1$ and $n\rho_n \to \infty$. If $|\eta_n| \leq 1$, then \cref{eq:E small chi} follows from \cref{eq:Exponent gap} and \cref{eq:Denominator ratio small}, where $2 c_{10}$ is the implied constant in \cref{eq:Denominator ratio small}.
\end{proof}

\noindent
We can now prove (i). Suppose first that $\eta_n = o\big( (n\rho_n)^{-1/2} \big)$. Then $\eta_n \to 0$, as $n\rho_n \to \infty$, and by \cref{eq:E small chi}, $|\mathcal{E}_n| \ll n\rho_n \eta_n^2 + |\eta_n| = o(1)$. By \cref{eq:Ratio to saddle form}, \cref{eq:Perturbed saddle point} holds.

Conversely, assume \cref{eq:Perturbed saddle point}. Then $\mathcal{E}_n \to 0$ by \cref{eq:Ratio to saddle form}. By \cref{eq:E large chi}, $|\eta_n| < 1$ for large $n$. If $|\eta_n| \geq 4c_{10}/(c_8 n\rho_n)$, then $c_{10} |\eta_n| \leq \frac{c_8}{4} n\rho_n \eta_n^2$, and \cref{eq:E small chi} gives $\frac{c_8}{4} n\rho_n \eta_n^2 \leq \mathcal{E}_n$. Otherwise $n\rho_n \eta_n^2 < 16 c_{10}^2/(c_8^2 n\rho_n)$. In both cases $n\rho_n \eta_n^2 \to 0$, which is the condition in (i) by \cref{eq:Size of n rho}.

We turn to (ii). For $k \geq 2$ the expansion of $\log \rho_n^{-1}$ in powers of $1/\log n$ and $\log \log n/\log n$ does not terminate, as for the Lambert function, and a finite sum of monomials leaves a term of size at least a power of $1/\log n$, which is much larger than $(n\rho_n)^{-1/2}$. The expansion is given in \Cref{thm:Structure of saddle point}, and \Cref{thm:Monomials} gives the size of a finite sum of monomials. For the first lemma we write $Q_1(u) = \frac{C}{1 + \alpha} u^{k - 1} \widehat{Q}_1(1/u)$, where $\widehat{Q}_1$ is a polynomial with $\widehat{Q}_1(0) = 1$, by \cref{eq:Leading coefficient}.

\begin{lemma}
\label{thm:Structure of saddle point}
There are $\kappa_0 > 0$ and a function $h(X, Y)$, analytic in a neighbourhood of $(0, 0)$, such that
\begin{equation}
\label{eq:Structure of saddle point}
\log \rho_n^{-1} = \frac{\log n}{1 + \alpha} - \frac{k - 1}{1 + \alpha} \log \log n + h\left( \frac{1}{\log n}, \frac{\log \log n}{\log n} \right) + O(n^{-\kappa_0})
\end{equation}
and
\begin{equation}
\label{eq:Theta on axis}
h(0, Y) = \frac{k \log(1 + \alpha) - \log C}{1 + \alpha} - \frac{k - 1}{1 + \alpha} \log\left( 1 - (k - 1) Y \right).
\end{equation}
In particular, if $h(X, Y) = \sum_{i, m \geq 0} e_{i, m} X^i Y^m$, then
\begin{equation}
\label{eq:Theta coefficients}
e_{0, m} = \frac{(k - 1)^{m + 1}}{(1 + \alpha)\, m} \neq 0, \qquad m \geq 1.
\end{equation}
\end{lemma}

\begin{proof}
By \Cref{thm:Psi decomposition}, \Cref{thm:Dominant saddle point location} and \Cref{thm:Size of rho}, the saddle-point equation $n = -\Psi'(\rho_n)$ gives
\[
\Omega(u_n) = \log n + O(n^{-\kappa_0}), \qquad \Omega(u) := (1 + \alpha) u + \log Q_1(u), \qquad \kappa_0 := \frac{\alpha - \sigma_0}{2(1 + \alpha)}.
\]
With $X = 1/\log n$, $Y = \log \log n/\log n$ and $u = \frac{\log n}{1 + \alpha} - \frac{k - 1}{1 + \alpha} \log \log n + Z$, a direct computation gives $\Omega(u) = \log n + H(Z, X, Y)$, where
\[
\begin{aligned}
H(Z, X, Y) := {} & (1 + \alpha) Z + \log C - k \log(1 + \alpha) + (k - 1) \log\left( 1 - (k - 1) Y + (1 + \alpha) Z X \right) \\
& + \log \widehat{Q}_1\left( \frac{(1 + \alpha) X}{1 - (k - 1) Y + (1 + \alpha) Z X} \right).
\end{aligned}
\]
This function is analytic near $(Z_0, 0, 0)$, where $Z_0 := (k \log(1 + \alpha) - \log C)/(1 + \alpha)$, and $H(Z_0, 0, 0) = 0$, $\frac{\partial H}{\partial Z}(Z_0, 0, 0) = 1 + \alpha$. The analytic implicit function theorem gives $h$ with $h(0, 0) = Z_0$ and $H(h(X, Y), X, Y) = 0$, and since $\widehat{Q}_1(0) = 1$, the equation $H(Z, 0, Y) = 0$ gives \cref{eq:Theta on axis}. Moreover, $\Omega'(u) \to 1 + \alpha$ as $u \to \infty$, so that $\Omega' \geq 1$ for large $u$, and the mean value theorem gives \cref{eq:Structure of saddle point}. Finally, expanding $-\log(1 - w) = \sum_{m \geq 1} w^m/m$ in \cref{eq:Theta on axis} gives \cref{eq:Theta coefficients}.
\end{proof}

\noindent
We call the functions $\mathsf{m}(n) = n^{a_0} (\log n)^{a_1} \cdots (\log_r n)^{a_r}$ with real $a_i$ {\it monomials}, and $(a_0, \ldots, a_r)$ their exponent vectors; taking into account zero coordinates, we may lexicographically compare exponent vectors of the same length, so if the exponent vectors of two monomials differ, the one whose vector is larger in the first coordinate where they differ is of larger order.

\begin{lemma}
\label{thm:Monomials}
  A finite linear combination of distinct monomials with nonzero coefficients, 
  is asymptotic to its term with the largest monomial. 
\end{lemma}

\begin{proof}
Let $\mathsf{m}$ and $\mathsf{m}'$ have the exponent vectors $(a_j)$ and $(a'_j)$, and let $i$ be the first coordinate with $a_i \neq a'_i$, say $a_i > a'_i$. With $\log_0 n := n$ we have
\[
\frac{\mathsf{m}(n)}{\mathsf{m}'(n)} = (\log_i n)^{a_i - a'_i} \prod_{j > i} (\log_j n)^{a_j - a'_j}.
\]
Since $\log_j n = o\big( (\log_i n)^{\varepsilon} \big)$ for every $j > i$ and every $\varepsilon > 0$, this quotient tends to infinity. In a finite linear combination, every term other than the one with the largest monomial is therefore $o$ of that term.
\end{proof}

\noindent
We can now prove (ii). Let $k \geq 2$. We may assume $r \geq 2$ and that the monomials in $L(n)$ are distinct. Choose $m^{*} \geq 1$ such that $\mathsf{m}^{*} := Y^{m^{*}} = (\log n)^{-m^{*}} (\log \log n)^{m^{*}}$ does not occur in $L(n)$. The monomials $X^i Y^m = (\log n)^{-i - m} (\log \log n)^{m}$ are distinct for distinct $(i, m)$, and since $h$ is analytic and $0 < X \leq Y$ for large $n$, \cref{eq:Structure of saddle point} gives
\[
\log \rho_n^{-1} = \frac{\log n}{1 + \alpha} - \frac{k - 1}{1 + \alpha} \log \log n + \sum_{i + m \leq m^{*}} e_{i, m} X^i Y^m + O\big( Y^{m^{*} + 1} \big)
\]
So $L(n) - \log \rho_n^{-1}$ is a finite linear combination of distinct monomials plus $o(\mathsf{m}^{*}(n))$, in which $\mathsf{m}^{*}$ has the coefficient $-e_{0, m^{*}} \neq 0$ by \cref{eq:Theta coefficients}. By \Cref{thm:Monomials}, its largest monomial with a nonzero coefficient is at least $\mathsf{m}^{*}$ for large $n$, and we obtain $|L(n) - \log \rho_n^{-1}| \gg \mathsf{m}^{*}(n) \geq (\log n)^{-m^{*} - 1}$ for large $n$. This is the first claim, with $c = m^{*} + 1$.

Now $|\eta_n| = |L(n) - \log \rho_n^{-1}| \geq (\log n)^{-c}$. If $|\eta_n| \geq 1$, then $\mathcal{E}_n \geq \frac{c_8}{4} n\rho_n - O(1)$ by \cref{eq:E large chi}. If $|\eta_n| < 1$, then $\mathcal{E}_n \geq \frac{c_8}{2} n\rho_n (\log n)^{-2c} - c_{10}$ by \cref{eq:E small chi}. In both cases $\mathcal{E}_n \to \infty$, because $n\rho_n$ grows like a positive power of $n$ by \cref{eq:Size of n rho}. By \cref{eq:Ratio to saddle form}, the right-hand side of \cref{eq:Perturbed saddle point} divided by $p_f(n)$ tends to infinity. \qed

\subsection{Motivation for \texorpdfstring{\Cref{thm:No closed form}}{Conjecture 2.6}}
\label{sec:Motivation}

\noindent
Let $k \geq 2$. By \Cref{thm:Accuracy of saddle point}(i) and \cref{eq:Size of n rho}, a function $L(n)$ gives an asymptotic formula in \cref{eq:Perturbed saddle point}, with $\tilde{\rho} = e^{-L(n)}$, if and only if
\[
L(n) - \log \rho_n^{-1} = o\big( (n\rho_n)^{-1/2} \big), \qquad (n\rho_n)^{-1/2} \asymp n^{-\frac{\alpha}{2(1 + \alpha)}} (\log n)^{-\frac{k - 1}{2(1 + \alpha)}}.
\]
Since $(n\rho_n)^{-1/2}$ is a negative power of $n$ up to a power of $\log n$, such $L(n)$ satisfies $L(n) - \log \rho_n^{-1} = O(n^{-c})$ for some $c > 0$. We show that \Cref{thm:No closed form} excludes this for every LE-function $L$.

Let $G(m)$ be the largest real solution $u$ of $e^{(1 + \alpha) u} Q_1(u) = m$, so that $G(n) = \log \rho_\alpha^{-1} = w_n/(1 + \alpha)$, where $w_n$ is the largest real solution of $e^{w} P(w) = \lambda n$ with $P(w) = \prod_{i = 1}^{k - 1} (w - w_i)$ and $\lambda = (1 + \alpha)^{k}/C$, as in \Cref{thm:Explicit saddle point}(iii). 
The function $G$ is the inverse of $u \mapsto e^{(1 + \alpha) u} Q_1(u)$, whose derivative is $e^{(1 + \alpha) u} Q_2(u)$. Hence every derivative $G^{(j)}(m)$ is a rational function of $u = G(m)$ and $e^{(1 + \alpha) u}$, and $G^{(j)}(m) \ll_j m^{-j}$. Now let $L$ be an LE-function with $L(n) - \log \rho_n^{-1} = O(n^{-c})$ for some $c > 0$, and put
\[
M(n) := e^{(1 + \alpha) L(n)} Q_1(L(n))
\]
Then $M$ is an LE-function, and $G(M(n)) = L(n)$ for large $n$. By \Cref{thm:Psi decomposition} and \cref{eq:Remainder bounds}, the saddle-point equation $\Psi'(\rho_n) + n = 0$ reads $\rho_n^{-\alpha - 1} Q_1(\log \rho_n^{-1}) = n + \Delta'(\rho_n) + O(\rho_n^{\varepsilon - 1})$, and $\Delta'(\rho_n) \ll \rho_n^{-\sigma_0 - 1}$. Since $\rho_n^{-\alpha - 1} \asymp n (\log n)^{-(k - 1)}$, it follows that $M(n) - n \ll n^{1 - c'}$ for some $c' > 0$. Taylor's formula for $G$ at $M(n)$ gives, for every fixed $J$,
\[
G(n) = \sum_{j = 0}^{J - 1} \frac{G^{(j)}(M(n))}{j!} \big( n - M(n) \big)^{j} + O\big( n^{-c' J} \big).
\]
The sum is an LE-function of $n$, since $G^{(j)}(M(n))$ is a rational function of $L(n)$ and $e^{(1 + \alpha) L(n)}$. Multiplying by $1 + \alpha$, we obtain an LE-function which approximates $w_n$ with an error $O(n^{-c' J})$, and this contradicts \Cref{thm:No closed form}. Note that this argument does not need any condition on $\sigma_0$, and that it does not need to know the poles in $\Delta$: their effect is contained in the difference $n - M(n)$.

Next we explain why \Cref{thm:Accuracy of saddle point}(ii) does not settle the question. \Cref{thm:Structure of saddle point} expresses $\log \rho_n^{-1}$, up to $O(n^{-\kappa_0})$, by a convergent double series in $1/\log n$ and $\log \log n/\log n$, and by \cref{eq:Theta coefficients} infinitely many of its coefficients are nonzero. A truncation of this series has an error of the size of a power of $\log \log n/\log n$, far larger than the negative power of $n$ that is needed; this is \Cref{thm:Accuracy of saddle point}(ii). An LE-function $L$ with $\log \rho_n^{-1} - L(n) = O(n^{-c})$ would have to reproduce the whole series, and the remainder as well. Finite sums of monomials cannot do this, but LE-functions can reproduce infinitely many terms: by \cref{eq:Theta on axis}, all terms $e_{0, m} Y^{m}$ add up to the LE-function $-\frac{k - 1}{1 + \alpha} \log(1 - (k - 1) Y)$. The conjecture says that the whole series cannot be reproduced in this way.

Of course $\log \rho_n^{-1}$ can be computed by iteration, but the number of steps has to grow with $n$. Let $\Omega(u) = (1 + \alpha) u + \log Q_1(u)$ be as in the proof of \Cref{thm:Structure of saddle point}, and let $u^{*}$ be the solution of $\Omega(u^{*}) = \log n$, so that $u^{*} = \log \rho_n^{-1} + O(n^{-\kappa_0})$. Newton's method
\[
L_0(n) := \frac{\log n}{1 + \alpha}, \qquad L_{j + 1}(n) := L_j(n) - \frac{\Omega(L_j(n)) - \log n}{\Omega'(L_j(n))}
\]
produces LE-functions $L_j(n)$. Since $\Omega'(u) = 1 + \alpha + O(1/u)$ and $\Omega''(u) = -(k - 1)/u^2 + O(1/u^3)$, Taylor's formula gives $L_{j + 1} - u^{*} = \frac{\Omega''(\xi_j)}{2\Omega'(L_j)} (L_j - u^{*})^2$ for some $\xi_j \asymp \log n$, and $|L_0 - u^{*}| \asymp \log \log n$. By induction, for each fixed $j$,
\[
\left| L_j(n) - u^{*} \right| \asymp_j (\log \log n)^{2^j} (\log n)^{2 - 2^{j + 1}}.
\]
For each fixed $j$ this is only a negative power of $\log n$, so by \Cref{thm:Accuracy of saddle point}(i) the formula \cref{eq:Perturbed saddle point} fails for $\tilde{\rho} = e^{-L_j(n)}$. The required precision is reached only when $2^{j}$ is roughly of order $\log n/\log \log n$. Thus the number of steps must tend to infinity with $n$, although it grows only like $\log_2 \log n$. In practice a few steps suffice: for $f(n) = d(n)$, where $\alpha = 1$ and $k = 2$, and $n = 10^{200}$, four Newton steps give $|L_4(n) - u^{*}| \approx 5 \cdot 10^{-73}$, far below $n^{-1/4} = 10^{-50}$. The simpler iteration $L_{j + 1} = (\log n - \log Q_1(L_j))/(1 + \alpha)$ gains only a factor $\asymp 1/\log n$ per step. For $k = 1$, in contrast, $\Omega$ is linear, and $L_1(n) = u^{*}$ exactly, in accordance with \Cref{thm:Explicit saddle point}(i).



\section{Examples and applications}
\label{sec:Examples}
We now apply the results of \Cref{sec:Main result} and \Cref{sec:Dominant saddle point} to a number of concrete examples of weights, as also mentioned in the introduction. In each case we determine $D(s)$, the dominant pole $\alpha$ and its order $k$, the constant $C$, and the range of $\sigma_0$ for which (A2) holds; this range decides which part of \Cref{thm:Dominant saddle point} applies and, for part (iii), how many Newton steps are needed. \Cref{thm:Main theorem} and \Cref{thm:Log asymptotics} apply to all examples except the ordered factorizations in \Cref{sec:Ordered factorizations}, for which \cref{eq:Complex poles to the left} fails.

The examples are ordered by the result they need. In \Crefrange{sec:Partitions}{sec:Logarithm}, the only poles in $\Re(s) \geq -\varepsilon$ are $\alpha$ and possibly $0$, so \Cref{thm:Explicit saddle point} applies directly, and these examples show the three cases $k = 1$, $k = 2$ and $k \geq 3$. In \Crefrange{sec:Divisor of squares}{sec:Unitary}, a factor $\zeta(as + b)^{-1}$ produces complex poles from the zeros of $\zeta$, and their position decides between \Cref{thm:Explicit saddle point} and \Cref{thm:Dominant saddle point}(ii), sometimes only under the Riemann hypothesis. In \Crefrange{sec:Commuting}{sec:Representations}, some pole $p \neq \alpha$ satisfies $\Re(p) \geq \alpha/2$, and parts (ii) and (iii) of \Cref{thm:Dominant saddle point} and \Cref{thm:Exponent expansion} are needed. Finally, \Cref{sec:Ordered factorizations} shows where the method stops.

To apply the results to a given weight $f$, one proceeds as follows.
\begin{enumerate}
    \item[(1)] Compute $D(s)$, its meromorphic continuation to $\Re(s) \geq -\varepsilon$ and its poles, and check \cref{eq:Complex poles to the left}. The largest real pole is $\alpha$, its order is $k$, and $C$, $Q_0$, $Q_1$ and $Q_2$ are obtained from the principal part at $\alpha$ as in \Cref{sec:Main result}.
    \item[(2)] Find the range of $\sigma_0$ for which (A2) holds. Its lower end is at least the largest real part of a pole $p \neq \alpha$. For products of shifted zeta functions and their reciprocals, \Cref{thm:Zeta quotients} below gives (A2) on lines to the right of the poles, unconditionally or under the Riemann hypothesis.
    \item[(3)] Check (A3). Usually $f(n) \geq 1$ for all $n$, or $\Lambda \setminus p\mathbb{N}$ is infinite for every prime $p$.
    \item[(4)] Compute the residues of $\Gamma(z) \zeta(z + 1) D(z) \rho^{-z}$ at the other poles and at $z = 0$ by \Cref{thm:Res at 0} and \Cref{thm:Res at p}. Their sum is $\Delta(\rho)$, and finitely many poles with $\Re(p) < 0$ may be omitted.
    \item[(5)] Compute $\rho_\alpha$ by \Cref{thm:Explicit saddle point}: it is a power of $n$ for $k = 1$, given by the Lambert function for $k = 2$, and by a generalized Lambert function for $k \geq 3$.
    \item[(6)] Choose the formula according to $\sigma_0$: \Cref{thm:Dominant saddle point}(i) if $\sigma_0 < \alpha/2$, part (ii) if $\sigma_0 < 2\alpha/3$, and otherwise part (iii), with the least $J$ such that $2^{J} (\alpha - \sigma_0) > \alpha/2$, or \Cref{thm:Exponent expansion} for a formula in which only $\rho_\alpha$ occurs.
\end{enumerate}

\Cref{thm:Main theorem} and \Cref{thm:Log asymptotics} need only (1)--(3). \Cref{tab:Examples} summarizes the examples of this section.

The last column of \Cref{tab:Examples} gives the ratio of $p_f(n)$ to the right-hand side of the asymptotic formula of the corresponding subsection for $n = 3000$, for the parameters in brackets, and for $2^{\omega(n)}$, $d(n^2)$ and $d(n)^2$ in this order. For $d_3(n)$ this is the formula of \Cref{thm:Dominant saddle point}(i), and for $d(n^2)$ and $d(n)^2$ that of \Cref{thm:Dominant saddle point}(ii), with $\Delta$ as given there. The sums over $\varrho$ are taken over the first $40$ zeros in the upper half-plane and their conjugates. For the subrings of $\mathbb{Z}^3$ the two values are without and with the full denominator of \cref{eq:Expanded formula}, for $SU(3)$ the ratio is for Romik's formula, for $SU(3)^2$ for \cref{eq:Critical form}, and for $H(n)$ it is the ratio of $\log p_f(n)$ to the right-hand side of the logarithmic asymptotics. The larger deviations for $2^{\omega(n)}$, for the $r$-free parts and for $\sigma^{*}$ come from the omitted residues at the poles with $\Re(z) < 0$, which tend to $0$ only like a power of $\rho_\alpha$. With the residues at the poles in $-3 < \Re(z) < 0$ added to the exponent, the ratios become $0.997$, $1.001$ ($r = 2$), $1.000$ ($r = 3$) and $0.998$.

\begin{table}[tbp]
\centering
\footnotesize
\setlength{\tabcolsep}{3pt}
\renewcommand{\arraystretch}{1.45}
\begin{tabular}{@{}c >{\raggedright\arraybackslash}p{1.7cm} l c c l >{\raggedright\arraybackslash}p{1.75cm} >{\raggedright\arraybackslash}p{2.3cm} >{\raggedright\arraybackslash}p{1.3cm}@{}}
\toprule
Sec. & $f(n)$ & $D(s)$ & $\alpha$ & $k$ & $C$ & (A2) for $\sigma_0 \in$ & Result & Ratio at $n = 3000$ \\
\midrule
\ref{sec:Partitions} & $1$ & $\zeta(s)$ & $1$ & $1$ & $\pi^2/3$ & $(0, 1)$ & Cor.~\ref{thm:Explicit saddle point}(i) & $0.992$ \\
\ref{sec:Partitions} & $n$ & $\zeta(s - 1)$ & $2$ & $1$ & $6\zeta(3)$ & $(0, 2)$ & Cor.~\ref{thm:Explicit saddle point}(i) & $0.999$ \\
\ref{sec:Divisor functions} & $d_r(n)$ & $\zeta(s)^{r}$ & $1$ & $r$ & $\frac{\pi^2}{3 (r - 1)!}$ & $(0, 1)$ & Cor.~\ref{thm:Explicit saddle point}(ii), (iii) & $0.9995$ ($r = 2$), $0.999$ ($r = 3$) \\
\ref{sec:Logarithm} & $(\log n)^{j}$ & $(-1)^{j} \zeta^{(j)}(s)$ & $1$ & $j + 1$ & $\pi^2/3$ & $(0, 1)$ & Cor.~\ref{thm:Explicit saddle point}(ii), (iii) & $1.004$ ($j = 1$) \\
\ref{sec:Divisor of squares} & $2^{\omega(n)}$, $d(n^2)$, $d(n)^2$ & $\zeta(s)^{a}/\zeta(2s)$ & $1$ & $a$ & $\frac{2}{(a - 1)!}$ & $(\frac{1}{2}, 1)$; RH: $(\frac{1}{4}, 1)$ & Thm.~\ref{thm:Dominant saddle point}(ii); RH: (i) & $0.979$, $1.002$, $0.999$ \\
\ref{sec:r-free} & $r$-free indicator & $\zeta(s)/\zeta(rs)$ & $1$ & $1$ & $\frac{\pi^2}{3 \zeta(r)}$ & $(\frac{1}{r}, 1)$ & Cor.~\ref{thm:Explicit saddle point}(i) ($r \geq 3$); Thm.~\ref{thm:Dominant saddle point}(ii) ($r = 2$) & $1.080$ ($r = 2$), $1.512$ ($r = 3$) \\
\ref{sec:Jordan} & $J_r(n)$ & $\zeta(s - r)/\zeta(s)$ & $r + 1$ & $1$ & $\frac{(r + 2)!\, \zeta(r + 2)}{\zeta(r + 1)}$ & $(1, r + 1)$ & Cor.~\ref{thm:Explicit saddle point}(i) ($r \geq 2$); Thm.~\ref{thm:Dominant saddle point}(ii) ($r = 1$) & $0.9995$ ($r = 1$), $0.990$ ($r = 2$) \\
\ref{sec:Jordan} & $\psi_r(n)$ & $\frac{\zeta(s) \zeta(s - r)}{\zeta(2s)}$ & $r + 1$ & $1$ & $\frac{(r + 2)!\, \zeta(r + 1) \zeta(r + 2)}{\zeta(2r + 2)}$ & $(1, r + 1)$ & Cor.~\ref{thm:Explicit saddle point}(i) ($r \geq 2$); Thm.~\ref{thm:Dominant saddle point}(ii) ($r = 1$) & $0.991$ ($r = 1$), $1.005$ ($r = 2$) \\
\ref{sec:Unitary} & $\sigma^{*}(n)$ & $\frac{\zeta(s) \zeta(s - 1)}{\zeta(2s - 1)}$ & $2$ & $1$ & $\pi^2$ & $(1, 2)$; RH: $(\frac{3}{4}, 2)$ & Thm.~\ref{thm:Dominant saddle point}(ii); RH: Cor.~\ref{thm:Explicit saddle point}(i) & $0.864$ \\
\ref{sec:Commuting} & $a_n(\mathbb{Z}^d)$ & $\prod_{i = 0}^{d - 1} \zeta(s - i)$ & $d$ & $1$ & $(d + 1)! \prod_{j = 2}^{d + 1} \zeta(j)$ & $(d - 1, d)$ & Thm.~\ref{thm:Dominant saddle point}(ii) ($d = 2$); Thms.~\ref{thm:Dominant saddle point}(iii), \ref{thm:Exponent expansion} ($d \geq 3$) & $0.990$ ($d = 3$), $0.989$ ($d = 4$) \\
\ref{sec:Heisenberg} & $a_n(\mathcal{N})$ & see \cref{eq:Heisenberg zeta} & $2$ & $2$ & $\pi^4/12$ & $(\frac{3}{2}, 2)$ & Thms.~\ref{thm:Dominant saddle point}(iii), \ref{thm:Exponent expansion} & $0.9992$ \\
\ref{sec:Subrings} & subrings of $\mathbb{Z}^3$ & $\frac{\zeta(3s - 1) \zeta(s)^3}{\zeta(2s)^2}$ & $1$ & $3$ & $1$ & $(\frac{2}{3}, 1)$ & Thms.~\ref{thm:Dominant saddle point}(iii), \ref{thm:Exponent expansion} & $0.990$; $0.998$ \\
\ref{sec:Representations} & $r_{SU(3)}(n)$ & $2^{s} \omega(s)$ & $\frac{2}{3}$ & $1$ & $8.9599\ldots$ & $(\frac{1}{2}, \frac{2}{3})$ & Thms.~\ref{thm:Dominant saddle point}(iii), \ref{thm:Exponent expansion} & $0.86$ \\
\ref{sec:Representations} & $r_{SU(3)^r}(n)$ & $2^{rs} \omega(s)^{r}$ & $\frac{2}{3}$ & $r$ & $\frac{10\, \Gamma(\frac{2}{3}) \zeta(\frac{5}{3})}{9 (r - 1)!} R^{r}$ & $(\frac{1}{2}, \frac{2}{3})$ & Thms.~\ref{thm:Dominant saddle point}(iii), \ref{thm:Exponent expansion} & $0.998$ ($r = 2$) \\
\ref{sec:Ordered factorizations} & $H(n)$ & $1/(2 - \zeta(s))$ & $1.7286\ldots$ & $1$ & $3.0034\ldots$ & (A1) fails & $\log p_f(n)$ only & $0.978$ \\
\bottomrule
\end{tabular}
\caption{The summary of examples in \Cref{sec:Examples}. The column ``Result'' gives asymptotic formulas for $p_f(n)$ in terms of $\rho_\alpha$, some of which are under Riemann Hypothesis (RH). The last column is explained at the beginning of \Cref{sec:Examples}. In the row for $SU(3)^r$, $R$ is as in \Cref{sec:Representations}.}
\label{tab:Examples}
\end{table}

When $D(s)$ contains a factor $\zeta(as + b)^{-1}$, every non-trivial zero $\varrho$ of $\zeta$ gives a pole $p = (\varrho - b)/a$ of $D(s)$, and $\Delta(\rho)$ contains the residues of $\Gamma(z) \zeta(z + 1) D(z) \rho^{-z}$ at these poles. If $\varrho$ is a simple zero, this residue is
\[
\frac{\Gamma(p) \zeta(p + 1) \widetilde{D}(p)}{a \zeta'(\varrho)} \rho^{-p}, \qquad \widetilde{D}(s) := \zeta(as + b) D(s)
\]
In the formulas below we write these residues out, with $\rho$ expressed in terms of $n$. The sums $\sum_{\varrho}$ run over the distinct non-trivial zeros of $\zeta$, and, in accordance with \Cref{thm:Infinite B}, they are the limits as $j \to \infty$ of the sums over $|\Im(\varrho)| < a T_j$, with $T_j$ as in \Cref{thm:Zeta quotients}(ii). The terms are written for simple zeros; at a multiple zero, the term is replaced by the residue of $\Gamma(z) \zeta(z + 1) D(z) \rho^{-z}$ at $z = p$. As noted in \Cref{thm:Remark on poles}, the factor $\Gamma(p)$ makes these terms numerically small, but their modulus is of order $\rho^{-\Re(p)}$, so they cannot be left out when $\Re(p) > 0$.

Most of the Dirichlet series below are products of shifted zeta functions and their reciprocals, and for these (A2) can be checked with the following lemma.

\begin{lemma}
\label{thm:Zeta quotients}
Let $D(s) = \prod_{i = 1}^{l} \zeta(a_i s + b_i)^{e_i}$ with $a_i > 0$, $b_i \in \mathbb{R}$ and $e_i \in \mathbb{Z} \setminus \{0\}$.
\begin{enumerate}
    \item[(i)] Let $\sigma \in \mathbb{R}$ be such that the line $\Re(s) = \sigma$ contains no pole of $D(s)$, and $a_i \sigma + b_i \notin [0, 1]$ for every $i$ with $e_i < 0$. Then $D(\sigma + it) \ll_{\sigma} (1 + |t|)^{A}$ for some $A = A(\sigma)$, so \cref{eq:Exponential bound} holds on this line for every $a > 0$.
    \item[(ii)] For every $c > 0$ there is a sequence $T_j \to \infty$ such that the segments $\{\sigma \pm iT_j : -\varepsilon \leq \sigma \leq c\}$ contain no pole of $D(s)$, and $D(\sigma \pm iT_j) \ll e^{a T_j}$ uniformly for $-\varepsilon \leq \sigma \leq c$, for every $a > 0$.
    \item[(iii)] Under the Riemann hypothesis, (i) holds with $[0, 1]$ replaced by $[0, \frac{1}{2}]$.
\end{enumerate}
\end{lemma}

\begin{proof}
For (i), the factors with $e_i > 0$ grow at most polynomially in $|t|$ on vertical lines \cite[Chapter~V]{Titchmarsh}. For a factor with $e_i < 0$, put $u = a_i \sigma + b_i$. If $u > 1$, then $|\zeta(u + iv)| \geq \zeta(2u)/\zeta(u) > 0$. If $u < 0$, the functional equation gives $|\zeta(u + iv)| \asymp |v|^{1/2 - u} |\zeta(1 - u - iv)| \gg |v|^{1/2 - u}$ for $|v| \geq 1$. Hence $\zeta(u + iv)^{-1} \ll 1$ for $|v| \geq 1$, and the bounded range of $t$ is handled by continuity.

For (ii), consider a factor with $e_i < 0$. For $-1 \leq u \leq 2$ and $v \geq 2$, we have $\frac{\zeta'}{\zeta}(u + iv) = \sum_{|v - \gamma| < 1} (u + iv - \varrho)^{-1} + O(\log v)$, where $\varrho = \beta + i\gamma$ runs over the zeros of $\zeta$ \cite[Chapter~15]{Davenport}. There are $O(\log V)$ zeros with $|\gamma - V| < 2$. Excluding intervals of length $\delta/\log V$ around their ordinates, we obtain a set of $v \in [V, V + 1]$ of measure at least $1 - O(\delta)$ on which $\frac{\zeta'}{\zeta}(u + iv) \ll_\delta \log^2 V$ for $-1 \leq u \leq 2$. Integrating from $u = 2$ gives $|\zeta(u + iv)|^{-1} \leq \exp(O(\log^2 V))$ there. Since there are finitely many factors, we can choose $T_j$ so that this holds for all of them simultaneously, with $v = a_i T_j$. For $u > 2$ or $u < -1$ we use (i). The factors with $e_i > 0$ grow polynomially, so $D(\sigma \pm iT_j) \ll \exp(O(\log^2 T_j)) \ll e^{a T_j}$.

For (iii), under the Riemann hypothesis $\zeta(u + iv)^{-1} \ll_{u, \eta} |v|^{\eta}$ for every fixed $u > \frac{1}{2}$ and $\eta > 0$ \cite[Chapter~XIV]{Titchmarsh}, and the line $\Re(s) = \sigma$ contains no pole of the factor.
\end{proof}

\noindent
In all examples below, (A3) holds because $f(n) \geq 1$ for all $n$, or, for $(\log n)^j$ and the $r$-free indicator, because $\Lambda \setminus p\mathbb{N}$ is infinite for every prime $p$. All Dirichlet series are meromorphic in $\mathbb{C}$, and $\varepsilon$ is chosen so small that no pole lies on the line $\Re(s) = -\varepsilon$.

In several examples below only part (iii) of \Cref{thm:Dominant saddle point} applies. There we use \Cref{thm:Exponent expansion} to write the formulas in terms of $\rho_\alpha$, and hence of $n$.

\subsection{Partitions and plane partitions}
\label{sec:Partitions}
For $f = 1$ we have $D(s) = \zeta(s)$, $\alpha = 1$, $k = 1$, $C = \pi^2/3$, and the only other pole of $\Gamma(z) \zeta(z + 1) D(z)$ is $z = 0$, with $D(0) = -\frac{1}{2}$ and $D'(0) = -\frac{1}{2} \log 2\pi$. (A2) holds for every $\sigma_0 \in (0, 1)$, so \Cref{thm:Explicit saddle point}(i) applies. Since $\Delta(\rho) = -\frac{1}{2} \log \rho^{-1} - \frac{1}{2} \log 2\pi$, \cref{eq:Simple pole formula} becomes the Hardy--Ramanujan formula \cite{HardyRamanujan}
\[
p(n) \sim \frac{1}{4\sqrt{3}\, n} e^{\pi \sqrt{2n/3}}, \quad n \to \infty. 
\]
For $f(n) = n$, the coefficient $p_f(n)$ is the number of plane partitions of $n$, and $D(s) = \zeta(s - 1)$. Here $\alpha = 2$, $k = 1$, $C = 6\zeta(3)$, $D(0) = \zeta(-1) = -\frac{1}{12}$ and $D'(0) = \zeta'(-1)$, and (A2) holds for every $\sigma_0 \in (0, 2)$. \Cref{thm:Explicit saddle point}(i) gives Wright's formula \cite{Wright1931}
\[
p_f(n) \sim \frac{\zeta(3)^{7/36} e^{\zeta'(-1)}}{2^{11/36} \sqrt{3\pi}\, n^{25/36}} \exp\left\{ 3 \zeta(3)^{1/3} \left( \frac{n}{2} \right)^{2/3} \right\}, \quad n \to \infty. 
\]

\subsection{Generalized divisor functions and representations of $SU(2)^r$}
\label{sec:Divisor functions}
Let $f(n) = d_r(n)$, the number of ways to write $n$ as an ordered product of $r \geq 2$ positive integers. These weights also have representation theoretic meaning: $d_r(n)$ counts the number of irreducible representations of $SU(2)^r$ of dimension $n$. 
By complete reducibility, $p_f(n)$ is the number of $n$-dimensional representations of $SU(2)^r$ up to isomorphism. Then $D(s) = \zeta(s)^r$, with the single pole $\alpha = 1$ of order $k = r$ and $A(1, r) = 1$, so that $C = \pi^2/(3(r - 1)!)$. The only other pole of $\Gamma(z) \zeta(z + 1) D(z)$ is $z = 0$, with $h_0 = 0$, $D(0) = (-\frac{1}{2})^r$ and $D'(0) = r (-\frac{1}{2})^r \log 2\pi$. (A2) holds for every $\sigma_0 \in (0, 1)$ by \Cref{thm:Zeta quotients}. Hence \Cref{thm:Dominant saddle point}(i) applies with
\[
\Delta(\rho) = \left( -\tfrac{1}{2} \right)^r \left( \log \rho^{-1} + r \log 2\pi \right)
\]
and \Cref{thm:Log asymptotics} gives $\log p_f(n) \sim \frac{2\pi}{\sqrt{3 \cdot 2^r (r - 1)!}} \sqrt{n}\, (\log n)^{\frac{r - 1}{2}}$.

For $r = 2$, $f(n) = d(n)$, we have $F(1, 1) = \frac{\pi^2}{6}$ and $F(1, 0) = \frac{\gamma \pi^2}{6} + \zeta'(2)$, so $\kappa = 1 + \gamma + \frac{\zeta'(2)}{\zeta(2)}$. With $v_n = W(12 e^{2\kappa} n/\pi^2)$, \Cref{thm:Explicit saddle point}(ii) gives
\[
p_f(n) \sim \frac{\exp\left\{ \frac{\pi}{\sqrt{3}} (v_n - 1) \sqrt{\frac{n}{v_n}} + \frac{1}{8} \log \frac{12 n}{\pi^2 v_n} + \frac{1}{2} \log 2\pi \right\}}{\sqrt{2\pi \cdot \frac{\pi^2 (v_n + 1)}{6} \left( \frac{12 n}{\pi^2 v_n} \right)^{3/2}}}.
\]
Since $SU(2) \times SU(2) \cong \mathrm{Spin}(4)$, we obtain an explicit asymptotic formula for the number of $n$-dimensional representations of $\mathrm{Spin}(4)$. For $r \geq 3$, $\rho_\alpha$ is given by the generalized Lambert function of \Cref{thm:Explicit saddle point}(iii), and by \Cref{thm:No closed form} we do not expect a closed form for it; see \Cref{sec:Motivation}.

\subsection{Powers of the logarithm}
\label{sec:Logarithm}
Let $f(n) = (\log n)^j$ with $j \geq 1$. Then $D(s) = (-1)^j \zeta^{(j)}(s)$. Since $\zeta(s) - \frac{1}{s - 1}$ is entire, the only pole is $\alpha = 1$, of order $k = j + 1$, with $A(1, j + 1) = j!$ and $A(1, i) = 0$ for $1 \leq i \leq j$. Hence $C = \pi^2/3$, and $h_0 = 0$, $D(0) = (-1)^j \zeta^{(j)}(0)$, $D'(0) = (-1)^j \zeta^{(j + 1)}(0)$. By Cauchy's estimates, $\zeta^{(j)}$ satisfies the same polynomial bounds on vertical lines as $\zeta$, so (A2) holds for every $\sigma_0 \in (0, 1)$, as in \Cref{thm:Zeta quotients}. Hence \Cref{thm:Dominant saddle point}(i) applies, with $\Delta(\rho) = D(0) \log \rho^{-1} + D'(0)$, and 
$$\log p_f(n) \sim \frac{2\pi}{\sqrt{3 \cdot 2^{j + 1}}} \sqrt{n}\, (\log n)^{j/2}.$$
For $j = 1$ we have $F(1, 1) = \frac{\pi^2}{6}$ and $F(1, 0) = \zeta'(2) - \frac{\gamma \pi^2}{6}$, so $\kappa = 1 - \gamma + \frac{\zeta'(2)}{\zeta(2)}$, and $D(0) = \frac{1}{2} \log 2\pi$, $D'(0) = -\zeta''(0)$. With $v_n = W(12 e^{2\kappa} n/\pi^2)$, \Cref{thm:Explicit saddle point}(ii) gives
\[
p_f(n) \sim \frac{\exp\left\{ \frac{\pi}{\sqrt{3}} (v_n - 1) \sqrt{\frac{n}{v_n}} + \frac{\log 2\pi}{4} \log \frac{12 n}{\pi^2 v_n} - \zeta''(0) \right\}}{\sqrt{2\pi \cdot \frac{\pi^2 (v_n + 1)}{6} \left( \frac{12 n}{\pi^2 v_n} \right)^{3/2}}}.
\]

\subsection{Unitary divisors, squares of divisors, and divisors of squares}
\label{sec:Divisor of squares}
The functions $2^{\omega(n)}$ (where $\omega(n)$ is the number of distinct prime factors of $n$), $d(n^2)$ and $d(n)^2$ (where $d(n)$ is the number of divisors of $n$) have the following Dirichlet series \cite{Ramanujan}, \cite[Chapter~XVII]{HardyWright}: 
\[
\sum_{n = 1}^{\infty} \frac{2^{\omega(n)}}{n^s} = \frac{\zeta(s)^2}{\zeta(2s)}, \qquad \sum_{n = 1}^{\infty} \frac{d(n^2)}{n^s} = \frac{\zeta(s)^3}{\zeta(2s)}, \qquad \sum_{n = 1}^{\infty} \frac{d(n)^2}{n^s} = \frac{\zeta(s)^4}{\zeta(2s)}.
\]
We treat $D(s) = \zeta(s)^{a}/\zeta(2s)$ with $a \in \{2, 3, 4\}$. The dominant pole $\alpha = 1$ has order $k = a$ and $A(1, a) = 6/\pi^2$, so that $C = 2/(a - 1)!$. Moreover $h_0 = 0$, $D(0) = (-\frac{1}{2})^{a - 1}$ and $D'(0) = (a - 2)(-\frac{1}{2})^{a - 1} \log 2\pi$. The other poles in $\Re(s) \geq -\varepsilon$ are the points $\varrho/2$, with $0 < \Re(\varrho/2) < \frac{1}{2}$. By \Cref{thm:Zeta quotients}, (A2) holds for every $\sigma_0 \in (\frac{1}{2}, 1)$, and, under the Riemann hypothesis, for every $\sigma_0 \in (\frac{1}{4}, 1)$. Hence \Cref{thm:Dominant saddle point}(ii) applies unconditionally, and \Cref{thm:Dominant saddle point}(i) applies under the Riemann hypothesis, with
\[
\Delta(\rho) = \left( -\tfrac{1}{2} \right)^{a - 1} \left( \log \rho^{-1} + (a - 2) \log 2\pi \right) + \sum_{\varrho} \frac{\Gamma(\frac{\varrho}{2}) \zeta(1 + \frac{\varrho}{2}) \zeta(\frac{\varrho}{2})^{a}}{2 \zeta'(\varrho)} \rho^{-\varrho/2}.
\]
By \Cref{thm:Log asymptotics},
\[
\log p_f(n) \sim \frac{2^{(3 - a)/2}}{\sqrt{(a - 1)!}} \sqrt{n}\, (\log n)^{\frac{a - 1}{2}}
\]
that is, $\sqrt{2 n \log n}$, $\frac{1}{\sqrt{2}} \sqrt{n} \log n$ and $\frac{1}{2\sqrt{3}} \sqrt{n} (\log n)^{3/2}$ for $a = 2, 3, 4$.

For $f(n) = 2^{\omega(n)}$ we have $F(1, 1) = 1$ and $F(1, 0) = \gamma - \frac{\zeta'(2)}{\zeta(2)}$, so $\kappa = 1 + \gamma - \frac{\zeta'(2)}{\zeta(2)}$. With $v_n = W(2 e^{2\kappa} n)$, the computation in the proof of \Cref{thm:Explicit saddle point}(ii), which does not use $\sigma_0 < \alpha/2$, gives $\rho_\alpha = \sqrt{v_n/(2n)}$ and $Q_2(u_\alpha) = v_n + 1$, and \Cref{thm:Dominant saddle point}(ii) gives
\begin{multline*}
p_f(n) \sim \frac{1}{\sqrt{2\pi (v_n + 1) \left( \frac{2n}{v_n} \right)^{3/2}}} \exp\Bigg\{ \sqrt{2} (v_n - 1) \sqrt{\frac{n}{v_n}} - \frac{1}{4} \log \frac{2n}{v_n} \\
+ \sum_{\varrho} \frac{\Gamma(\frac{\varrho}{2}) \zeta(1 + \frac{\varrho}{2}) \zeta(\frac{\varrho}{2})^{2}}{2 \zeta'(\varrho)} \left( \frac{2n}{v_n} \right)^{\frac{\varrho}{4}} \\
- \frac{1}{2 (v_n + 1)} \sqrt{\frac{v_n}{2n}} \bigg( \frac{1}{2} - \sum_{\varrho} \frac{\varrho\, \Gamma(\frac{\varrho}{2}) \zeta(1 + \frac{\varrho}{2}) \zeta(\frac{\varrho}{2})^{2}}{4 \zeta'(\varrho)} \left( \frac{2n}{v_n} \right)^{\frac{\varrho}{4}} \bigg)^{2} \Bigg\}.
\end{multline*}
The last line is the term $-\rho_\alpha^{3} \Delta'(\rho_\alpha)^2/(2 Q_2(u_\alpha))$ of \cref{eq:Dominant formula corrected}. Under the Riemann hypothesis it may be omitted, by \Cref{thm:Explicit saddle point}(ii).

\subsection{Partitions into \texorpdfstring{$r$}{r}-free parts}
\label{sec:r-free}
Fix $r \geq 2$, and let $f(n) = 1$ if $n$ is $r$-free, that is, not divisible by $p^r$ for any prime $p$, and $f(n) = 0$ otherwise. For $r = 2$, we have $f(n) = \mu(n)^2$ where $\mu(n)$ is the M\"obius function. Then $D(s) = \zeta(s)/\zeta(rs)$. The dominant pole $\alpha = 1$ is simple with residue $1/\zeta(r)$, so $C = \pi^2/(3\zeta(r))$, and $h_0 = 0$, $D(0) = 1$, $D'(0) = -(r - 1) \log 2\pi$. The other poles are the points $\varrho/r$, with $0 < \Re(\varrho/r) < 1/r$, and the real points $-2/r, -4/r, \ldots$ coming from the trivial zeros. By \Cref{thm:Zeta quotients}, (A2) holds for every $\sigma_0 \in (\frac{1}{r}, 1)$.

For $r \geq 3$ we can take $\sigma_0 < \frac{1}{2}$, and \Cref{thm:Explicit saddle point}(i) gives, unconditionally,
\begin{multline*}
p_f(n) \sim \frac{1}{\sqrt{2 \pi \frac{\pi^2}{3\zeta(r)} \left( \frac{6\zeta(r) n}{\pi^2} \right)^{3/2}}} \exp\Bigg\{ \pi \sqrt{\frac{2n}{3\zeta(r)}} + \frac{1}{2} \log \frac{6\zeta(r) n}{\pi^2} - (r - 1) \log 2\pi \\
+ \sum_{\varrho} \frac{\Gamma(\frac{\varrho}{r}) \zeta(\frac{\varrho}{r}) \zeta(1 + \frac{\varrho}{r})}{r \zeta'(\varrho)} \left( \frac{6\zeta(r) n}{\pi^2} \right)^{\frac{\varrho}{2r}} \Bigg\}.
\end{multline*}
The residues at the poles $-2/r, -4/r, \ldots$ may be omitted, but they tend to $0$ only like a power of $n^{-1/r}$, with large coefficients. For moderate $n$ the approximation is much better if the residues at $z = -2/r$, $z = -1$ and $z = -4/r$ are included.

For $r = 2$, \Cref{thm:Dominant saddle point}(ii) applies unconditionally. Under the Riemann hypothesis, (A2) holds for $\sigma_0 \in (\frac{1}{4}, 1)$, and \Cref{thm:Explicit saddle point}(i) gives
\[
p_f(n) \sim \frac{1}{4 \pi^{3/2} n^{1/4}} \exp\bigg\{ 2\sqrt{n} + \sum_{\varrho} \frac{\Gamma(\frac{\varrho}{2}) \zeta(\frac{\varrho}{2}) \zeta(1 + \frac{\varrho}{2})}{2 \zeta'(\varrho)} n^{\frac{\varrho}{4}} \bigg\}.
\]
Each term of the sums over $\varrho$ has modulus of order $n^{\Re(\varrho)/(2r)}$, which tends to infinity with $n$, so these terms cannot be omitted. For $r = 2$, Kot\v{e}\v{s}ovec \cite{OEIS} stated the formula that $p_f(n)$ is asymptotic to 
\[
K(n) := \frac{e^{2\sqrt{n}}}{4\pi^{3/2} n^{1/4}}, 
\]
which is the formula above but with the sum over $\varrho$ omitted. This formula is false, which can be shown even without the Riemann hypothesis. Suppose that $p_f(n)/K(n)$ is bounded between two positive constants. By applying our saddle-point method to $K(n)$, we have $\sum_{n \geq 1} K(n) e^{-n\rho} \sim e^{1/\rho}/(2\pi\rho)$ as $\rho \to 0^{+}$, and hence $\Psi(\rho) = \rho^{-1} + \log \rho^{-1} + O(1)$. Then the Mellin transform $\Gamma(s) \zeta(s + 1) \zeta(s)/\zeta(2s)$ of $\Psi$ continues analytically to the half-plane $\Re(s) > 0$, apart from the pole at $s = 1$. This is a contradiction, since it has a pole at $s = \varrho_1/2$, where $\varrho_1 = \frac{1}{2} + 14.1347\ldots i$ is the first zero of $\zeta$: there $\zeta(2s)$ vanishes, while $\Gamma(s)$, $\zeta(1 + s)$ and $\zeta(s)$ do not, as $\zeta$ has no zeros with imaginary part in $(0, 14.13)$. The same argument shows that $\log(p_f(n)/K(n)) \neq O(n^{1/8 - \delta})$ for every $\delta > 0$. In the computable range, however, the oscillating terms are negligible; for instance, $p_f(10^4)/K(10^4) = 1.047$.

\subsection{Jordan's and Dedekind's functions}
\label{sec:Jordan}
For Jordan's function $J_r(n) = n^r \prod_{p \mid n} (1 - p^{-r})$, with $J_1 = \varphi$, we have $D(s) = \zeta(s - r)/\zeta(s)$, where $r \ge 1$. The dominant pole $\alpha = r + 1$ is simple with residue $1/\zeta(r + 1)$, and
\[
C = \frac{(r + 2)!\, \zeta(r + 2)}{\zeta(r + 1)}, \qquad D(0) = -2\zeta(-r), \qquad D'(0) = -2\zeta'(-r) + 2\zeta(-r) \log 2\pi.
\]
The other poles are the zeros $\varrho$, with $0 < \Re(\varrho) < 1$, and $h_0 = 0$. By \Cref{thm:Zeta quotients}, (A2) holds for $\sigma_0 \in (1, r + 1)$. If $r \geq 2$, then $1 < \alpha/2$, and \Cref{thm:Explicit saddle point}(i) gives, unconditionally,
\begin{multline*}
p_f(n) \sim \frac{1}{\sqrt{2\pi C \left( \frac{(r + 2) n}{C} \right)^{\frac{r + 3}{r + 2}}}} \exp\Bigg\{ \frac{r + 2}{r + 1} \left( \frac{C}{r + 2} \right)^{\frac{1}{r + 2}} n^{\frac{r + 1}{r + 2}} - \frac{2\zeta(-r)}{r + 2} \log \frac{(r + 2) n}{C} \\
+ D'(0) + \sum_{\varrho} \frac{\Gamma(\varrho) \zeta(1 + \varrho) \zeta(\varrho - r)}{\zeta'(\varrho)} \left( \frac{(r + 2) n}{C} \right)^{\frac{\varrho}{r + 2}} \Bigg\}
\end{multline*}
For even $r$, $\zeta(-r) = 0$. For $r = 1$, where $\alpha/2 = 1$, \Cref{thm:Dominant saddle point}(ii) applies unconditionally, and, since (A2) holds for $\sigma_0 \in (\frac{1}{2}, 2)$ under the Riemann hypothesis, the formula above holds with $r = 1$ under the Riemann hypothesis. Then $C = 36\zeta(3)/\pi^2$, and the main term of the exponent is $\frac{3}{2} \left( \frac{12\zeta(3)}{\pi^2} \right)^{1/3} n^{2/3}$.

For Dedekind's function $\psi_r(n) = n^r \prod_{p \mid n} (1 + p^{-r})$ we have $D(s) = \zeta(s) \zeta(s - r)/\zeta(2s)$. The dominant pole $\alpha = r + 1$ is simple with residue $\zeta(r + 1)/\zeta(2r + 2)$, so $C = (r + 2)!\, \zeta(r + 1) \zeta(r + 2)/\zeta(2r + 2)$, and $D(0) = \zeta(-r)$, $D'(0) = \zeta'(-r) - \zeta(-r) \log 2\pi$. The other poles are $s = 1$, where the residue of $\Gamma(z) \zeta(z + 1) D(z) \rho^{-z}$ is $\zeta(1 - r) \rho^{-1}$, and the points $\varrho/2$. (A2) holds for $\sigma_0 \in (1, r + 1)$. If $r \geq 2$, then $1 < \alpha/2$, and \Cref{thm:Explicit saddle point}(i) gives, unconditionally,
\begin{multline*}
p_f(n) \sim \frac{1}{\sqrt{2\pi C \left( \frac{(r + 2) n}{C} \right)^{\frac{r + 3}{r + 2}}}} \exp\Bigg\{ \frac{r + 2}{r + 1} \left( \frac{C}{r + 2} \right)^{\frac{1}{r + 2}} n^{\frac{r + 1}{r + 2}} + \zeta(1 - r) \left( \frac{(r + 2) n}{C} \right)^{\frac{1}{r + 2}} \\
+ \frac{\zeta(-r)}{r + 2} \log \frac{(r + 2) n}{C} + D'(0) + \sum_{\varrho} \frac{\Gamma(\frac{\varrho}{2}) \zeta(\frac{\varrho}{2}) \zeta(\frac{\varrho}{2} - r) \zeta(1 + \frac{\varrho}{2})}{2 \zeta'(\varrho)} \left( \frac{(r + 2) n}{C} \right)^{\frac{\varrho}{2(r + 2)}} \Bigg\}
\end{multline*}
For $r = 1$, the pole $s = 1$ lies on the line $\Re(s) = \alpha/2$, so only \Cref{thm:Dominant saddle point}(ii) applies, also under the Riemann hypothesis. Here $C = 90\zeta(3)/\pi^2$, $\rho_\alpha = (C/(3n))^{1/3}$, and
\[
\Delta(\rho) = -\frac{1}{2} \rho^{-1} - \frac{1}{12} \log \rho^{-1} + D'(0) + \sum_{\varrho} \frac{\Gamma(\frac{\varrho}{2}) \zeta(\frac{\varrho}{2}) \zeta(\frac{\varrho}{2} - 1) \zeta(1 + \frac{\varrho}{2})}{2 \zeta'(\varrho)} \rho^{-\varrho/2}
\]
By the argument of \Cref{thm:Psi decomposition}, the derivative of the sum over $\varrho$ is $\ll \rho^{-\sigma - 1}$ for every $\sigma \in (\frac{1}{2}, 1)$. Hence $\rho_\alpha^{4} \Delta'(\rho_\alpha)^2 \to \frac{1}{4}$, and the last term of the exponent in \cref{eq:Dominant formula corrected} tends to the constant $-\frac{1}{8C} = -\frac{\pi^2}{720 \zeta(3)}$. Since $D'(0) = \zeta'(-1) + \frac{1}{12} \log 2\pi$, we obtain
\begin{multline*}
p_f(n) \sim \frac{1}{\sqrt{2\pi C \left( \frac{3n}{C} \right)^{4/3}}} \exp\Bigg\{ \frac{3}{2} \left( \frac{C}{3} \right)^{\frac{1}{3}} n^{\frac{2}{3}} - \frac{1}{2} \left( \frac{3n}{C} \right)^{\frac{1}{3}} - \frac{1}{36} \log \frac{3n}{C} + \zeta'(-1) + \frac{1}{12} \log 2\pi \\
- \frac{\pi^2}{720 \zeta(3)} + \sum_{\varrho} \frac{\Gamma(\frac{\varrho}{2}) \zeta(\frac{\varrho}{2}) \zeta(\frac{\varrho}{2} - 1) \zeta(1 + \frac{\varrho}{2})}{2 \zeta'(\varrho)} \left( \frac{3n}{C} \right)^{\frac{\varrho}{6}} \Bigg\}
\end{multline*}
with $C = 90\zeta(3)/\pi^2$.

\subsection{The unitary divisor sum}
\label{sec:Unitary}
Let $f(n) = \sigma^{*}(n) = \prod_{p^b \parallel n} (1 + p^b)$, the sum of the unitary divisors of $n$ (i.e. divisors $d$ of $n$ relative prime with $n/d$). Then 
$$D(s) = \frac{\zeta(s) \zeta(s - 1)}{\zeta(2s - 1)}.$$ 
The dominant pole $\alpha = 2$ is simple with residue $\zeta(2)/\zeta(3)$, so $C = \pi^2$, and $D(0) = -\frac{1}{2}$, $D'(0) = -\frac{1}{2} \log 2\pi - 6\zeta'(-1)$. At $s = 1$ the poles of $\zeta(s)$ and $\zeta(2s - 1)$ cancel. The other poles are the points $(1 + \varrho)/2$, with $\frac{1}{2} < \Re(s) < 1$, and $-\frac{1}{2}, -\frac{3}{2}, \ldots$. By \Cref{thm:Zeta quotients}, (A2) holds for $\sigma_0 \in (1, 2)$, so \Cref{thm:Dominant saddle point}(ii) applies. Under the Riemann hypothesis, (A2) holds for $\sigma_0 \in (\frac{3}{4}, 2)$, and \Cref{thm:Explicit saddle point}(i) gives
\begin{multline*}
p_f(n) \sim \frac{1}{\sqrt{2\pi \cdot \pi^2 \left( \frac{3n}{\pi^2} \right)^{4/3}}} \exp\Bigg\{ \frac{3}{2} \left( \frac{\pi^2}{3} \right)^{\frac{1}{3}} n^{\frac{2}{3}} - \frac{1}{6} \log \frac{3n}{\pi^2} - \frac{1}{2} \log 2\pi - 6\zeta'(-1) \\
+ \sum_{\varrho} \frac{\Gamma(\frac{1 + \varrho}{2}) \zeta(\frac{1 + \varrho}{2}) \zeta(\frac{\varrho - 1}{2}) \zeta(\frac{3 + \varrho}{2})}{2 \zeta'(\varrho)} \left( \frac{3n}{\pi^2} \right)^{\frac{1 + \varrho}{6}} \Bigg\}.
\end{multline*}

\subsection{Subgroups of free abelian groups and commuting tuples of permutations}
\label{sec:Commuting}
For a finitely generated group $G$, let $a_n(G)$ be the number of subgroups of index $n$ in $G$, and 
$\zeta_G(s) := \sum_{n \geq 1} a_n(G) n^{-s}$ its {\it subgroup zeta function} \cite{LubotzkySegal}. For $f(n) = a_n(\mathbb{Z}^d)$ we have 
$$
D(s) = \zeta(s) \zeta(s - 1) \cdots \zeta(s - d + 1),
$$ 
and $p_f(n) = |\mathrm{Hom}(\mathbb{Z}^{d + 1}, S_n)|/n!$ is the number of commuting $(d + 1)$-tuples of permutations of the symmetric group $S_n$, divided by $n!$ \cite{BryanFulman}. For $d = 1$ this is $p(n)$. The dominant pole $\alpha = d$ is simple with residue $\prod_{j = 2}^{d} \zeta(j)$, so $C = (d + 1)! \prod_{j = 2}^{d + 1} \zeta(j)$. The other poles are $1, \ldots, d - 1$, and $D(0) = \prod_{i = 0}^{d - 1} \zeta(-i)$, which vanishes for $d \geq 3$. (A2) holds for $\sigma_0 \in (d - 1, d)$.

For $d = 2$, we have $f(n) = a_n(\mathbb{Z}^2) = \sigma(n)$, and $p_f(n)$ counts commuting triples of permutations. Here $C = \pi^2 \zeta(3)$, the pole $s = 1$ lies on the line $\Re(s) = \alpha/2$, and (A2) holds with $\sigma_0 \in (1, \frac{4}{3})$, so \Cref{thm:Dominant saddle point}(ii) applies. We have $\rho_\alpha = (C/(3n))^{1/3}$ and
\[
\Delta(\rho) = -\frac{\pi^2}{12} \rho^{-1} + \frac{1}{24} \log \rho^{-1} + \frac{1}{24} \log 2\pi - \frac{1}{2} \zeta'(-1)
\]
so that $\rho_\alpha^{4} \Delta'(\rho_\alpha)^2 = \frac{\pi^4}{144} + O(\rho_\alpha)$, and the last term of the exponent in \cref{eq:Dominant formula corrected} is $-\frac{\pi^2}{288 \zeta(3)} + o(1)$. Hence
\[
p_f(n) \sim \frac{\exp\left\{ \frac{3}{2} \left( \frac{C}{3} \right)^{\frac{1}{3}} n^{\frac{2}{3}} - \frac{\pi^2}{12} \left( \frac{3n}{C} \right)^{\frac{1}{3}} + \frac{1}{72} \log \frac{3n}{C} + \frac{1}{24} \log 2\pi - \frac{1}{2} \zeta'(-1) - \frac{\pi^2}{288\zeta(3)} \right\}}{\sqrt{2\pi C \left( \frac{3n}{C} \right)^{4/3}}}
\]
with $C = \pi^2 \zeta(3)$. This is the formula of Bringmann, Franke and Heim \cite{BringmannFrankeHeim} for commuting triples, including the constant $-\pi^2/(288\zeta(3))$, which here comes from the last term of the exponent in \cref{eq:Dominant formula corrected}. For $d \geq 3$, the pole $s = d - 1$ satisfies $d - 1 \geq 2d/3$, so parts (i) and (ii) of \Cref{thm:Dominant saddle point} do not apply. Part (iii) applies with $\sigma_0$ close to $d - 1$ and $J = \lfloor \log_2 d \rfloor$, for instance $J = 1$ for $d = 3$. All poles $1, \ldots, d$ are simple and $D(0) = 0$, so that the residues are $r_m x^{-m}$ at $m = 1, \ldots, d - 1$ and $D'(0)$ at $0$, where
\[
r_m := \Gamma(m) \zeta(m + 1) \prod_{\substack{0 \leq i \leq d - 1 \\ i \neq m - 1}} \zeta(m - i)
\]
Thus \cref{eq:Delta bound} holds with $\sigma_1 = d - 1$ and $h_1 = 0$, and \Cref{thm:Exponent expansion} applies with $N = d$. Since $\rho_\alpha = (C/((d + 1) n))^{1/(d + 1)}$, the exponent in \cref{eq:Expanded formula} is a finite sum of powers of $n^{1/(d + 1)}$, and the terms that tend to $0$ may be omitted. For commuting quadruples, $d = 3$, we have $C = 2\pi^6 \zeta(3)/45$, $r_2 = -\zeta(2) \zeta(3)/2$, $r_1 = \zeta(2)/24$ and $D'(0) = \zeta'(-2)/24$. Computing $T_2$ and $T_3$, we obtain
\begin{multline*}
p_f(n) \sim \frac{\pi^{1/4} \zeta(3)^{1/8}}{2^{13/8}\, 3^{1/4}\, 5^{1/8}\, n^{5/8}} \exp\Bigg\{ \frac{4}{3} \left( \frac{\pi^6 \zeta(3)}{90} \right)^{\frac{1}{4}} n^{\frac{3}{4}} - \frac{\sqrt{10 \zeta(3)}}{4\pi}\, n^{\frac{1}{2}} \\
+ \frac{\pi^4 - 45 \zeta(3)}{144 \pi^2} \left( \frac{90}{\pi^6 \zeta(3)} \right)^{\frac{1}{4}} n^{\frac{1}{4}} + \frac{\zeta'(-2)}{24} + \frac{5 (\pi^4 - 60 \zeta(3))}{192 \pi^6} \Bigg\}.
\end{multline*}
For commuting $5$-tuples, $d = 4$, we have $C = 2\pi^6 \zeta(3) \zeta(5)/9$ and $D'(0) = \zeta'(-2)/2880$, and $T_2$, $T_3$ and $T_4$ give
\begin{multline*}
p_f(n) \sim \frac{\pi^{1/10} (\zeta(3) \zeta(5))^{1/10}}{2^{2/5}\, 3^{1/5}\, 5^{3/5}\, n^{3/5}} \exp\Bigg\{ \frac{5}{4} \left( \frac{2\pi^6 \zeta(3) \zeta(5)}{45} \right)^{\frac{1}{5}} n^{\frac{4}{5}} - \frac{\pi^6 \zeta(3)}{540} \left( \frac{45 n}{2\pi^6 \zeta(3) \zeta(5)} \right)^{\frac{3}{5}} \\
+ \frac{\pi^2 \zeta(3) (100 \zeta(5) - \pi^4)}{14400\, \zeta(5)} \left( \frac{45 n}{2\pi^6 \zeta(3) \zeta(5)} \right)^{\frac{2}{5}} + \frac{\pi^2 \zeta(3) (100 \zeta(5) - \pi^4)}{288000\, \zeta(5)^2} \left( \frac{45 n}{2\pi^6 \zeta(3) \zeta(5)} \right)^{\frac{1}{5}} \\
+ \frac{\zeta'(-2)}{2880} - \frac{\zeta(3) (\pi^8 - 120 \pi^4 \zeta(5) + 2400 \zeta(5)^2)}{5529600\, \pi^2 \zeta(5)^3} \Bigg\}.
\end{multline*}
For $d = 3$, the leading term of the exponent and the factor $n^{-5/8}$ agree with \cite[Theorem~3.3]{BringmannFrankeHeim}, where the constants of the terms $n^{1/2}$, $n^{1/4}$ and $n^{0}$ in the exponent are described as computable constants; the formula above gives them in closed form. Bringmann, Franke and Heim also obtained complete asymptotic expansions in powers of $n^{-1/(d + 1)}$. The relative error of the two formulas above is of order $n^{-1/(d + 1)}$. For $n = 3000$ the ratio of $p_f(n)$ to the right-hand side is $0.990$ for $d = 3$ and $0.989$ for $d = 4$, whereas \cref{eq:Newton formula} with $J = 1$ gives a ratio within $10^{-4}$ of $1$ for $d = 3$. \Cref{thm:Log asymptotics} gives $\log p_f(n) \sim (1 + \frac{1}{d}) \left( \frac{C}{1 + d} \right)^{\frac{1}{1 + d}} n^{\frac{d}{1 + d}}$.

\subsection{Subgroups of the Heisenberg group}
\label{sec:Heisenberg}
Let $\mathcal{N}$ be the discrete Heisenberg group of upper unitriangular $3 \times 3$ matrices with integer entries. By Grunewald, Segal and Smith \cite{GrunewaldSegalSmith},
\begin{equation}
\label{eq:Heisenberg zeta}
D(s) = \zeta_{\mathcal{N}}(s) = \frac{\zeta(s) \zeta(s - 1) \zeta(2s - 2) \zeta(2s - 3)}{\zeta(3s - 3)}
\end{equation}
and $p_f(n)$, with $f(n) = a_n(\mathcal{N})$, counts multisets of subgroups of finite index in $\mathcal{N}$ whose indices add up to $n$. The factors $\zeta(s - 1)$ and $\zeta(2s - 3)$ both have a simple pole at $s = 2$, so $\alpha = 2$ is a pole of order $k = 2$, with $A(2, 2) = \pi^4/(72\zeta(3))$ and $C = \pi^4/12$. The other poles in $\Re(s) \geq -\varepsilon$ are the real points $\frac{3}{2}$, $1$ and $\frac{1}{3}$, the last coming from the trivial zero $\zeta(-2) = 0$ of $\zeta(3s - 3)$, and the complex points $1 + \varrho/3$, with $1 < \Re(s) < \frac{4}{3}$. Moreover $h_0 = 0$ and $D(0) = 0$. By \Cref{thm:Zeta quotients}, (A2) holds for $\sigma_0 \in (\frac{3}{2}, 2)$. Since the pole $\frac{3}{2}$ exceeds $2\alpha/3 = \frac{4}{3}$, parts (i) and (ii) of \Cref{thm:Dominant saddle point} do not apply, but part (iii) applies with $\sigma_0$ close to $\frac{3}{2}$ and $J = 2$. Near $s = 2$, $D(s) = A(2, 2)(s - 2)^{-2} + A(2, 1)(s - 2)^{-1} + O(1)$ with
\[
A(2, 1) = \frac{3\zeta(2)^2}{2\zeta(3)} \left( \gamma + \frac{\zeta'(2)}{\zeta(2)} - \frac{\zeta'(3)}{\zeta(3)} \right), \qquad \kappa = \frac{3}{2} + 2\gamma + 3\frac{\zeta'(2)}{\zeta(2)} - 2\frac{\zeta'(3)}{\zeta(3)}
\]
The computation in the proof of \Cref{thm:Explicit saddle point}(ii), which does not use $\sigma_0 < \alpha/2$, gives $\rho_\alpha = 1/y_n$ with
\[
y_n := \left( \frac{108 n}{\pi^4 v_n} \right)^{\frac{1}{3}}, \qquad v_n := W\left( \frac{108 e^{3\kappa}}{\pi^4}\, n \right)
\]
and $Q_2(u_\alpha) = \pi^4 (v_n + 1)/36$, $Q_3(u_\alpha) = \pi^4 (4 v_n + 7)/36$. The poles $p \neq 2$ with $\Re(p) > 0$ are simple, and $D(0) = 0$. Hence the residue at $0$ is $\zeta'(-2)/12$, and the residue at each of the poles $p = \frac{3}{2}$, $1$, $\frac{1}{3}$ and $1 + \varrho/3$ is $r_p x^{-p}$, where
\[
r_{3/2} = -\frac{\sqrt{\pi}\, \zeta(\frac{1}{2}) \zeta(\frac{5}{2})}{8}, \qquad r_1 = \frac{\pi^2}{144}, \qquad r_{1/3} = \frac{\Gamma(\frac{1}{3}) \zeta(\frac{1}{3}) \zeta(\frac{4}{3}) \zeta(-\frac{2}{3}) \zeta(-\frac{4}{3}) \zeta(-\frac{7}{3})}{3 \zeta'(-2)}
\]
\[
r_{1 + \varrho/3} = \frac{\Gamma(1 + \frac{\varrho}{3}) \zeta(2 + \frac{\varrho}{3}) \zeta(1 + \frac{\varrho}{3}) \zeta(\frac{\varrho}{3}) \zeta(\frac{2\varrho}{3}) \zeta(\frac{2\varrho}{3} - 1)}{3 \zeta'(\varrho)}
\]
By \Cref{thm:Zeta quotients}(i), \cref{eq:Exponential bound} holds on every line $\Re(s) = \sigma'$ with $\frac{4}{3} < \sigma' < \frac{3}{2}$, so \cref{eq:Delta bound} holds with $\sigma_1 = \frac{3}{2}$ and $h_1 = 0$, and \Cref{thm:Exponent expansion} applies with $N = 3$, the condition on $N$ being $(\log \rho^{-1})^{-3} \to 0$. For $j \geq 0$ put
\[
b_j(n) := \sum_{p} p (p + 1) \cdots (p + j - 1)\, r_p\, y_n^{p}
\]
where $p$ runs over the poles above, so that $\Delta(\rho_\alpha) = b_0(n) + \zeta'(-2)/12$ and $b_j(n) = b_j$ for $j \geq 1$. Then \cref{eq:Expanded formula} becomes
\begin{multline*}
p_f(n) \sim \frac{1}{y_n^2 \sqrt{2\pi \left( \frac{\pi^4 (v_n + 1)}{36} + \frac{1}{y_n^2} \left( b_2(n) - \frac{4 v_n + 7}{v_n + 1}\, b_1(n) \right) \right)}} \\
\times \exp\Bigg\{ \frac{1}{2} \left( \frac{\pi^4}{4} \right)^{\frac{1}{3}} \left( v_n - \frac{1}{2} \right) \left( \frac{n}{v_n} \right)^{\frac{2}{3}} + b_0(n) + \frac{\zeta'(-2)}{12} \\
- \frac{18\, b_1(n)^2}{\pi^4 (v_n + 1)\, y_n^2} + \frac{216\, b_1(n)^2 \left( 3 (v_n + 1) b_2(n) - (4 v_n + 7) b_1(n) \right)}{\pi^8 (v_n + 1)^3\, y_n^4} \Bigg\}.
\end{multline*}
The last two terms of the exponent are $T_2$ and $T_3$. Both tend to infinity, like $y_n/v_n$ and $y_n^{1/2}/v_n^2$, while $T_4$ is of order $v_n^{-3}$. The terms with $b_1(n)$ and $b_2(n)$ in the denominator may be omitted, but they improve the approximation considerably: for $n = 3000$ the ratio of $p_f(n)$ to the right-hand side is $0.9992$, and $1.0074$ without these terms. \Cref{thm:Log asymptotics} gives
\[
\log p_f(n) \sim \frac{3}{2} \left( \frac{\pi^4}{108} \right)^{\frac{1}{3}} n^{\frac{2}{3}} (\log n)^{\frac{1}{3}}
\]
Here a double pole at $\alpha$, further real poles and infinitely many complex poles occur together, so this example needs \Cref{thm:Main theorem} in its full generality.

\subsection{Subrings of \texorpdfstring{$\mathbb{Z}^3$}{Z3}}
\label{sec:Subrings}
A subring of $\mathbb{Z}^m$ is a sublattice of finite index that is closed under componentwise multiplication and contains $(1, \ldots, 1)$. For $m = 2$ there is exactly one subring of each index $n$, namely $\{(a, b) : a \equiv b \bmod n\}$. Let $f(n)$ be the number of subrings of $\mathbb{Z}^3$ of index $n$, so that $p_f(n)$ is the number of multisets of subrings of $\mathbb{Z}^3$ whose indices add up to $n$. By Datskovsky and Wright \cite{DatskovskyWright} (see also \cite{Liu, AtanasovKaplanKrakoffMenzel}),
\[
D(s) = \frac{\zeta(3s - 1) \zeta(s)^3}{\zeta(2s)^2}
\]
The dominant pole $\alpha = 1$ has order $k = 3$ and $A(1, 3) = 1/\zeta(2)$, so that $C = 1$. The other poles in $\Re(s) \geq -\varepsilon$ are $s = \frac{2}{3}$, coming from $\zeta(3s - 1)$, and the points $\varrho/2$, with $0 < \Re(\varrho/2) < \frac{1}{2}$. Moreover $h_0 = 0$ and $D(0) = \zeta(-1) \zeta(0) = \frac{1}{24}$. By \Cref{thm:Zeta quotients}, (A2) holds for every $\sigma_0 \in (\frac{2}{3}, 1)$, and (A3) holds because $f(n) \geq 1$ for all $n$, the subring $\{x : x_1 \equiv x_2 \bmod n\}$ having index $n$. The pole $\frac{2}{3} = 2\alpha/3$ lies exactly on the boundary of \Cref{thm:Dominant saddle point}(ii), which therefore does not apply, but part (iii) applies with $\sigma_0$ close to $\frac{2}{3}$ and $J = 1$. To make the formula explicit we use \Cref{thm:Exponent expansion}. By \Cref{thm:Zeta quotients}(i), \cref{eq:Exponential bound} holds on every line $\Re(s) = \sigma'$ with $\frac{1}{2} < \sigma' < \frac{2}{3}$, so \cref{eq:Delta bound} holds with $\sigma_1 = \frac{2}{3}$ and $h_1 = 0$, and $N = 2$ is admissible, the condition on $N$ being $(\log \rho^{-1})^{-4} \to 0$. At $s = 1$ we find $F(1, 2) = 1$, $F(1, 1) = 2\gamma$ and
\[
F(1, 0) = \frac{\pi^2}{12} + \frac{\gamma^2}{2} - 3\gamma_1 + \frac{\zeta''(2)}{\zeta(2)} - \left( \frac{\zeta'(2)}{\zeta(2)} \right)^2 = 2.0919\ldots
\]
where $\gamma_1$ is the Stieltjes constant in $\zeta(s) = \frac{1}{s - 1} + \gamma - \gamma_1 (s - 1) + O(|s - 1|^2)$. Hence $Q_0(u) = \frac{1}{2} u^2 + 2\gamma u + F(1, 0)$, and $\rho_\alpha = e^{-u_n}$, where $u_n$ is the largest solution of
\[
e^{2u} \left( u^2 + 2 (2\gamma + 1) u + 4\gamma + 2 F(1, 0) \right) = 2n
\]
As noted after \Cref{thm:Explicit saddle point}, $2u_n$ is the generalized Lambert function $w_n$ of part (iii) there, and by \Cref{thm:No closed form} we do not expect a closed form for it. The residues at $\frac{2}{3}$ and at $0$ are
\[
r_{2/3}\, x^{-2/3}, \quad r_{2/3} := \frac{\Gamma(\frac{2}{3}) \zeta(\frac{5}{3}) \zeta(\frac{2}{3})^3}{3 \zeta(\frac{4}{3})^2} = -1.0838\ldots, \qquad \frac{1}{24} \log x^{-1} - \frac{3}{2} \zeta'(-1) - \frac{1}{24} \log 2\pi
\]
With $\widetilde{D}(s) := \zeta(2s)^2 D(s) = \zeta(3s - 1) \zeta(s)^3$, the residue at the double pole $\varrho/2$, for a simple zero $\varrho$, is $x^{-\varrho/2} (A_\varrho + B_\varrho \log x^{-1})$, where
\[
B_\varrho := \frac{\Gamma(\frac{\varrho}{2}) \zeta(1 + \frac{\varrho}{2}) \widetilde{D}(\frac{\varrho}{2})}{4 \zeta'(\varrho)^2}, \qquad A_\varrho := \frac{1}{4 \zeta'(\varrho)^2} \frac{d}{dz} \Big( \Gamma(z) \zeta(z + 1) \widetilde{D}(z) \Big) \Big|_{z = \varrho/2} - \frac{2 \zeta''(\varrho)}{\zeta'(\varrho)} B_\varrho
\]
Put
\[
b_1(n) := \frac{2 r_{2/3}}{3} e^{2u_n/3} + \frac{1}{24} + \sum_{\varrho} e^{\varrho u_n/2} \left( \frac{\varrho}{2} (A_\varrho + B_\varrho u_n) + B_\varrho \right)
\]
so that $b_1(n) = b_1$ in the notation of \Cref{thm:Exponent expansion}, and note that $Q_2(u) = u^2 + (4\gamma + 3) u + 2F(1, 0) + 6\gamma + 1$. Then \cref{eq:Expanded formula}, with the term in the denominator omitted, gives
\begin{multline*}
p_f(n) \sim \frac{1}{\sqrt{2\pi e^{3u_n} Q_2(u_n)}} \exp\Bigg\{ e^{u_n} \left( u_n^2 + (4\gamma + 1) u_n + 2\gamma + 2F(1, 0) \right) + r_{2/3}\, e^{2u_n/3} \\
+ \frac{u_n}{24} - \frac{3}{2} \zeta'(-1) - \frac{\log 2\pi}{24} + \sum_{\varrho} e^{\varrho u_n/2} (A_\varrho + B_\varrho u_n) - \frac{b_1(n)^2}{2 e^{u_n} Q_2(u_n)} \Bigg\}
\end{multline*}
The last term of the exponent is $T_2$, and it tends to infinity like $e^{u_n/3}/u_n^2$. For $n = 3000$ the ratio of $p_f(n)$ to the right-hand side is $0.990$, and $0.998$ with the full denominator of \cref{eq:Expanded formula}. \Cref{thm:Log asymptotics} gives
\[
\log p_f(n) \sim \frac{1}{\sqrt{2}} \sqrt{n} \log n
\]
the same leading term as for $d(n^2)$ in \Cref{sec:Divisor of squares}.

\subsection{Representations of products of compact groups}
\label{sec:Representations}
For a compact Lie group $G$, let $r_G(n)$ be the number of irreducible representations of $G$ of dimension $n$, up to isomorphism. By complete reducibility, $p_f(n)$ with $f = r_G$ is the number of $n$-dimensional representations of $G$ up to isomorphism, and $D(s) = \sum_{n \geq 1} r_G(n) n^{-s}$ is the {\it Witten zeta function} of $G$ \cite{Witten}. For $G = SU(2)^r$ we have $r_G = d_r$, as in \Cref{sec:Divisor functions}. For $G = SU(3)$, the irreducible representations $W_{j, k}$, $j, k \geq 1$, have dimension $jk(j + k)/2$, so that $D(s) = 2^{s} \omega(s)$ with
\[
\omega(s) := \sum_{j, k \geq 1} \frac{1}{(jk(j + k))^{s}}.
\]
Romik \cite[Theorem~1.2]{Romik} proved that $\omega$ continues meromorphically to $\mathbb{C}$, with simple poles at $s = \frac{2}{3}$, with residue $\Gamma(\frac{1}{3})^3/(2\pi\sqrt{3})$, and at $s = \frac{1}{2} - m$ for $m \geq 0$, and no other singularities, and that for every bounded interval $I$ there is $M$ with $\omega(\sigma + it) \ll |t|^{M}$ for $\sigma \in I$ and $|t| \geq 1$. He also determined the asymptotics of the number of $n$-dimensional representations of $SU(3)$ \cite{Romik}, and Bringmann and Franke \cite{BringmannFranke2023} gave an asymptotic expansion; see also \cite{BridgesBringmannFranke2024}. In our notation $\alpha = \frac{2}{3}$ and $k = 1$, the residue of $D(s)$ at $\frac{2}{3}$ is $R := 2^{2/3}\Gamma(\frac{1}{3})^3/(2\pi\sqrt{3})$, and $C = \frac{10}{9}\Gamma(\frac{2}{3})\zeta(\frac{5}{3}) R$. The pole at $\frac{1}{2}$ exceeds $2\alpha/3 = \frac{4}{9}$, so parts (i) and (ii) of \Cref{thm:Dominant saddle point} do not apply, but part (iii) applies with $\sigma_0$ close to $\frac{1}{2}$ and $J = 2$. Since $\omega(0) = \frac{1}{3}$ and $\omega'(0) = \log 2\pi$ \cite[Theorems~1.3 and~1.5]{Romik}, the residues at $\frac{1}{2}$ and at $0$ are $\sqrt{2\pi}\, \zeta(\frac{1}{2}) \zeta(\frac{3}{2})\, x^{-1/2}$ and $\frac{1}{3} \log x^{-1} + \frac{1}{3} \log 2 + \log 2\pi$. By Romik's bound, \cref{eq:Exponential bound} holds on every line $\Re(s) = \sigma'$ with $0 < \sigma' < \frac{1}{2}$, so \cref{eq:Delta bound} holds with $\sigma_1 = \frac{1}{2}$ and $h_1 = 0$, and \Cref{thm:Exponent expansion} applies with $N = 4$. Here $\rho_\alpha = (3C/(5n))^{3/5}$, and expanding $T_2$, $T_3$ and $T_4$ in powers of $n^{1/10}$, we obtain
\[
p_f(n) \sim \frac{K}{n^{3/5}} \exp\left( A_1 n^{2/5} - A_2 n^{3/10} - A_3 n^{1/5} - A_4 n^{1/10} \right)
\]
where, with $X = (\frac{1}{9}\Gamma(\frac{1}{3})^2 \zeta(\frac{5}{3}))^{3/10}$ and $Y = -\sqrt{\pi}\, \zeta(\frac{1}{2}) \zeta(\frac{3}{2})$,
\[
A_1 = 5X^2, \quad A_2 = \frac{Y}{X}, \quad A_3 = \frac{3 Y^2}{80 X^4}, \quad A_4 = \frac{11 Y^3}{3200 X^7}, \quad K = \frac{2\sqrt{3\pi}}{\sqrt{5}}\, X^{1/3} \exp\left( -\frac{Y^4}{2560 X^{10}} \right)
\]
This is Romik's theorem \cite[Theorem~1.1]{Romik}. The terms $A_3 n^{1/5}$ and $A_4 n^{1/10}$ come from $T_2$ and $T_3$, and the exponential factor in $K$ comes from $T_4$. The error terms are powers of $n^{-1/10}$, and the convergence is slow: for $n = 3000$ the ratio of $p_f(n)$ to the right-hand side is $0.86$. For comparison, \cref{eq:Newton formula} with $J = 3$ gives the ratio $1.007$, and \cref{eq:Critical form}, with the critical point $\hat{\rho}$ of $\Phi_\Delta$ from the proof of \Cref{thm:Exponent expansion}, gives $0.996$. \Cref{thm:Log asymptotics} gives Romik's leading term $\log p_f(n) \sim 5 X^2 n^{2/5}$.

Now let $G = SU(3)^r$ with $r \geq 2$. Then $r_G$ is the $r$-fold Dirichlet convolution of $r_{SU(3)}$, and $D(s) = 2^{rs} \omega(s)^r$. The dominant pole $\alpha = \frac{2}{3}$ has order $k = r$, with $A(\frac{2}{3}, r) = R^{r}$, so that
\[
C = \frac{10\, \Gamma(\frac{2}{3}) \zeta(\frac{5}{3})}{9\, (r - 1)!} R^{r}
\]
The only other pole in $\Re(s) \geq -\varepsilon$ is $s = \frac{1}{2}$, now of order $r$, and $h_0 = 0$, $D(0) = 3^{-r}$ by \cite[Theorem~1.3]{Romik}. Romik's theorem gives (A1), and (A2) with every $\sigma_0 \in (\frac{1}{2}, \frac{2}{3})$, where the bounds on horizontal segments hold for every sequence $T_j \to \infty$ since all poles are real. (A3) holds because $r_G(1) = 1$. Hence \Cref{thm:Dominant saddle point}(iii) applies with $J = 2$, as for $SU(3)$. \Cref{thm:Log asymptotics} gives
\[
\log p_f(n) \sim \frac{5}{2} \left( \frac{3^{r} C}{5^{r}} \right)^{\frac{3}{5}} n^{\frac{2}{5}} (\log n)^{\frac{3(r - 1)}{5}}
\]
Since the dominant pole has order $r \geq 2$, these Dirichlet series are not covered by \cite{GranovskyStark2012, GranovskyStark2015, Bridges2023}. The convergence in \Cref{thm:Log asymptotics} is slow here: the relative error is $O(\log \log n/\log n)$, and the residue at the pole $\frac{1}{2}$, which is close to $\alpha$, is still large for computable $n$. By Romik's bound, \cref{eq:Delta bound} holds with $\sigma_1 = \frac{1}{2}$ and $h_1 = r - 1$, and \Cref{thm:Exponent expansion} applies with $N = 4$, the condition on $N$ being $\rho^{1/6} (\log \rho^{-1})^{r - 1} \to 0$.

For $r = 2$ we write the result out. Let
\[
2^{s} \omega(s) = \frac{R}{s - \frac{2}{3}} + R_0 + O(|s - \tfrac{2}{3}|), \qquad 2^{s} \omega(s) = \frac{R'}{s - \frac{1}{2}} + R_0' + O(|s - \tfrac{1}{2}|)
\]
where $R' = \sqrt{2}\,\zeta(\frac{1}{2})$ by \cite[Theorem~1.2]{Romik}. We do not know closed forms for $R_0$ and $R_0'$; from a Mellin--Barnes integral for $\omega(s)$ we computed $R_0 = -11.9040\ldots$ and $R_0' = -16.3581\ldots$. Then $C = 25.1269\ldots$, and $\rho_\alpha$ is given by the Lambert function as in the proof of \Cref{thm:Explicit saddle point}(ii), with
\[
\kappa = \frac{3}{2} + \frac{\Gamma'(\frac{2}{3})}{\Gamma(\frac{2}{3})} + \frac{\zeta'(\frac{5}{3})}{\zeta(\frac{5}{3})} + \frac{2R_0}{R} = -9.3363\ldots
\]
Namely, $\rho_\alpha = 1/y_n$, where
\[
y_n := \left( \frac{25 n}{9 C v_n} \right)^{\frac{3}{5}}, \qquad v_n := W\left( \frac{25 e^{5\kappa/3}}{9 C}\, n \right)
\]
and the numbers $q_j$ of \Cref{thm:Exponent expansion} are $q_2 = \frac{3}{5} C (v_n + 1)$, $q_3 = \frac{1}{5} C (8 v_n + 13)$ and $q_4 = \frac{1}{15} C (88 v_n + 183)$. The residues of $\Gamma(z) \zeta(z + 1) D(z) x^{-z}$ at the double pole $\frac{1}{2}$ and at $0$ are $x^{-1/2} (S_0 + S_1 \log x^{-1})$ and $\frac{1}{9} \log x^{-1} + \frac{2}{9} \log 2 + \frac{2}{3} \log 2\pi$, where
\[
S_1 = 2\sqrt{\pi}\, \zeta(\tfrac{1}{2})^2 \zeta(\tfrac{3}{2}), \qquad S_0 = 2\sqrt{\pi}\, \zeta(\tfrac{1}{2})^2 \Big( \zeta'(\tfrac{3}{2}) - (\gamma + 2 \log 2)\, \zeta(\tfrac{3}{2}) \Big) + 2\sqrt{2\pi}\, \zeta(\tfrac{1}{2}) \zeta(\tfrac{3}{2}) R_0'
\]
With $S(n) := S_0 + S_1 \log y_n$, the numbers $b_j$ of \Cref{thm:Exponent expansion} are
\[
b_1 = \sqrt{y_n} \left( \tfrac{1}{2} S(n) + S_1 \right) + \tfrac{1}{9}, \qquad b_2 = \sqrt{y_n} \left( \tfrac{3}{4} S(n) + 2 S_1 \right) + \tfrac{1}{9}, \qquad b_3 = \sqrt{y_n} \left( \tfrac{15}{8} S(n) + \tfrac{23}{4} S_1 \right) + \tfrac{2}{9}
\]
and \cref{eq:Expanded formula} becomes
\begin{multline*}
p_f(n) \sim \frac{1}{y_n^{4/3} \sqrt{2\pi \left( q_2 + y_n^{-2/3} (b_2 - b_1 q_3/q_2) \right)}} \\
\times \exp\Bigg\{ \frac{3}{2} \left( \frac{3}{5} \right)^{\frac{1}{5}} C^{\frac{3}{5}} \left( v_n - \frac{3}{2} \right) \left( \frac{n}{v_n} \right)^{\frac{2}{5}} + \sqrt{y_n}\, S(n) \\
+ \frac{1}{9} \log y_n + \frac{2}{9} \log 2 + \frac{2}{3} \log 2\pi + T_2 + T_3 + T_4 \Bigg\},
\end{multline*}
where $T_2$, $T_3$ and $T_4$ are given by \Cref{thm:Exponent expansion} with $\rho_\alpha^{\alpha} = y_n^{-2/3}$. Unlike in the previous examples, $T_4$ does not tend to $0$; it grows like $\log n$. Since $\kappa$ is large and negative, $v_n$ is small unless $n$ is very large, and this formula becomes accurate only for very large $n$: the ratio of its right-hand side to the right-hand side of \cref{eq:Critical form} is $e^{-1.09}$ for $n = 10^{20}$ and $0.997$ for $n = 10^{40}$. For computable $n$, $\rho_\alpha$ is far from $\rho_n$. For $n = 1000$, $\rho_\alpha = 8.8 \cdot 10^{-5}$ while $\hat{\rho} = 0.056$, and the Newton steps of \Cref{thm:Dominant saddle point}(iii) approach $\hat{\rho}$ only slowly. In this range one should use \cref{eq:Critical form} itself: the ratio of $p_f(n)$ to its right-hand side is $0.9965$ for $n = 1000$ and $0.9980$ for $n = 3000$. For $r \geq 3$ the dominant pole has order $r \geq 3$, and $\rho_\alpha$ is the generalized Lambert function of \Cref{thm:Explicit saddle point}(iii).

\subsection{Ordered factorizations}
\label{sec:Ordered factorizations}
Let $f(n) = H(n)$ be the number of ordered factorizations of $n$, that is, of representations $n = m_1 m_2 \cdots m_{\ell}$ with $\ell \geq 0$ factors $m_i \geq 2$, where the order of the factors is taken into account, so that $H(1) = 1$. Then $p_f(n)$ is the number of partitions of $n$ in which every part $m$ carries one of the $H(m)$ ordered factorizations of $m$. Since 
$$
\sum_{n \geq 1} \frac{H(n)}{n^{s}} = \sum_{\ell \geq 0} (\zeta(s) - 1)^{\ell},
\text{ we have }
D(s) = \frac{1}{2 - \zeta(s)},
$$
and so $D(s)$ is meromorphic in $\mathbb{C}$, and its poles are the solutions of $\zeta(s) = 2$. The dominant pole is the real solution $\alpha = 1.7286\ldots$ of $\zeta(\alpha) = 2$ \cite{Kalmar}. It is simple with residue $-1/\zeta'(\alpha)$, so $k = 1$ and
$
C = -\frac{(1 + \alpha)\, \Gamma(\alpha + 1)\, \zeta(\alpha + 1)}{\zeta'(\alpha)} = 3.0034\ldots.
$
Condition (A3) holds because $H(n) \geq 1$ for all $n$.

However, \cref{eq:Complex poles to the left} fails. If $\Re(s) > 1$ and $\Im(s) \neq 0$, then $|\zeta(s)| < \zeta(\Re(s))$, so every complex solution of $\zeta(s) = 2$ satisfies $\Re(s) < \alpha$. But the real parts of these solutions accumulate at $\alpha$. To see this, let
\[
\zeta(s; \theta) := \frac{1 - 2^{-s}}{1 - e^{i\theta} 2^{-s}}\, \zeta(s), \qquad \theta \in \mathbb{R}
\]
which is the Euler product of $\zeta$ with the factor at $p = 2$ rotated by $e^{i\theta}$. Since $\zeta'(\alpha) \neq 0$, for every small $\theta \neq 0$ the equation $\zeta(s; \theta) = 2$ has a simple solution $s_\theta$, with $s_\theta \to \alpha$ as $\theta \to 0$. Moreover $\Re(s_\theta) < \alpha$, because $|\zeta(s; \theta)| < \zeta(\Re(s))$ for $\Re(s) > 1$ and $\theta \notin 2\pi\mathbb{Z}$. The numbers $\log p$ are linearly independent over $\mathbb{Q}$, so by Kronecker's theorem \cite[Chapter~XXIII]{HardyWright} there are arbitrarily large $t$ such that $2^{-it}$ is close to $e^{i\theta}$ and $p^{-it}$ is close to $1$ for all primes $3 \leq p \leq P$. For such $t$ and large $P$, $\zeta(s + it)$ is uniformly close to $\zeta(s; \theta)$ on a small disc around $s_\theta$, and Rouché's theorem gives a solution of $\zeta(s) = 2$ in the disc translated by $it$. Hence every strip $\alpha - \delta < \Re(s) < \alpha$ contains infinitely many poles of $D(s)$, and $\sup_{\beta \in \mathcal{B}} \Re(\beta) = \alpha$. Since \cref{eq:Complex poles to the left} is used to show that the pole $\alpha$ dominates $\Psi$, none of \Cref{thm:Main theorem}, \Cref{thm:Dominant saddle point}, \Cref{thm:Explicit saddle point} and \Cref{thm:Log asymptotics} can be applied.

The logarithmic asymptotics can nevertheless be obtained directly. By Kalmár's theorem \cite{Kalmar}, $\sum_{n \leq x} H(n) \sim -x^{\alpha}/(\alpha \zeta'(\alpha))$. By partial summation, $\sum_{n \geq 1} H(n) e^{-n\rho} \sim -\Gamma(\alpha) \rho^{-\alpha}/\zeta'(\alpha)$, and since $\Psi(\rho) = \sum_{m \geq 1} m^{-1} \sum_{n \geq 1} H(n) e^{-mn\rho}$, dominated convergence gives
\[
\Psi(\rho) \sim -\frac{\Gamma(\alpha) \zeta(\alpha + 1)}{\zeta'(\alpha)} \rho^{-\alpha} = \frac{C}{\alpha (1 + \alpha)} \rho^{-\alpha}.
\]
Since $f(1) = 1$, the generating function is $(1 - q)^{-1}$ times a power series with non-negative coefficients, so $p_f(n)$ is non-decreasing, and $\log \sum_{m \leq n} p_f(m) = \log p_f(n) + O(\log n)$. Kohlbecker's Tauberian theorem \cite{Kohlbecker} therefore gives
\[
\log p_f(n) \sim \left( 1 + \frac{1}{\alpha} \right) \left( \frac{C}{1 + \alpha} \right)^{\frac{1}{1 + \alpha}} n^{\frac{\alpha}{1 + \alpha}} = 1.6349\ldots\, n^{0.6335\ldots},
\]
which is the conclusion of \Cref{thm:Log asymptotics} with $k = 1$. An asymptotic formula for $p_f(n)$ itself would have to take into account infinitely many poles close to the line $\Re(s) = \alpha$, and it is unknown. 

\section{Some further remarks}
\label{sec:Conclusion}



\subsection{}
A pole of order $k \geq 2$, as for $d_r(n)$, $(\log n)^j$, $2^{\omega(n)}$, $d(n^2)$, $d(n)^2$, the Heisenberg group, the representations of $SU(3)^r$ and the subrings of $\mathbb{Z}^3$, adds the factor $(\log n)^{(k - 1)/(1 + \alpha)}$ to the growth and makes the saddle point harder to describe. For $k = 2$ it is still given by the Lambert function, as in the formulas for $d(n)$, which counts the representations of $\mathrm{Spin}(4)$, for $\log n$ and for $2^{\omega(n)}$. For $k \geq 3$, $\rho_\alpha$ is the inverse of the elementary function $\mathcal{H}_k$ of \Cref{sec:Motivation}, and if \Cref{thm:No closed form} holds, it can only be computed numerically, for instance by Newton's method with a number of steps that grows with $n$.

\subsection{}
A pole on the line $\Re(s) = \alpha/2$ does not affect the growth of $p_f(n)$, but when $\rho_n$ is replaced by $\rho_\alpha$ it leaves a constant factor, coming from the last term of the exponent in \cref{eq:Dominant formula corrected}. For the weight $\sigma(n)$, where $p_f(n)$ is the number of commuting triples of permutations of $S_n$ divided by $n!$, this factor is $\exp(-\pi^2/(288\zeta(3)))$, in accordance with \cite{BringmannFrankeHeim}; for $\psi_1(n)$ it is $\exp(-\pi^2/(720\zeta(3)))$.

\subsection{}
Our method depends essentially on the assumption \cref{eq:Complex poles to the left}. Complex poles on the line $\Re(s) = \alpha$ contribute oscillating terms of the same size as the main term, as in the example $f(n) = 1 + \cos(t_0 \log n)$ of \Cref{sec:Contributions} and for the weighted $m$-ary partitions, where $D(s) = (1 - m^{j - s})^{-1}$; for $m$-ary partitions these oscillations are known to produce periodic factors in the asymptotic formula \cite{Mahler, deBruijn}. If complex poles accumulate at the line $\Re(s) = \alpha$, as for the ordered factorizations $H(n)$ \cite{Kalmar, Hille, ChorLemkeMador} discussed in \Cref{sec:Ordered factorizations}, they cannot be separated from $\alpha$ by shifting the line of integration. In that example only the logarithmic asymptotics is available, and an asymptotic formula for the partitions generated by $\prod_{n \geq 1} (1 - q^n)^{-H(n)}$ is unknown. 

\subsection{}
For other natural weights the Dirichlet series does not continue to $\Re(s) \geq -\varepsilon$. For $d(n)^j$ and $d(n^j)$ with $j \geq 3$, Estermann's theorem \cite{Estermann} gives the natural boundary $\Re(s) = 0$, and the subgroup zeta functions of nilpotent groups are in general known to continue only slightly to the left of their abscissa of convergence \cite{duSautoyGrunewald}; the Heisenberg group is an exception, since its zeta function is a quotient of Riemann zeta functions \cite{GrunewaldSegalSmith}. The same applies to the number of isomorphism classes of $G$-sets with $n$ elements for a non-abelian nilpotent group $G$. Here the weights are the numbers $c_n(G)$ of conjugacy classes of subgroups of index $n$, with Dirichlet series $\zeta_{\mathbb{Z} \times G}(s + 1)/\zeta(s + 1)$. By Burnside's lemma this number equals $|\mathrm{Hom}(\mathbb{Z} \times G, S_n)|/n!$, and for finitely generated torsion-free nilpotent $G$ the asymptotics of its logarithm follows from the work of M\"uller and Schlage-Puchta \cite{MullerSchlagePuchta}, who determined the asymptotics of $\log(|\mathrm{Hom}(\Gamma, S_n)|/n!)$ for finitely generated torsion-free nilpotent groups $\Gamma$. The estimates of \Cref{sec:Contributions} up to \Cref{thm:Saddle point conclusion} use $D(s)$ only in a half-plane $\Re(s) \geq \sigma_0$ with $\sigma_0 < \alpha$. They remain valid, together with \Cref{thm:Main theorem} and \Cref{thm:Log asymptotics}, whenever $D(s)$ continues meromorphically to such a half-plane, has no pole there other than $\alpha$, and satisfies \cref{eq:Exponential bound} on the line $\Re(s) = \sigma_0$ and \cref{eq:Horizontal bound} for $\sigma_0 \leq \sigma \leq c$. 
We do not know 
whether the Dirichlet series above have these properties, 
and the residues at poles further to the left are not available. 

\subsection{}
Several counting problems for algebraic structures lead to Dirichlet series with poles of higher order and are natural candidates for further applications. For subrings of $\mathbb{Z}^m$ with $m = 4, 5$, the pole at $s = 1$ has order $\binom{m}{2}$ \cite{KaplanMarcinekTaklooBighash, AtanasovKaplanKrakoffMenzel}, but the local factors of the subring zeta function have non-trivial numerators, and a meromorphic continuation to $\Re(s) \geq -\varepsilon$ does not seem to be known; for $m \geq 6$ only bounds for the number of subrings are known \cite{AtanasovKaplanKrakoffMenzel}. Subgroup zeta functions of nilpotent groups often have a pole of higher order at their abscissa of convergence \cite{duSautoyGrunewald}, and many of them are computed explicitly in \cite{duSautoyWoodward}. Whenever such a zeta function is a finite product of shifted Riemann zeta functions and their reciprocals, as for the Heisenberg group, \Cref{thm:Zeta quotients} verifies (A2), and \Cref{thm:Main theorem} and \Cref{thm:Dominant saddle point}(iii) apply. For representations of products $G^r$ of other compact groups, 
the meromorphic continuation and the growth of the Witten zeta function of $G$ on vertical lines would have to be established first.

\subsection{}
The dominant singularity may also be a branch point rather than a pole. This is the case for partitions into prime powers \cite{Gafni2021}, where the Dirichlet series behaves like $\log \zeta(js)$, and for Dirichlet series behaving like $(s - \alpha)^{-w}$ with non-integer $w$, which by Delange's theorem \cite{Delange1954} lead to fractional powers of $\log x$. Finally, signed weights such as Möbius convolutions \cite{BasakRoblesZaharescu} lie outside our setting, since Landau's theorem and the positivity of $\Psi''(\rho)$ are no longer available. It would be interesting to extend the method to these cases. 

\subsection*{Note on AI tools}
We used AI tools for assisting with computations (especially for \Cref{sec:Examples}), references, and proofreading.


\end{document}